\documentclass[11pt]{amsart}
\usepackage{latexsym, amssymb, mathrsfs, graphics, fullpage}
\newtheorem{theorem}{Theorem}
\newtheorem{cor}[theorem]{Corollary}
\newtheorem{lemma}[theorem]{Lemma}
\newtheorem{prop}[theorem]{Proposition}

\def\og{\leavevmode\raise.3ex\hbox{$\scriptscriptstyle\langle\!\langle\,$}}
\def\fg{\leavevmode\raise.3ex\hbox{$\scriptscriptstyle\,\rangle\!\rangle\ $}}

\begin{document}

\title[Group-theoretic compactification of Bruhat-Tits buildings]
{Group-theoretic compactification of Bruhat-Tits buildings}
\author{Yves Guivarc'h and Bertrand R\'emy}
\maketitle

\vspace{1cm}Ê

{\footnotesize
{\bf Abstract:}~
Let $G_F$ denote the rational points of a semisimple group $G$ over a non-archimedean local field $F$, with Bruhat-Tits building $X$. 
This paper contains five main results. 
We prove a convergence theorem for sequences of parahoric subgroups of $G_F$ in the Chabauty topology, which enables to compactify the vertices of $X$. 
We obtain a structure theorem showing that the Bruhat-Tits buildings of the Levi factors all lie in the boundary of the compactification. 
Then we obtain an identification theorem with the polyhedral compactification (previously defined in analogy with the case of symmetric spaces). 
We finally prove two parametrization theorems extending the Bruhat-Tits dictionary between maximal compact subgroups and vertices of $X$: one is about Zariski connected amenable subgroups, and the other is about subgroups with distal adjoint action. 


}

\vspace{5mmplus1mmminus1mm}Ê

\section*{Introduction}
\label{s - intro}

There exist lots of deep motivations to construct compactifications of symmetric spaces and Euclidean buildings. 
One of them is to determine the cohomological properties of arithmetic groups. 
When the ambient algebraic group is defined over a global field of characteristic zero, this was done in \cite{BorelSerre1} and \cite{BorelSerre2}. 
In positive characteristic, there are still important open questions \cite[\S VII]{Brown}, \cite{Behr}. 
Another related motivation is the compactification of locally symmetric manifolds \cite{Satake2}, in particular those carrying a natural complex structure. 
Some compactifications \cite{BailyBorel} are useful tools in number theory, and also provide nice examples of complex projective varieties of general type \cite{AMRT} or moduli spaces \cite{AllcockCarlsonToledo}. 
We refer to \cite{BorelJi} for a recent review of compactifications of symmetric spaces and of their quotients by lattices of the isometry group. 

\vspace{2mmplus1mmminus1mm} 

In this paper we are interested in compactifying Euclidean buildings by group-theoretic techniques. 
In return, we obtain group-theoretic results, e.g. geometric parametrizations of classes of remarkable closed subgroups in non-archimedean semisimple Lie groups. 
The analogy with the case of symmetric spaces is of course highly relevant. 
In fact, our general project is to generalize to Bruhat-Tits buildings all the compactification procedures described in \cite{GJT} or \cite{Guivarc'h} in the real case. 
Furstenberg (i.e. measure-theoretic) and Martin compactifications are included in the project, but will not appear in the present paper. 
Here, we are interested in the simplest approach: in the real case, it consists in seeing each point of a symmetric space as a maximal compact subgroup of the isometry group (via the Bruhat-Tits fixed-point lemma) and reminding that the space $\mathscr{S}(G)$ of all closed subgroups of a given locally compact group $G$ has a natural compact topology. 
The latter topology is the Chabauty topology \cite{Integration78}, and the compactification under consideration is the closure of the image of the map attaching to each point its isotropy subgroup. 
Of course, one has to check that this map is a topological embedding onto its image, and this is done by proving the Chabauty convergence of 
a suitable class of sequences in the symmetric space. 
This compactification is equivariantly homeomorphic to the maximal Satake compactification \cite{Satake1}, which itself was identified by C.C. Moore with the maximal Furstenberg compactification \cite{MooreCompactification}. 

\vspace{2mmplus1mmminus1mm} 

In the non-archimedean case, an additional subtlety is that the Bruhat-Tits building, i.e. the analogue of the symmetric space in this situation, is bigger than the set of maximal compact subgroups (which corresponds to the vertices of the building). 
The statement below deals with sequences of maximal compact subgroups only, but our most general result takes into account
sequences of parahoric subgroups (see Theorem \ref{th - CV groups} for a precise version). 

\vspace{2mmplus1mmminus1mm} 

{\bf Convergence theorem.}\it~
Let $G$ be a semisimple group over a local field $F$, and let $X$ be its Bruhat-Tits building. 
Let $\{v_n\}_{n\geq 1}$ be a sequence of vertices in some closed Weyl chamber $\overline{\mathscr{Q}}^X$. 
By passing to stabilizers in $G_F$ we obtain a sequence of maximal compact subgroups $\{K_{v_n}\}_{n\geq1}$. 
We make the following assumption: 

\begin{enumerate}
\item [] for each codimension one sector panel $\Pi$ of $\overline{\mathscr{Q}}^X$, the distance $d_X(v_n,\Pi)$ has a, possibly infinite, limit as $n\to\infty$. 
\end{enumerate}

Then $\{K_{v_n}\}_{n\geq1}$ is a convergent sequence in the space of closed subgroups $\mathscr{S}(G_F)$, endowed with the compact Chabauty topology. 
The limit group $D$ is Zariski dense in some parabolic $F$-subgroup $Q$ fixing a face of the chamber at infinity $\partial_\infty\mathscr{Q}$. 
Moreover $D$ can be written as a semi-direct product $K\ltimes \mathscr{R}_u(Q)_F$, where $K$ is an explicit maximal compact subgroup of some reductive Levi factor of $Q$, and $\mathscr{R}_u(Q)_F$ is the unipotent radical of $Q_F$. 
\rm\vspace{2mmplus1mmminus1mm} 

As already mentioned, this convergence is the key fact to define a compact space $\overline V_X^{\rm gp}$ with a natural $G_F$-action. We call it the group-theoretic compactification of $X$. 
The next step then is to understand the geometry of $\overline V_X^{\rm gp}$ by means of the structure of $G_F$. 
For instance, in the Borel-Serre compactification, the boundary reflects the combinatrorics of the parabolic subgroups defined over the ground field of the isometry group; a single spherical building is involved. 
In our case, the group-theoretic compactification of the Euclidean building of each Levi factor of $G_F$ appears, as in Satake's compactifications of symmetric spaces \cite{Satake1}.  
The result below sums up Theorem \ref{th - BT in boundary}. 

\vspace{2mmplus1mmminus1mm} 

{\bf Structure theorem.}\it~
For any proper parabolic $F$-subgroup $Q$, the group-theoretic compactification of the Bruhat-Tits building of the semisimple $F$-group $Q/\mathscr{R}(Q)$ naturally sits in the boundary of $\overline V_X^{\rm gp}$. 
Let $P$ be a minimal parabolic $F$-subgroup of $G$. 
We set $D_\varnothing=K\ltimes\mathscr{R}_u(P)_F$, where $K$ is the maximal compact subgroup of some reductive Levi factor of $P$. 
Then the conjugacy class of $D_\varnothing$ is $G_F$-equivariantly homeomorphic to the maximal Furstenberg boundary $\mathscr{F}$ and is the only closed $G_F$-orbit in  $\overline V_X^{\rm gp}$. 
In fact, for any closed subgroup $D\!\in\!\overline{V}^{\rm gp}_X$ there is a sequence $\{g_n\}_{n\geq 1}$ in $G_F$ such that 
$\displaystyle \lim_{n\to\infty}g_nDg_n{}^{-1}$ exists and lies in $\mathscr{F}$. 
\rm\vspace{2mmplus1mmminus1mm} 

There are then two ways to exploit this result. 
The first one is to use it to compare $\overline V_X^{\rm gp}$ with previous compactifications of Bruhat-Tits buildings. 
One compactification was defined in \cite{Landvogt} by compactifying apartments first (the convergence there is the same as in flats of maximal Satake's compactifications), and then by extending the Bruhat-Tits gluing procedure which defines $X$ out of $G_F$ and the model of an apartment. 
We call the so-obtained compactification the polyhedral compactification of $X$ and we denote by $\overline V_X^{\rm pol}$ the closure of the vertices in the latter space. 
We fill a gap in [loc. cit.] pointed to us by A. Werner and then compare $\overline V_X^{\rm gp}$ and $\overline V_X^{\rm pol}$, see Theorem \ref{th - identification polyhedral} for a more precise version. 

\vspace{2mmplus1mmminus1mm} 

{\bf Identification theorem.}\it~
Let $G$ be a semisimple simply connected group defined over a non-archimedean local field $F$. 
Let $X$ be the corresponding Bruhat-Tits building. 
Then there exists a natural $G_F$-equivariant homeomorphism $\overline{V}^{\rm pol}_X\simeq\overline{V}^{\rm gp}_X$. 
\rm\vspace{2mmplus1mmminus1mm} 

The second use of the compactification is to parametrize remarkable classes of closed subgroups. 
By taking stabilizers, we can extend to the non-archimedean case a theorem of C.C. Moore's \cite{MooreAmenable} which answers a question of H. Furstenberg's and unifies the classification of maximal compact and minimal parabolic subgroups in the same geometric framework. 
The proofs of the result below (made more precise in Theorem \ref{th - parametrization amenable}) and of the next one use subtle results due to Ph. Gille on unipotent elements in algebraic groups, in order to cover the case of a local ground field of characteristic $p>0$. 

\vspace{2mmplus1mmminus1mm} 

{\bf Parametrization theorem (amenable case).}\it~
Any closed, amenable, Zariski connected subgroup of $G_F$ fixes a facet in $\overline{X}^{\rm pol}$. 
The closed amenable Zariski connected subgroups of $G_F$, maximal for these properties, are the vertex fixators for
the $G_F$-action on the polyhedral compactification $\overline{X}^{\rm pol}$. 
\rm\vspace{2mmplus1mmminus1mm} 

To state our last main result, we need to go back to the very definition of the group-theoretic compactification. 
Since it is the closure of the maximal compact subgroups in the compact space $\mathscr{S}(G_F)$ of closed subgroups of $G_F$, it is natural to ask for an intrinsic characterization of the groups appearing by taking the closure. 
A satisfactory answer is given by the notion of distality, which formalizes the fact that a group action on a metric space has no contraction. 
Here the linear action under consideration is given by the adjoint representation. 
We refer to Theorem \ref{th - parametrization distal} and its preliminaries for a more precise version. 

\vspace{2mmplus1mmminus1mm} 

{\bf Parametrization theorem (distal case).}\it~
Any subgroup of $G_F$ with distal adjoint action is contained in a point of $\overline V_X^{\rm gp}$. 
The maximal distal subgroups of $G_F$ are the groups of $\overline V_X^{\rm gp}$, i.e. they are the maximal compact subgroups and the limits of sequences of maximal compact subgroups. 
In particular, they are all closed and Zariski connected. 
\rm\vspace{2mmplus1mmminus1mm} 

Let us finish the presentation of our results by mentioning that our proofs may simplify some arguments in the real case. 
We also plan to compare our approach to some concrete compactifications defined by A. Werner in the ${\rm SL}_n$ case \cite{WernerSublattices}, \cite{WernerSeminorms}. 
In the present paper, the latter case is presented in the last section as an illustration of the general semisimple case. 
Unfortunately, it would have been too long to develop it completely, but we think that checking the details is useful to have a good intuition of the geometry of Euclidean buildings. 
At last, since algebraic group theory uses a lot of notation, we found useful to collect some of it below. 

\vspace{2mmplus1mmminus1mm} 

{\bf Notation.}~
In all this paper, we assume we are considering the following objects: 
\begin{enumerate}
\item[--] a locally compact non-archimedean local field $F$, with valuation ring $\mathscr{O}_F$, uniformizer $\varpi_F$ 
and residue field $\kappa_F=\mathscr{O}_F/\varpi_F$. 
The absolute value is denoted by $\mid \cdot \mid_F$ and the discrete valuation by $v_F : F \to {\bf Z}$. 
We set: $q_F=\,\mid\!\kappa_F\!\mid$; 
\item[--] a simply connected semisimple algebraic $F$-group $G$; 
\item[--] the Bruhat-Tits building $X$ of $G_{/F}$, whose set of vertices is denoted by $V_X$. 
\end{enumerate}
We let $\overline F$ be an algebraic closure of the field $F$. 
There is a unique valuation $v_{\overline F}:{\overline F}\to{\bf Q}$ (resp. absolute value $\mid \cdot \mid_{\overline F}$) extending 
$v_F$ (resp. $\mid \cdot \mid_F$), and we denote by $\mathscr{O}_{\overline F}$ the valuation ring of $\overline F$. 

\vspace{2mmplus1mmminus1mm} 

In general, given an algebraic group $H$ over $F$, we denote by ${\rm Lie}(H)$ or $\mathfrak{h}$ its Lie algebra, by $H_F$ its $F$-rational points and by  ${\rm Lie}(H)_F$ or $\mathfrak{h}_F$ the $F$-rational points of the Lie algebra of $H$. 
We denote by $\mathscr{R}(H)$ (resp. by $\mathscr{R}_u(H)$) the radical (resp. the unipotent radical) of $H$.  

\vspace{2mmplus1mmminus1mm}

{\bf Structure of the paper.}~
Section \ref{s - BT buildings} fixes notation and recalls basic facts on algebraic groups and Bruhat-Tits theory; it also introduces the class of fundamental sequences in Euclidean buildings. 
Section \ref{s - groups} is mainly devoted to studying convergence of fundamental sequences of parahoric subgroups for the Chabauty topology on closed subgroups in the semisimple group $G_F$; this is the main step to define the group-theoretic compactification $\overline{V}_X^{\rm gp}$ of the Bruhat-Tits building $X$. 
Section \ref{s - description} describes $\overline{V}_X^{\rm gp}$ and in particular shows that, as a $G_F$-space, 
$\overline{V}_X^{\rm gp}$ contains a single closed $G_F$-orbit; moreover the compactification of the Bruhat-Tits building of any Levi factor lies in the boundary of $\overline{V}_X^{\rm gp}$. 
We also prove the identification with the polyhedral compactification $\overline{X}^{\rm pol}$.  
Section \ref{s - trees} deals with compactifications of trees in a slightly more general context than rank one algebraic groups over local fields; it can be seen both as an illustration of the previous sections and the first step of induction arguments in the next section. 
Section \ref{s - parametrization} contains the proofs of the two parametrization theorems in terms of the geometry of the compactification $\overline{X}^{\rm pol}$; the two parametrized classes of subgroups are that of maximal Zariski connected amenable and of maximal distal subgroups. 
Section \ref{s - GL} provides examples of arbitrary positive $F$-rank since it deals with special linear groups; we recall Goldman-Iwahori's concrete definition of the Bruhat-Tits building of ${\rm SL}_n(F)$ and we try to illustrate as many previous notions as possible. 

\vspace{2mmplus1mmminus1mm}

{\bf Acknowledgements.}~We thank H. Abels, H. Behr, M. Brion, Ph. Gille, N. Monod, F. Paulin, G. Rousseau and A. Werner for useful discussions during the preparation of this paper. 

\vspace{1cm}Ê

\section{Bruhat-Tits buildings. Levi factors. Unbounded sequences}
\label{s - BT buildings}

We introduce some algebraic subgroups and notation, and we recall the geometric meaning of the valuated root datum axioms for the rational points $G_F$. 
We also recall a more technical point: the Bruhat-Tits building of a Levi factor naturally sits in the Bruhat-Tits building of the ambient group $G_{/F}$. 
We finally use the Cartan decomposition with respect to a suitable maximal compact subgroup, in order to distinguish a class of sequences in buildings which will become convergent in the group-theoretic compactification of the next section. 

\subsection{Bruhat-Tits buildings}
\label{ss - Bruhat-Tits}
We choose once and for all a maximal $F$-split torus $T$ in $G$, to which is associated an apartment $A$ in the Bruhat-Tits building $X$ \cite[2.8.11]{BrT1}. 
We denote by $\Phi=\Phi(T,G)$ the corresponding root system \cite[8.17]{Borel}. 
It is a (possibly non-reduced) root system \cite[5.8]{BorelTits65} in the sense of \cite[VI.1]{Lie456}. 
In order to avoid confusions and to emphasize the analogy with symmetric spaces, a {\it sector}~of $X$ \cite[VI.7]{Brown} (in French \og quartier\fg \cite[7.1.4]{BrT1}) is called a {\it Weyl chamber}~in this article.
We also use the terminology {\it alcove} \cite[1.3.8]{BrT1}, so that the word \og chamber\fg alone is meaningless in the present paper.  

\subsubsection{}
\label{sss - apartment and standard subgroups}
Let us pick in $A$ a Weyl chamber $\mathscr{Q}$ with tip a special vertex which we call $o$. 
Let us denote by $c$ the alcove contained in $\mathscr{Q}$ whose closure contains $o$. 
We refer to $A$ (resp. $\mathscr{Q}$, $o$, $c$) as the {\it standard}~apartment (resp. Weyl chamber, vertex, alcove) of $X$. 
The fixator $K_o={\rm Fix}_{G_F}(o)$ is called the {\it standard maximal compact subgroup}~and its subgroup $\mathscr{B}={\rm Fix}_{G_F}(c)$ is called the {\it standard Iwahori subgroup}~of $G_F$. 
The choice of $\mathscr{Q}$ corresponds to the choice of a subset of positive roots $\Phi^+$, or equivalently to the choice of a system of simple roots $\{a_s\}_{s\in S}$ which we identify with its indexing set $S$. 
We set: $\Phi^-=-\Phi^+$. 
By definition of $\Phi$, we have a decomposition of the Lie algebra $\mathfrak{g}$ as a $T$-module via the adjoint representation: 
$\mathfrak{g}=\mathfrak{g}_0\oplus\bigoplus_{a\in\Phi}\mathfrak{g}_a$, where $\mathfrak{g}_0$ is the fixed-point set of $T$. 
The group $U_a$ with Lie algebra $\mathfrak{g}_a=\{X\!\in\!\mathfrak{g}:{\rm Ad}(t).X=a(t)X$ for each $t\!\in\!T\}$ is the {\it root group}~associated to $a$. 
We denote by $P$ the minimal parabolic $F$-subgroup determined by $\Phi^+$, i.e. such that 
$\mathfrak{p}=\mathfrak{g}_0\oplus\bigoplus_{a\in\Phi^+}\mathfrak{g}_a$. 
For a subset $I$ of $S$, we denote by $\Phi_I$ the subset of roots which are linear combinations of simple roots indexed by $I$. 
We also set: $\Phi_I^\pm=\Phi^\pm\cap\Phi_I$, and $\Phi^{I,\pm}=\Phi^\pm-\Phi_I$. 
We also introduce the following $F$-subgroups \cite[\S 4]{BorelTits65}: 

\begin{enumerate}
\item[] the {\it standard parabolic subgroup}~$P_I$ of type $I$, defined by $\mathfrak{p}_I=\mathfrak{g}_0\oplus\bigoplus_{a\in\Phi_I}\mathfrak{g}_a\oplus\bigoplus_{a\in\Phi^{I,+}}\mathfrak{g}_a$; 
\item[] the {\it standard reductive Levi factor}~$M_I$ of $P_I$, defined by 
$\mathfrak{m}_I=\mathfrak{g}_0\oplus\bigoplus_{a\in\Phi_I}\mathfrak{g}_a$;  
\item[] the {\it standard semisimple Levi factor}~$G_I=[M_I,M_I]$ of $P_I$; 
\item[] the {\it standard unipotent radical}~$U^I=\mathscr{R}_u(P_I)$ of type $I$, also defined by 
$\mathfrak{u}^I=\bigoplus_{a\in\Phi^{I,+}}\mathfrak{g}_a$. 
\end{enumerate}

The parabolic $F$-subgroup {\it opposite $P_I$}~with respect to $T$ \cite[14.20]{Borel} is the $F$-subgroup with Lie algebra 
$\mathfrak{m}_I\oplus\bigoplus_{a\in\Phi^{I,-}}\mathfrak{g}_a$. 
When $I=\varnothing$, we simply omit the index ${}_\varnothing$, e.g. we denote $M=M_\varnothing$. 
Note that in the classification theory of semisimple $F$-groups, the ($G_F$-conjugacy class of the) semisimple Levi factor $[M,M]$ is usually called the {\it anisotropic kernel}~of $G_{/F}$ \cite[16.2.1]{Springer}. 
The $F$-points $[M,M]_F$ form a compact group \cite[5.2.3]{Rousseau}, and we denote by $Z_G(T)_{\rm cpt}$ the commutative product 
$[M,M]_F \cdot T_{\rm cpt}$, where $T_{\rm cpt}$ is the unique maximal compact subgroup of $T_F$. 
The group $Z_G(T)_{\rm cpt}$ is the pointwise fixator of the apartment $A$ and the maximal compact subgroup of $M_F=Z_{G_F}(T_F)$. 

\subsubsection{}
\label{sss - valuated root data and half spaces} 
A substantial part of Borel-Tits theory (i.e. the theory of reductive groups over an arbitrary ground field \cite{BorelTits65}) can be summed up in combinatorial terms \cite[21.15]{Borel}. 
The most refined version of this approach is provided by the notion of a {\it generating root datum}~\cite[6.1]{BrT1}. 
This is relevant to the case of an arbitrary ground field. 
The combinatorics becomes richer when the ground field is a local field $F$ as in our paper. 
One of the main results of Bruhat-Tits theory is the existence of a {\it valuation}~on the generating root datum of $G_F$ associated to the choice of the maximal $F$-split torus $T$ \cite[5.1.20]{BrT2}. 
Each root group $U_a$ is unipotent, abelian or metabelian, and roughly speaking the latter notion corresponds to the existence of a filtration on each group $(U_a)_F$; any such filtration comes from the filtration of the additive group $(F,+)$ given by the preimages of $v_F$. 
Further compatibilities (e.g. with respect to the action of the normalizer of $T_F$, to taking some commutators) are required, but we are only interested in the geometric interpretation of this valuation \cite[1.4]{TitsCorvallis}. 

\vspace{2mmplus1mmminus1mm} 

Let $b\!\in\!\Phi$ be a root and let $\{U_{b,m}\}_{m\in{\bf Z}}$ be the decreasing filtration of $(U_b)_F$ given by Bruhat-Tits theory \cite[6.2]{BrT1}. 
To the pair of opposite roots $\{\pm b\}$ is attached a parallelism class $\partial b$ of affine hyperplanes in $A$ and a family of affine hyperplanes $\{H_{\partial b,m}\}_{m\in{\bf Z}}$ in this class, which are called the {\it walls}~directed by $\partial b$. 
The family $\{H_{\partial b,m}\}_{m\in{\bf Z}}$ provides a useful exhaustion of $A$ by the fixed-point sets of the groups $U_{b,m}$ \cite[2.1]{TitsCorvallis}. 
More precisely, if $b\!\in\!\Phi^+$ we denote by $D_{b,m}$ the half-space of $A$ bounded by $H_{\partial b,m}$ which contains a translate of the Weyl chamber $\mathscr{Q}$; otherwise, we choose the other half-space to be $D_{b,m}$. 
The family $\{H_{\partial b,m}\}_{m\in{\bf Z}}$ is characterized by the fact that the fixed-point set of $U_{b,m}$ in $A$ is equal to the half space 
$D_{b,m}$. 
In other words, we have an increasing exhaustion $A=\bigcup_{m\in{\bf Z}} D_{b,m}$. 
Geometrically, the bigger $m\!\in\!{\bf Z}$ is, the smaller $U_{b,m}$ is, and the bigger $D_{b,m}=A^{U_{b;,m}}$ is. 
We suggest the reader to have a look at the third paragraph of \ref{sss - parametrization for GL} which deals with the example of 
${\rm SL}_n(F)$. 

\subsubsection{}
\label{sss - parabolic and asymptotic} 
As a combinatorial Euclidean building, $X$ can be endowed with a distance $d$ of non-positive curvature, unique up to homothety on each irreducible factor \cite[2.5]{BrT1}. 
More precisely, the distance $d$ makes $X$ a CAT$(0)$-space \cite[II]{BriHae}; we fix once and for all such a metric $d$. 
The  boundary at infinity $\partial_\infty X$, i.e. the space of geodesic rays modulo the relation of being at finite Hausdorff distance from one another \cite[II.8]{BriHae}, is a geometric realization of the spherical building of parabolic subgroups in $G_F$ \cite[VI.9E]{Brown}. 
Therefore, we can also define the $F$-points of the parabolic subgroups of $G_{/F}$ to be the stabilizers of the facets at infinity in $\partial_\infty X$. 
For instance, the standard sector $\mathscr{Q}$ defines a chamber at infinity $\partial_\infty\mathscr{Q}$, and we have: 
$P_F={\rm Stab}_{G_F}(\partial_\infty\mathscr{Q})={\rm Fix}_{G_F}(\partial_\infty\mathscr{Q})$. 
The inclusion $\partial_\infty\mathscr{Q}\subset \partial_\infty A$ corresponds to the Levi decomposition $P=M\ltimes U$: the group $M$ is characterized by the fact that $M_F$ is the fixator of the union of the facet $\partial_\infty\mathscr{Q}$ and its opposite in $\partial_\infty A$. 
There is a similar geometric interpretation for each standard Levi decomposition $P_I=M_I\ltimes U^I$. 

\vspace{2mmplus1mmminus1mm} 

Recall that a {\it Furstenberg boundary}~for a topological group $G$ is a compact metrizable $G$-space $Y$ whose continuous $G$-action is {\it minimal}~(i.e. any orbit in $Y$ is dense) and {\it strongly proximal}~(i.e. there is a Dirac mass in the closure of any $G$-orbit in the space $\mathscr{M}^1(Y)$ of probability measures on $Y$) \cite[VI.1.5]{Margulis}. 
It is a classical fact that the family of Furstenberg boundaries of a given semisimple Lie group coincides with the family of its flag varieties; in the non-archimedean case, this is checked for instance in \cite[Proposition 5.1]{BurMozCAT}. 
Moreover, if $Q$ denotes a parabolic $F$-subgroup of $G$, the $F$-rational points of $G/Q$ form a homogeneous space under $G_F$; in fact, we have: $(G/Q)_F=G_F/Q_F$ \cite[4.13]{BorelTits65}. 
In this paper the quotient spaces $G_F/Q_F$, for $Q$ a parabolic $F$-subgroup of $G_{/F}$, are indifferently called 
{\it Furstenberg boundaries}~or {\it flag varieties}. 

\vspace{2mmplus1mmminus1mm} 

{\bf Definition.}~\it 
We denote by $\mathscr{F}$ the {\rm (maximal) Furstenberg boundary}~$G_F/P_F$ of $G_F$. 
For each subset $I$ of simple roots of $S$, we denote by $\mathscr{F}^I$ the Furstenberg boundary $G_F/(P_I)_F$. 
\rm\vspace{2mmplus1mmminus1mm} 

There is an obvious $G_F$-equivariant map $\pi_I:\mathscr{F}\to\mathscr{F}^I$ between Furstenberg boundaries. 
We denote by $\omega$ (resp. $\omega_I$) the class of the identity in $\mathscr{F}$ (resp. in $\mathscr{F}^I$). 
The preimage of $\omega_I$ by $\pi_I$ is $(P_I)_F.\omega$. 
It is a copy in $\mathscr{F}$ of the (maximal) Furstenberg boundary of the Levi factor $(M_I)_F$; we denote it by $\mathscr{F}_I$. 
We denote by $U^{I,-}$ the unipotent radical of the parabolic subgroup opposite $P_I$ with respect to $T$; we simply write $U^-$ for $U^{\varnothing,-}$. 
Note that there is a unique $G_F$-invariant class of measures on $\mathscr{F}$ \cite[\S 1]{Mackey}, and that the $(U^-)_F$-orbit of $\omega$ is conegligible for this class. 
In the algebraic terminology, $(U^-)_F.\omega$ is called the {\it big cell}~of $\mathscr{F}$. 
Let us finally recall that there is a natural way to glue $G_F$-equivariantly $\partial_\infty X$ to $X$ \cite[II.8]{BriHae}. 
The so-obtained space is called the {\it geometric compactification}~of $X$; we denote it by $\overline{X}^{\rm geom}$. 
The partition of the boundary of $\overline{X}^{\rm geom}$ under its natural $K_o$-action is connected to flag varieties by the fact that each $K_o$-orbit is isomorphic, as a $K_o$-space, to some flag variety of $G_F$. 

\subsection{Levi factors}
\label{ss - Levi}
We recall that the Bruhat-Tits building $X$ of $G_{/F}$ contains inessential realizations of the buildings of the Levi factors in $G$ \cite[7.6]{BrT1}. 
We also introduce some remarkable closed subgroups of $G_F$. 
They will turn out to form the boundary of our compactification of $X$, or to be the stabilizers of the points in this boundary. 

\subsubsection{}
\label{sss - Levi BT buildings}
We simply recall here (with our notation) the facts of \cite[7.6]{BrT1} we will be using later. 
Recall that the apartment $A$ is the vector space $X_*(T)_F\otimes_{\bf Z}{\bf R}$ endowed with a suitable simplicial structure and a natural action of $N_{G_F}(T_F)$ ($X_*(T)_F\otimes_{\bf Z}$ is the group of $F$-cocharacters of $T$) \cite[1.2]{TitsCorvallis}. 
Let $I$ be a proper subset of simple roots in $S$. 
Let $L_I$ denote the affine subspace of $A$ obtained as the intersection of the walls passing through the vertex $o$ and directed by the simple roots in $I$. 
The semisimple Levi factor $G_I$ is simply connected \cite[8.4.6, exercise 6]{Springer} and we denote by $A_I$ the standard apartment of its Bruhat-Tits building $X_I$. 
We are interested in the subset $(G_I)_F.A$ of the $(G_I)_F$-transforms of the points in $A$, which we want to compare to $X_I$. 
Intuitively, the idea is that the direction $L_I$ is not relevant to the combinatorics of the semisimple Levi factor $G_I$: it corresponds to the cocharacters of $T$ which centralize $G_I$. 
But after shrinking $(G_I)_F.A$ along $L_I$, we obtain a realization of $X_I$. 
This is formalized by \cite[Proposition 7.6.4]{BrT1} which provides a unique extension $\tilde p_I:(G_I)_F.A\to X_I$ of the natural affine map $p_I:A\to A_I$ between apartments, with the following properties: 

\begin{enumerate}
\item [(i)]Êthe map $\tilde p_I$ is $(G_I)_F$-equivariant; 
\item [(ii)] the preimage of $A_I$ by $\tilde p_I$ is $A$, and in fact the preimage of any apartment, wall, half-apartment in $X_I$ is an apartment, a wall or a half-apartment in $X$, respectively; 
\item [(iii)] there is an $L_I$-action on $(G_I)_F.A$ extending that on $A$ with the following compatibility with the $(G_I)_F$-action: 
$g.(x+v)=g.x+v$ for any $g\!\in\!(G_I)_F$, $x\!\in\!(G_I)_F.A$ and $v\!\in\!L_I$; 
\item [(iv)] the factor map $\displaystyle {(G_I)_F.A\over L_I}\to X_I$ is a $(G_I)_F$-equivariant bijection. 
\end{enumerate}

The choice of positive roots in $\Phi$ corresponding to the Weyl chamber $\mathscr{Q}$ induces a choice of positive roots in the subroot system of $G_I$, and the corresponding  Weyl chamber of $A_I$ is $p_I(\mathscr{Q})$. 

\subsubsection{}
\label{sss - simplicial cone}
We denote by $\overline{\mathscr{Q}}^X$ the (non-compact) closure of the Weyl chamber $\mathscr{Q}$ in the building $X$. 
It is a simplicial cone. 
Any of its codimension one faces is equal to the intersection of $\overline{\mathscr{Q}}^X$ with the wall directed by some simple root $s\!\in\!S$ and passing through $o$. 
We denote by $\Pi^s$ the latter cone, and we call it a {\it sector panel}~of $\overline{\mathscr{Q}}^X$. 
For any non-empty subset $I$ of simple roots, we denote: $\mathscr{Q}^I=\bigcap_{s\in I}\Pi^s$; for instance $\mathscr{Q}^S=\{o\}$. 
We set: $\mathscr{Q}^\varnothing=\mathscr{Q}$, so that $\overline{\mathscr{Q}}^X=\bigsqcup_{I\subset S}\mathscr{Q}^I$. 
Note also that $L_I$ above is the affine subspace generated by $\mathscr{Q}^I$. 

\vspace{2mmplus1mmminus1mm} 

The points in $\overline{\mathscr{Q}}^X$ are parametrized by the distances to the sector panels $\Pi^s$ when $s$ ranges over $S$. 
Given a family $\underline d=\{ d_s \}_{s\in S}$ of non-negative real numbers, we denote by $x_{\underline d}$ the corresponding point in 
$\overline{\mathscr{Q}}^X$. 
Similarly, in the building $X_I$ we parametrize the closed Weyl chamber $\overline{p_I(\mathscr{Q})}^{X_I}$ by the set of finite sequences $\underline d=\{ d_s \}_{s\in I}$ of non-negative real numbers (corresponding to the distances to the sector panels $p_I(\Pi^s)$, $s\!\in\! I$). 
The point defined by the parameters $\underline d$ is denoted by $x_{I,\underline d}$. 
The preimage $p_I^{-1}(x_{I,\underline d})$ is an affine subspace of $A$ parallel to $\langle\mathscr{Q}^I\rangle$, 
which we denote by $L_{I,{\underline d}}$. 
Each space $L_{I,{\underline d}}$ has dimension ${\rm rk}_F(G)\,-\mid\!I\!\mid\,={\rm dim}(X)\,-\mid\!I\!\mid$, where ${\rm rk}_F$ denotes the $F$-rank of a semisimple $F$-group. 

\vspace{2mmplus1mmminus1mm} 

{\bf Definition.}~\it
We define $K_{I,\underline d}$ to be pointwise fixator of the affine subspace $L_{I,{\underline d}}$ in the reductive Levi factor $(M_I)_F$, i.e. 
$K_{I,\underline d}={\rm Fix}_{(M_I)_F}(L_{I,\underline d})$. 
\rm\vspace{2mmplus1mmminus1mm} 

Denoting by $\sigma$ the facet of $A_I$ containing $x_{I,\underline d}$, we have: 
$K_{I,\underline d}={\rm Fix}_{(M_I)_F}\bigl(p_I^{-1}(\sigma)\bigr)$, i.e. the group $K_{I,\underline d}$ only depends on the facet of $p_I(\overline{\mathscr{Q}}^X)$ containing $x_{\underline d}$. 
The group $G_I$ is generated by the root groups $U_a$, where $a$ is a root such that $\partial a$ contains $\langle\mathscr{Q}^I\rangle$. 
The group $K_{I,\underline d}$ is the parahoric subgroup of $(M_I)_F$ generated by $Z_G(T)_{\rm cpt}$ and the groups $U_{a,m}$ with $a$ as before and $m\!\in\!{\bf Z}$ such that the closed half-apartment $\overline{D_{a,m}}$ contains the affine subspace $L_{I,{\underline d}}$ \cite[6.4]{BrT1}. 
Note that $p_I(\langle\mathscr{Q}^I\rangle)$ is a special vertex in $A_I$, which we choose as origin in $A_I$.
We simply write $K_I$ when $L_{I,\underline d}=\langle\mathscr{Q}^I\rangle=L_I$. 
The group $K_I$ (resp. $K_I\cap G_I$) is a special maximal compact subgroup of the standard reductive Levi factor $(M_I)_F$ (resp. of the semisimple Levi factor $(G_I)_F$). 

\subsubsection{}
\label{sss - definition limit groups}
We denote by $T^I$ the subtorus of $T$ such that $M_I=Z_G(T^I)$ and by $T_I$ the subtorus defined as the identity component 
$(G_I \cap T)^\circ$ for the Zariski topology. 
The latter group is a maximal $F$-split torus of $G_I$. 
The $F$-points of the former group have a geometric interpretation: 
$(T^I)_F={\rm Stab}_{T_F}(\langle\mathscr{Q}^I\rangle)={\rm Fix}_{T_F}(\partial_\infty\mathscr{Q}^I)$. 
In fact, $(T^I)_F$ acts on each any affine subspace $L_{I,\underline d}$ in $A$ as a discrete cocompact group of translations. 
Intersecting the decomposition $M_I=G_I\cdot T^I$ along $K_{I,\underline d}$ provides a decomposition: 
$K_{I,\underline d}=(G_I \cap K_{I,\underline d})\cdot{\rm Fix}_{(T^I)_F}(L_{I,\underline d})$. 
The factor ${\rm Fix}_{(T^I)_F}(L_{I,\underline d})$ is the unique maximal compact subgroup of $(T^I)_F$, and $G_I \cap K_{I,\underline d}$ is a parahoric subgroup of $(G_I)_F$. 
Intersecting $M_I=G_I\cdot T^I$ along ${\rm Stab}_{(M_I)_F}(L_{I,\underline d})$ provides the decomposition: 
${\rm Stab}_{(M_I)_F}(L_{I,\underline d})=K_{I,\underline d}\cdot (T^I)_F$, with a well-understood action of each factor on $L_{I,\underline d}$. 
The groups $D_{I,\underline d}$ and $R_{I,\underline d}$ below will play an important role in the definition and the description of the group-theoretic compactification of $X$. 

\vspace{2mmplus1mmminus1mm} 

{\bf Definition.}~\it
Let $I$ be a proper subset of simple roots and let $\underline d$ be a sequence of non-negative real numbers indexed by $I$. 
\begin{enumerate}
\item [(i)] We define $D_{I,\underline d}$ to be the semidirect product $K_{I,\underline d} \ltimes (U^I)_F$. 
\item [(ii)] We define $R_{I,\underline d}$ to be the semidirect product $(K_{I,\underline d}\cdot(T^I)_F)\ltimes (U^I)_F$
\end{enumerate}
\rm\vspace{2mmplus1mmminus1mm} 

It follows from the above definitions that we have: $R_{I,\underline d}=D_{I,\underline d}\cdot(T^I)_F$. 

\subsection{Unbounded sequences}
\label{ss - unbounded sequences}
We define some classes of sequences in the Euclidean building $X$. 
These sequences will turn out to be convergent in the later defined compactifications. 

\subsubsection{}
\label{sss - canonical and fundamental} 
Recall that a sequence $\{x_n\}_{n\geq 1}$ in a topological space {\it goes to infinity}~if it eventually leaves any
compact subset of this space. 

\vspace{2mmplus1mmminus1mm} 

{\bf Definition.}~\it 
Let $\{x_n\}_{n\geq 1}$ be a sequence of points in the Euclidean building $X$. 
Let $\overline{\mathscr{Q}}^X$ be the closure of the Weyl chamber $\mathscr{Q}$. 
Let $I$ be a subset of the corresponding set $S$ of simple roots. 

\begin{enumerate}
\item [(i)] We say that $\{x_n\}_{n\geq 1}$ is {\rm $I$-canonical}~if the following three conditions are satisfied: 
\begin{enumerate}
\item[(i-a)] for each $n\geq 1$, we have: $x_n\!\in\!\overline{\mathscr{Q}}^X$; 
\item[(i-b)] for each $s\!\in\!S\setminus I$, we have: $\displaystyle \lim_{n\to+\infty}{\rm dist}_X(x_n,\Pi^s)=+\infty$; 
\item[(i-c)] there exists a facet $\sigma$ in $A_I \cap p_I(\overline{\mathscr{Q}}^X)$ such that $p_I(x_n)\!\in\!\sigma$ for $n>\!\!>1$. 
\end{enumerate}
\item [(ii)] We say that $\{x_n\}_{n\geq 1}$ is {\rm $I$-fundamental}~if there exists a converging sequence $\{ k_n \}_{n\geq 1}$ in $K_o$ such that $\{k_n.x_n\}_{n\geq 1}$ is an $I$-canonical sequence. 
\item [(iii)] We simply say that $\{x_n\}_{n\geq 1}$ is {\rm fundamental}~if it is $I$-fundamental for some $I \subset S$. 
\end{enumerate}
\rm\vspace{2mmplus1mmminus1mm} 

Note that an $I$-fundamental sequence is bounded if, and only if, we have: $I=S$. 
Note also that the condition for being $I$-canonical depends on the choice of the Weyl chamber $\mathscr{Q}$ while that for being fundamental doesn't. 
The condition of being $I$-canonical here is slightly more general than in the real case \cite{GJT} since we don't impose that the points be in the face of $\overline{\mathscr{Q}}^X$. 
We have to adopt this definition because $T_F$ doesn't act transitively on the intersections of walls in a given parallelism class. 

\subsubsection{}
\label{sss - discuss canonical and fundamental}
Moreover in the context of symmetric spaces, the definition of $I$-fundamental sequences is not exactly the same \cite[Definition 3.35]{GJT}. 
The plain analogue of the latter definition would be obtained by replacing condition (i-c) of \ref{sss - canonical and fundamental} by the condition: 

\vspace{2mmplus1mmminus1mm} 

\begin{enumerate}
\item [(i-c$'$)] for each $s\!\in\!I$, the distance ${\rm dist}_X(x_n,\Pi^s)$ converges as $n\to\infty$. 
\end{enumerate}

\vspace{2mmplus1mmminus1mm} 

On the one hand, let us pick two points $x$ and $y$ in the alcove $c$, such that for any $s\!\in\!S$ we have: 
${\rm dist}_X(x,\Pi^s)\neq{\rm dist}_X(y,\Pi^s)$. 
Then any sequence $\{ x_n \}_{n\geq1}$ taking infinitely many times each value $x$ and $y$ is $S$-fundamental in our sense, while it is not for the above plain translation from the case of symmetric spaces. 
On the other hand, an injective sequence of points $\{x'_n\}_{n\geq1}$ in the alcove $c$ converging to the tip $o$ is $S$-fundamental for the modified definition, but isn't for the definition we will use (\ref{sss - canonical and fundamental}). 
The sequences $\{x_n\}_{n\geq1}$ and $\{x'_n\}_{n\geq1}$ show that the two definitions of being $I$-fundamental are different for $I=S$. 
To see the same phenomenon for $I\subsetneq S$, it suffices to pick a non-trivial element $t\!\in\!T^I$ and to replace 
$\{ x_n \}_{n\geq1}$ (resp. $\{ x_n' \}_{n\geq1}$) by $\{ t^n.x_n \}_{n\geq1}$ (resp. $\{ t^n.x_n' \}_{n\geq1}$). 

\vspace{2mmplus1mmminus1mm} 

The reason why we use the definition of  \ref{sss - canonical and fundamental} is that we want the map $x\mapsto {\rm Stab}_{G_F}(x)=K_x$ to be continuous for a certain topology on closed subgroups of $G_F$ (\ref{ss - Chabauty}). 
In the sequence $\{ x_n' \}_{n\geq 1}$, all the elements belong to the alcove $c$, so the associated sequence of parahoric subgroups $\{K_{x_n'}\}_{n\geq1}$ is constant equal to the Iwahori subgroup $\mathscr{B}$. 
But the parahoric subgroup attached to $\displaystyle o=\lim_{n\to\infty} x_n'$ is the maximal compact subgroup $K_o$. 
Of course, this phenomenon occurs also for unbounded fundamental sequences; this is reformulated in terms of convergence of parahoric subgroups in \ref{sss - limit groups}. 

\subsubsection{}
\label{sss - fundamental subsequence}
It is clear that a sequence is unbounded if, and only if, it has a subsequence going to infinity. 
In our case, the existence of a Cartan decomposition of $G_F$ with respect to $\mathscr{Q}$ implies a more precise result. 
The following lemma eventually says that any sequence in the building $X$ has a convergent subsequence in suitable embeddings of $X$ (Theorem \ref{th - CV groups}). 

\begin{lemma}
\label{lemma - fundamental subsequence}
Any sequence $\{x_n\}_{n\geq 1}$ in the building $X$ has an $I$-fundamental subsequence for some $I \subset S$. 
Moreover we can choose $I$ to be proper in $S$ whenever $\{x_n\}_{n\geq 1}$ is unbounded. 
\end{lemma}

{\it Proof}.~
By Cartan decomposition with respect to $K_o$ \cite[4.4.3 (2)]{BrT1}, the closed  Weyl chamber $\overline{\mathscr{Q}}^X$ is a fundamental domain for the $K_o$-action on $X$: for each $n \geq 1$, there exist $k_n\!\in\!K_o$ and 
$q_n\!\in\!\overline{\mathscr{Q}}^X$ such that $x_n=k_n.q_n$. 
Since $K_o$ stabilizes each sphere centered at $o$, the sequence $\{x_n\}_{n\geq 1}$ is bounded if, and only if, so is $\{q_n\}_{n\geq 1}$. 
In this case, there exists a subsequence of $\{q_n\}_{n\geq 1}$ converging to some $q\!\in\!\overline{\mathscr{Q}}^X$, and since $q$ lies in finitely many closures of facets, we are done. 
From now on, $\{x_n\}_{n\geq 1}$ is assumed to be unbounded, so we may -- and shall -- assume that $\{x_n\}_{n\geq 1}$ goes to infinity. 
We set: 
$\displaystyle J_1=\{s\!\in\!S:\limsup_{n\to+\infty}\,{\rm dist}_X(q_n,\Pi^s)=+\infty\}$.
We have: $J_1\neq\varnothing$. 
We pick $s_1\!\in\!J_1$ and choose an increasing map $\psi_1:{\bf N}\to{\bf N}$ such
that 
$\displaystyle \lim_{n\to+\infty}{\rm dist}_X(q_{\psi_1(n)},\Pi^{s_1})=+\infty$. 
Then we set: 
$\displaystyle 
J_2=\{s\!\in\!S\setminus\{s_1\}:\limsup_{n\to+\infty}\,{\rm dist}_X(q_{\psi_1(n)},\Pi^s)=+\infty\}$.
If $J_2=\varnothing$, we set $I=S\setminus\{s_1\}$ and $\psi=\psi_1$; otherwise, we pick  $s_2\!\in\!S\setminus\{s_1\}$ and choose an increasing map $\psi_2:{\bf N}\to{\bf N}$ such that $\displaystyle \lim_{n\to+\infty}{\rm dist}_X(q_{\psi_1\circ\psi_2(n)},\Pi^{s_2})=+\infty$. 
After a finite number of iterations, we obtain a subset $J$ of $S$ and an increasing map $\psi:{\bf N}\to{\bf N}$ such that 
$\displaystyle \lim_{n\to+\infty}{\rm dist}_X(q_{\psi(n)},\Pi^s)=+\infty$ for $s\!\in\!J$ and 
$\displaystyle \limsup_{n\to+\infty}\,{\rm dist}_X(q_{\psi(n)},\Pi^s)<+\infty$ otherwise. 
It remains to set $I=S\setminus J$ and to pass to convergent subsequences for distances to $\Pi^s$, $s\!\in\! I$, to obtain an increasing map $\varphi:{\bf N}\to{\bf N}$ such that ${\rm dist}_X(q_{\varphi(n)},\Pi^s)$ converges for $s\!\in\!I$ and diverges otherwise.
\qed

\section{Group-theoretic compactification}
\label{s - groups}

We use the fact that for a locally compact group $G$, the set $\mathscr{S}(G)$ of closed subgroups of $G$ carries a natural compact topology with several equivalent descriptions \cite[VIII \S 5]{Integration78}. 
The starting point is to see the vertices of the building $X$ as the set of maximal compact subgroups of $G_F$, hence as a subset of $\mathscr{S}(G_F)$. 
In the context of semisimple real Lie groups, the idea is originally due to the first author. 

\subsection{Chabauty topology and geometric convergence}
\label{ss - Chabauty}
Let $G$ be a locally compact group. 
We denote by $\mathscr{C}(G)$ the set of closed subsets of $G$. 

\subsubsection{}
\label{sss - definition Chabauty}
The space $\mathscr{C}(G)$ can be endowed with a separated 
uniform structure \cite[II \S 1]{Topo2}, defined as follows \cite[VIII \S 5 6]{Integration78}. 
For any compact subset $C$ in $G$ and any neighborhood $V$ of $1_G$ in $G$, 
we define $P(C,V)$ to be the set of couples $(X,Y)$ in $\mathscr{C}(G)\times\mathscr{C}(G)$ such that: 

\vspace{2mmplus1mmminus1mm} 

\centerline{$X \cap C \subset V \cdot Y$ \hspace{1cm} and \hspace{1cm} $Y \cap C \subset V \cdot X$.}

\vspace{2mmplus1mmminus1mm} 

The sets  $P(C,V)$ form a fundamental system of {\it entourages~}of a uniform structure 
on $\mathscr{C}(G)$. 
The so-obtained topology on $\mathscr{C}(G)$ is called the {\it Chabauty topology}. 
The space $\mathscr{S}(G)$ of closed subgroups is a compact subset of 
$\mathscr{C}(G)$ for the Chabauty topology \cite[VIII \S 5 3, Th\'eor\`eme 1]{Integration78}. 

\vspace{2mmplus1mmminus1mm} 

We henceforth assume that $G$ is metrizable; then so is the Chabauty topology. 
Moreover we can define the topology of {\it geometric convergence~}on $\mathscr{C}(G)$ \cite{CEG}, in which 
a sequence $\{F_n \}_{n\geq 1}$ of closed subsets converge 
to $F\!\in\!\mathscr{C}(G)$ if, and only if, the two conditions below are satisfied: 

\vspace{2mmplus1mmminus1mm} 

\begin{enumerate}
\item[(i)]Let $\varphi:{\bf N}\to{\bf N}$ be an increasing map and let 
$\{x_{\varphi(n)}\}_{n \geq 1}$ be a sequence in $G$ such that 
$x_{\varphi(n)} \!\in\! F_{\varphi(n)}$ for any $n\geq 1$. 
If $\{x_{\varphi(n)} \}_{n\geq 1}$ converges to some $x$ in $G$, then $x\!\in\! F$. 
\item[(ii)]Any point in $F$ is the limit of a sequence $\{x_n\}_{n \geq 1}$ with 
$x_n\!\in\! F_n$ for each $n\geq 1$. 
\end{enumerate}  

In fact, both topologies coincide: 

\begin{lemma} 
\label{lemma - Chabauty/geom}
Let $\{F_n \}_{n \geq 1}$ be a sequence of closed subsets in $G$. 
Then we have geometric convergence $\displaystyle \lim_{n\to+\infty}F_n=F$ if, and only if, $F_n$
converges to $F$ in the Chabauty topology. 
\end{lemma}

{\it Proof}.~For the sake of completeness, we recall the proof of this probably well-known lemma. 

\subsubsection{}
\label{sss - geometric implies Chabauty}
{\it Geometric convergence implies Chabauty convergence}.~
Let us assume that we have geometric convergence: $\displaystyle \lim_{n\to+\infty}F_n=F$. 
Let $C$ be a compact subset in $G$ and let $\Omega$ be an open, relatively compact, 
neighborhood of $1_G$ in $G$. 

\vspace{2mmplus1mmminus1mm} 

Let us first prove that here is an index $M\geq 1$ such that $F_n\cap C\subset \Omega\cdot F$ for any 
$n\geq M$. 
If not, there exist an increasing map $\varphi:{\bf N}\to{\bf N}$ and a sequence 
$\{x_{\varphi(n)}\}_{n \geq 1}$ such that 
$x_{\varphi(n)}\!\in\!(F_{\varphi(n)}\cap C)\setminus\Omega\cdot F$ for each $n\geq 1$. 
By compactness of $C$ and up to extracting again, we may -- and shall -- assume
that $\{x_{\varphi(n)}\}_{n \geq 1}$ converges, say to $x$, in $C$. 
But by condition (i), we have: $x\!\in\!F$, so for $n>\!\!>1$ we will have: 
$x_{\varphi(n)}\!\in\! \Omega.x$, a contradiction. 

\vspace{2mmplus1mmminus1mm} 

It remains to prove that there is an index $M\geq 1$ such that 
$F\cap C\subset \Omega\cdot F_n$ for  $n\geq M$. 
Let $W$ be an open, symmetric, neighborhood of $1_G$ in $G$ 
such that $W\cdot W\subset \Omega$. 
By compactness of $F \cap C$, we can write: 
$F \cap C = \bigcup_{i=1}^l Wx^i$ for $x^1, x^2, ... \, x^l$ in $F \cap C$. 
By condition (ii), for each $i\!\in\!\{1;2;...\, l\}$ 
we can write: $x^i=\displaystyle \lim_{n\to\infty} x^i_n$, with 
$x^i_n\!\in\!F_n$ for each $n\geq 1$. 
For each $i\!\in\!\{1;2;...\, l\}$, there is an index $N_i$ such that 
$x^i_n\!\in\!Wx^i$ for any $n\geq N_i$. 
We set: $M=\max_{1\leq i\leq l} N_i$. 
For any $n\geq M$ and any $i\!\in\!\{1;2;...\, l\}$ we have: 
$x^i\!\in\!Wx^i_n$. 
Let $x\!\in\! F$. Then there is $i\!\in\!\{1;2;...\, l\}$ such that 
$x \!\in\! Wx^i$, so that for any $n\geq M$ we have: 
$x\!\in\! \Omega\cdot F_n$. 

\subsubsection{}
\label{sss - Chabauty implies geometric}
{\it Chabauty convergence implies geometric convergence}.~
Let us assume that $F_n$ converges to $F$ for the Chabauty topology. 

\vspace{2mmplus1mmminus1mm} 

Let $\varphi:{\bf N}\to{\bf N}$ be an increasing map and let 
$\{x_{\varphi(n)}\}_{n \geq 1}$ be a sequence in $G$ such that 
$x_{\varphi(n)} \!\in\! F_{\varphi(n)}$ for any $n\geq 1$, converging to $x\!\in\!G$. 
We choose a compact neighborhood $C$ of $x$ in $G$. 
There is an index $N_0\geq 1$ such that $x_{\varphi(n)}\!\in\!C$ for any $n\geq N_0$. 
Let us choose $\{\Omega_j\}_{j\geq 1}$ a decreasing sequence of 
compact symmetric neighborhoods of $1_G$ in $G$ 
such that $\bigcap_{j\geq 1} \Omega_j=\{1\}$. 
By definition of Chabauty convergence, there exists $N_1\geq N_0$ such that 
$(F,F_{\varphi(n)})\!\in\!P(C,\Omega_1)$ for any $n\geq N_1$. 
By induction, we find an increasing sequence of indices $\{N_j\}_{j\geq 1}$ such that 
$(F,F_{\varphi(n)})\!\in\!P(C,\Omega_j)$ for any $n\geq N_j$. 
Let $n\!\in\!{\bf N}$. 
There is a unique $j\geq 1$ such that $n\!\in\![N_J;N_{j+1}]$, and we can write: 
$x_{\varphi(n)}=\omega_n.y_n$ with $\omega_n\!\in\!\Omega_j$ and $y_n\!\in\!F$. 
Since $\Omega_j$ shrinks to $\{1\}$ as $j\to\infty$, we have: 
$\displaystyle \lim_{n\to\infty} \omega_n^{-1}.x_{\varphi(n)}=x$. 
Since $F$ is closed, this implies $x\!\in\!F$ and finally 
condition (i) for geometric convergence. 

\vspace{2mmplus1mmminus1mm} 

Let $x\!\in\!F$ and let $C$ and $\{\Omega_j\}_{j\geq 1}$ be as in the previous paragraph. 
As before, we can find an increasing sequence of indices $\{N_j\}_{j\geq 1}$ such that 
$(F,F_n)\!\in\!P(C,\Omega_j)$ for any $n\geq N_j$. 
Since $x\!\in\!F\cap C$, for any $n\!\in\![N_J;N_{j+1}]$ we can write: 
$x=\omega_n.y_n$, with $\omega_n\!\in\!\Omega_j$ and $y_n\!\in\!F_n$. 
Then since $\Omega_j$ shrinks to $\{1\}$ as $j\to\infty$, we have: 
$\displaystyle \lim_{n\to\infty} y_n=x$, which proves condition (ii) 
for geometric convergence. 
\qed

\subsection{Convergence of parahoric subgroups}
\label{ss - CV groups}
In this subsection, we prove that after taking stabilizers the fundamental sequences of points in the building $X$ (\ref{sss - canonical and fundamental}) lead to convergent sequences of parahoric subgroups in the Chabauty topology. 

\subsubsection{}
\label{sss - CV groups} 
The statement of the convergence result below uses the fact that the Euclidean buildings of the Levi factors of $G_F$ appear in the Bruhat-Tits building $X$ of $G_F$ (\ref{ss - Levi}). 
Recall that an $I$-canonical sequence $\{ x_n \}_{n\geq1}$ defines a facet in the Bruhat-Tits building $X_I$ of the Levi factor 
$G_I{}_{/F}$ (\ref{sss - canonical and fundamental}). 

\begin{theorem}
\label{th - CV groups} 
Let $I$ be a proper subset of $S$. 
Let $\{x_n\}_{n\geq 1}$ be an $I$-canonical sequence of points in the closed Weyl chamber $\overline{\mathscr{Q}}^X$ of $X$. 
Let $\sigma$ be the facet in the Bruhat-Tits building $X_I$ defined by $\{x_n\}_{n\geq 1}$. 
Let $\underline d=\{d_s\}_{s\in I}$ be any family of real non-negative numbers defining a point of $\sigma$. 
Then the sequence of parahoric subgroups $\{K_{x_n}\}_{n\geq1}$ converges in $\mathscr{S}(G_F)$ 
to the closed subgroup $D_{I,\underline d}$. 
\end{theorem}

Since $\overline{\mathscr{Q}}^X$ is a fundamental domain for the $K_o$-action on $X$ \cite[4.4.3 (2)]{BrT1}, we readily deduce the following consequence. 

\begin{cor}
\label{cor - CV groups}
Let $\{x_n\}_{n\geq 1}$ be a fundamental sequence in the Bruhat-Tits building $X$. 
Then the corresponding sequence of parahoric subgroups $\{K_{x_n}\}_{n\geq1}$ converges in $\mathscr{S}(G_F)$ to a $K_o$-conjugate of some subgroup $D_{I,\underline d}$. 
\qed\end{cor} 

The rest of the subsection is devoted to proving Theorem \ref{th - CV groups}. 
Let $\underline d=\{d_s\}_{s\in I}$ be a family of non-negative parameters whose associated point in the Weyl chamber 
$A_I \cap p_I(\overline{\mathscr{Q}}^X)$ lies inside the facet $\sigma$.  
By compactness of the Chabauty topology on $\mathscr{S}(G_F)$, it is enough to show that $D_{I,\underline d}$ is the only cluster value of $\{K_{x_n}\}_{n\geq1}$. 
Let $\{K_{x_{\psi(n)}}\}_{n\geq1}$ be a subsequence converging to some closed subgroup $D$ of $G_F$.

\subsubsection{}
\label{sss - CV probas}
Let us first prove some measure-theoretic results. 
Recall that if $Y$ is a compact and metrizable topological space, then so is the weak-$*$ topology on the space of probability measures $\mathscr{M}^1(Y)$ by the Banach-Alaoglu-Bourbaki theorem. 
We apply this to the case when $Y$ is a flag variety of some semisimple $F$-group. 
For $x$ in the apartment $A$, we denote by $\mu_x$ the unique $K_x$-invariant probability measure on the Furstenberg boundary $\mathscr{F}$, supported by the compact $K_x$-homogeneous subspace $K_x.\omega$. 

\begin{lemma} 
\label{lemma - big cell conegligible}
For each point $x\!\in\!A$, the big cell is conegligible, i.e. we have: $\mu_x\bigl((U^-)_F.\omega\bigr)=1$. 
\end{lemma} 

{\it Proof}.~
Let ${\rm d}k_x$ be the Haar measure of total mass 1 on the compact group $K_x$. 
Let $p:G_F\to\mathscr{F}$ denote the orbit map $g\mapsto g.\omega$. 
It is enough to show that the volume of 
$(p\!\mid_{K_x})^{-1}\bigl(K_x.\omega\setminus(U^-)_F.\omega\bigr)$ 
with respect to ${\rm d}k_x$ is zero \cite[Lemma 1.3]{Mackey}. 
On the one hand, since $K_x$ is an open subgroup of $G_F$, we have: ${\rm d}g\!\mid_{K_x}=C\cdot{\rm d}k_x$, where $C$ is a multiplicative constant $>0$ and ${\rm d}g$ is a  Haar measure on $G_F$. 
On the other hand, we have: 
$(p\!\mid_{K_x})^{-1}\bigl(K_x.\omega\setminus(U^-)_F.\omega\bigr)\subset G_F\setminus(U^-\cdot P)_F$. 
Therefore we finally obtain: 

\vspace{2mmplus1mmminus1mm} 

\centerline{
${\rm Vol}\bigl((p\!\mid_{K_x})^{-1}(K_x.\omega\setminus(U^-)_F.\omega),{\rm d}k_x\bigr) 
\leq C \cdot {\rm Vol}(G_F\setminus(U^-\cdot P)_F,{\rm d}g)=0$,}

\vspace{2mmplus1mmminus1mm} 

since small Bruhat cells are negligible for any Haar measure on $G_F$. 
This proves that the big cell $(U^-)_F.\omega$ is $\mu_x$-conegligible. 
\qed\vspace{2mmplus1mmminus1mm} 

Let us denote by  $E_{I,\sigma}$ the intersection of affine half-spaces $p_I^{-1}(\sigma)$ in the apartment $A$ (\ref{sss - Levi BT buildings}). 

\begin{prop}
\label{prop - CV probas} 
Let $\{x_n\}_{n\geq 1}$ be an $I$-canonical sequence defining the facet $\sigma$ of $X_I$. 
Let $\nu$ be a cluster value of $\{\mu_{x_n}\}_{n\geq1}$. 
\begin{enumerate}
\item [(i)] We have: $\overline{{\rm Supp}(\nu)}^Z=G_I.\omega$ and $R_{I,\underline d}<{\rm Stab}_{G_F}(\nu)$. 
In particular, we have: ${\rm supp}(\nu)\subset \mathscr{F}_I$, where $\mathscr{F}_I$ denotes the copy $(G_I)_F.\omega$ of the Furstenberg boundary of $(G_I)_F$. 
\item [(ii)] If $\sigma$ is a vertex of $X_I$ and if $\underline d=\{d_s\}_{s\in I}$ is the family of non-negative real numbers defining $\sigma$, then we have: ${\rm Stab}_{G_F}(\nu)=R_{I,\underline d}$. 
\end{enumerate}
\end{prop}

{\it Proof}.~
Let $\displaystyle \nu = \lim_{n\to\infty} \mu_{x_{\psi(n)}}$ be a cluster value as in the statement.  
We may -- and shall -- assume that each point $x_{\psi(n)}$ belongs to $E_{I,\sigma}$. 
The group $(T^I)_F$ acts as a cocompact translation group on $E_{I,\sigma}$, so there is a compact complete system of representatives $Y$ for the $(T^I)_F$-action on $E_{I,\sigma}$. 
We can write $x_{\psi(n)}=t_n.y_n$ with $t_n\!\in\!T^I$ and $y_n\!\in\!Y$ for each $n\geq 1$. 
Since $Y$ is contained in finitely many facets of $E_{I,\sigma}$, up to extracting again, we may -- and shall -- assume that there is a facet $\tau$ such that $\tau\subset Y$ and $y_n\!\in\!\tau$ for each $n\geq1$. 
Moreover, by uniqueness we have: $\mu_{t_n.y_n}=t_n{}_*\mu_\tau$, where $\mu_\tau$ is the unique probability measure invariant under $K_\tau={\rm Stab}_{G_F}(\tau)$ with ${\rm supp}(\mu_\tau)=K_\tau.\omega$. 

\vspace{2mmplus1mmminus1mm} 

Let us set: $R={\rm Stab}_{G_F}(\nu)$. 
The groups $K_{I,\underline d}$ and $T^I$ commute with one another. 
Moreover by \cite[6.4.9]{BrT1} applied to $K_{I,\underline d}$ in $(M_I)_F$ and to $K_\tau$ in $G_F$, we have: 
$K_{I,\underline d}<K_\tau$. 
Therefore for any $k\!\in\!K_{I,\underline d}$, we have: $k_*(t_n{}_*\mu_\tau)=t_n{}_*(k_*\mu_\tau)=t_n{}_*\mu_\tau$. 
By passing to the limit as $n\to+\infty$, the previous paragraph implies that we have: $K_{I,\underline d}<R$. 

\vspace{2mmplus1mmminus1mm} 

Let us pick $a\!\in\!\Phi^{I,-}$. 
The increasing family of affine half-apartments $\{t_n.D_{a,m}\}_{n\geq 1}$ exhausts $A$. 
Since $t_n.D_{a,m}$ is the set of fixed-points of $t_nU_{a,m}t_n^{-1}$ in $A$ (\ref{sss - valuated root data and half spaces}), we deduce that the valuation of elements in $t_nU_{a,m}t_n^{-1}$ goes to $+\infty$ as $n\to+\infty$. 
Let $U^{I,-}$ be the unipotent radical of the parabolic subgroup opposite $P_I$ with respect to $T$. 
The product map provides a $T$-equivariant isomorphism of algebraic varieties $\prod_{a\in\Phi^{I,+}}U_{-a} \simeq U^{I,-}$. 
By the previous remark, $\displaystyle\lim_{n\to+\infty}{\rm dist}_X(\Pi^s,x_m)=+\infty$ for each $s\!\in\!S\setminus I$ implies at the group level that $\{t_n\}_{n\geq1}$ is a contracting sequence on $(U^{I,-})_F$, i.e. for any compact subset $C$ in $(U^{I,-})_F$ and for any neighborhood $\Omega$ of the identity in $(U^{I,-})_F$, we have: $t_nCt_n^{-1} \subset \Omega$ for $n>\!\!>1$. 

\vspace{2mmplus1mmminus1mm} 

Let $f$ be a continuous function on $\mathscr{F}$ that is vanishing on $\mathscr{F}_I$. 
By definition of weak-$*$ convergence, we have: 

\vspace{2mmplus1mmminus1mm} 

\centerline{$\displaystyle \nu(f)=\lim_{n\to\infty} \int_{\mathscr{F}} f(t_nz)\,{\rm d}\mu_\tau(z)$.}

\vspace{2mmplus1mmminus1mm} 

But by Lemma \ref{lemma - big cell conegligible}, for each index $n\geq1$ we have: 

\vspace{2mmplus1mmminus1mm} 

\centerline{$\displaystyle \int_{\mathscr{F}} f(t_nz)\,{\rm d}\mu_\tau(z)=\int_{(U^-)_F.\omega} f(t_nz)\,{\rm d}\mu_\tau(z)$.}

\vspace{2mmplus1mmminus1mm} 

For each $z\!\in\!(U^-)_F.\omega$, the sequence $\{t_n.z\}_{n\geq1}$ converges to some point in $(U^-_I)_F.\omega$, so the sequence of functions $\{f\circ t_n\}_{n\geq 1}$ simply converges to 0 on the big cell. 
Since $\mu_\tau$ is of total mass 1, Lebesgue dominated convergence theorem implies: 
$\displaystyle \lim_{n\to\infty} \int_{(U^-)_F.\omega} f(t_nz)\,{\rm d}\mu_\tau(z)=0$, so $\nu(f)=0$ for each $f\!\in\!C(\mathscr{F})$ vanishing on $\mathscr{F}_I$. 
This finally implies: ${\rm supp}(\nu)\subset \mathscr{F}_I$. 

\vspace{2mmplus1mmminus1mm} 

At this stage, we already know that $\nu$ is fixed by $K_{I,\underline d}$ and supported on $\mathscr{F}_I$. 
This implies that ${\rm supp}(\nu)$ is a (finite) union of $K_{I,\underline d}$-orbits in $\mathscr{F}_I$. 
These orbits are finite in number and the pull-back of each of them under the orbit map $(G_I)_F\to \mathscr{F}_I$ of $\omega$ is Zariski dense in $G_I$ \cite[4.2.1]{BrT1}. 
This implies that $\overline{{\rm Supp}(\nu)}^Z=G_I.\omega$, which is the first statement of (i). 
We also have: $(T^I\ltimes U^I)_F<R$ because $(T^I\ltimes U^I)_F$ fixes pointwise $\mathscr{F}_I$. 
This implies that $R_{I,\underline d}=K_{I,\underline d}\cdot (T^I)_F\cdot (U^I)_F<R$, so that (i) is now proved. 

\vspace{2mmplus1mmminus1mm} 

Note that conversely, we have: $R<{\rm Stab}_{G_F}\bigl(\overline{{\rm Supp}(\nu)}^Z\bigr)$, which by (i) implies $R<(P_I)_F$. 

\vspace{2mmplus1mmminus1mm} 

(ii). We now assume that $\sigma$ is a vertex, i.e. that $G_I\cap K_{I,\underline d}$ is a maximal compact subgroup of $(G_I)_F$. 
We already know that: $R_{I,\underline d}<R<(P_I)_F$. 
Let us assume that there exists some $r\!\in\!R\setminus R_{I,\underline d}$. 
As an element of $(P_I)_F$, this element can be written $r=gtu$ with $g\!\in\!(G_I)_F$, $t\!\in\!(T^I)_F$ and $u\!\in\!(U^I)_F$. 
Since  $tu\!\in\!R$, we have: $g\!\in\!(G_I\cap R)\setminus R_{I,\underline d}$, which implies that $G_I\cap R$ is strictly bigger than the
parahoric subgroup $G_I \cap K_{I,\underline d}$ of $(G_I)_F$. 
In view of the lattice structure of parahoric subgroups in the affine Tits system of $(G_I)_F$, this implies that $G_I \cap K_{I,\underline d}$ contains a full simple factor of $(G_I)_F$. 
The projection of $\nu$ on the flag variety of this factor would lead to an invariant probability measure on the flag variety of a non-compact semisimple group. 
Since minimal parabolic subgroups are not unimodular, this is impossible \cite[Lemma 1.4]{Raghu}. 
We finally have: $R=R_{I,\underline d}$. 
\qed 

\subsubsection{}
\label{sss - proof CV groups} 
We now turn to the proof of the above convergence theorem (Theorem \ref{th - CV groups}). 
With the notation of \ref{sss - CV groups}, it is enough to show that $D=D_{I,\underline d}$, where $\displaystyle D=\lim_{n\to\infty} K_{x_{\psi(n)}}$. 
We start with a lower bound for $D$ with respect to the inclusion relation on closed subgroups in $G_F$. 

\begin{lemma}
\label{lemma - lower bound cluster value}
The cluster value $D$ contains $D_{I,\underline d}$. 
\end{lemma}

{\it Proof}.~We first show that the group $D$ necessarily contains the unipotent radical $(U^I)_F$. 
Since the product map provides a bijection: $\prod_{a\in\Phi^{I,+}}U_a\simeq U^I$ \cite[6.1.6]{BrT1}, it is enough to show that $u\!\in\!D$ for any  $a\!\in\!\Phi^{I,+}$ and any $u\!\in\!(U_a)_F$. 
Let $a\!\in\!\Phi^{I,+}$ and let $u\!\in\!(U_a)_F$. 
By definition of a valuated root datum \cite[6.2.1]{BrT1}, $u$ belongs to a subgroup $U_{a,\varphi_a(u)}$ of the filtration of $(U_a)_F$ given by $\varphi_a$. 
Moreover there is a half-space $D_{a,\varphi_a(u)}$ of the apartment $A$, containing a translate of $\mathscr{Q}$, fixed by  $U_{a,\varphi_a(u)}$ and bounded by a wall whose direction is transverse to $L_I=\langle\mathscr{Q}^I\rangle$. 
Since $\displaystyle \lim_{n\to+\infty}d_X(x_{\psi(n)},\Pi^s)=+\infty$ for each $s\!\in\!S\setminus I$, there is an index $N\geq 1$ such that for any $n\geq N$ we have: $x_{\psi(n)}\!\in\!D_{a,\varphi_a(u)}$. 
This is the geometric translation of the fact that $K_{x_{\psi(n)}}$ contains $U_{a,\varphi_a(u)}$, hence $u$, for $n\geq N$. 
This enables to see $u$ as the limit of the sequence $\{ g_{\psi(n)}\}_{n\geq N}$ with $g_{\psi(n)}=u\!\in\!K_{x_n}$ for each $n\geq N$. 
By definition of geometric convergence, this implies $u\!\in\!D$. 

\vspace{2mmplus1mmminus1mm} 

We now show that the closed subgroup $D$ necessarily contains the compact group $K_{I,\underline d}$. 
Let $a$ be a root in $\Phi(G_I,T_I)$: the direction $\partial a$ contains $L_I$. 
Let $m\!\in\!{\bf Z}$ be such that $\overline{D_{a,m}}\supset p_I^{-1}(\sigma)$. 
For $n>\!\!>1$ the point $x_{\psi(n)}$ lies in the closed half-apartment $\overline{D_{a,m}}$, so that $U_{a,m}<K_{v_\psi(n)}$. 
As in the previous paragraph, this shows that $D$ contains $U_{a,m}$. 
Since the groups $U_{a,m}$ with $a$ and $m$ as above generate the parahoric subgroup $G_I \cap K_{I,\underline d}$ of $G_I$ 
\cite[6.4.9]{BrT1}, we obtain: $G_I \cap K_{I,\underline d}<D$. 
Similarly, we see that: ${\rm Fix}_{T^I}(p_I^{-1}(\sigma))={\rm Fix}_{T^I}(L_I)$ lies in $D$. 
Finally, we have: $K_{I,\underline d}<D$ since 
$K_{I,\underline d}=(G_I \cap K_{I,\underline d})\cdot{\rm Fix}_{T^I}(L_I)$. 
\qed\vspace{2mmplus1mmminus1mm} 

Thanks to the measure-theoretic results of \ref{sss - CV probas}, we also have an upper bound for $D$ with respect to the inclusion relation on closed subgroups in $G_F$. 

\begin{lemma}
\label{lemma - upper bound cluster value}
The cluster value $D$ is contained in $R_{I,\underline d}$. 
\end{lemma}

{\it Proof}.~
Up to extracting again in order to have a convergent sequence of probability measures as in \ref{sss - CV probas}, we may -- and shall -- assume that: 
$\displaystyle \lim_{n\to\infty}K_{x_{\psi(n)}}=D$ and $\displaystyle \lim_{n\to\infty}\mu_{x_{\psi(n)}}=\nu$. 
Lebesgue's dominated convergence theorem then implies: $D<{\rm Stab}_{G_F}(\nu)$ \cite[Lemma 9.7]{GJT}. 
If $\sigma$ is a vertex in $X_I$, it remains to use Proposition  \ref{prop - CV probas} to conclude. 

\vspace{2mmplus1mmminus1mm} 

Otherwise, we note that $R_{I,\underline d}$ is the intersection of the groups $R_{I,\underline d'}$ where 
$\underline d'$ varies over the families of parameters defining a vertex in the closure of $\sigma$. 
Let us fix such a family of parameters $\underline d'$, defining a vertex $v\!\in\!\overline\sigma$. 
Then for each $n\geq 1$ there exists an element $x'_{\psi(n)}$ in the intersection of the closure of the facet containing $x'_{\psi(n)}$ and of $p_I^{-1}(v)$. 
For each $n\geq 1$, we have: $K_{x_{\psi(n)}}<K_{x'_{\psi(n)}}$, and up to extracting we may -- and shall -- assume that 
$\{ K_{x'_{\psi(n)}}\}_{n\geq1}$ converges for the Chabauty topology to some closed subgroup $D'<G_F$.
This group contains $D$, and by the first paragraph dealing with the case of vertices in $X_I$, we have: $D'<R_{I,\underline d'}$. 
The conclusion follows by intersecting over the vertices $v$ in $\overline\sigma$. 
\qed\vspace{2mmplus1mmminus1mm}

The previous two lemmas show that we have: 

\vspace{2mmplus1mmminus1mm} 

\centerline{$(*)$ \hspace{1cm} 
$D_{I,\underline d}=K_{I,\underline d}\ltimes (U^I)_F<D<(K_{I,\underline d}\cdot T^I)\ltimes (U^I)_F=R_{I,\underline d}$,}

\vspace{2mmplus1mmminus1mm}

so it remains to show that the cluster value $D=\displaystyle\lim_{n\to\infty}K_{v_{\psi(n)}}$ cannot be bigger than 
$K_{I,\underline d}\ltimes (U^I)_F$. 
By $(*)$, it is enough to show that $T^I \cap D<K_{I,\underline d}$. 
Let $G<{\rm GL}_m$ be an embedding of $F$-algebraic groups. 

\begin{lemma}
\label{lemma - limit is distal}
Let $t\!\in\!D$. 
Then any eigenvalue of $t$ has absolute value $1$. 
\end{lemma} 

{\it Proof}.~By \cite[Proposition 1.12]{PlatonovRapinchuk}, for any $n\geq1$ there exists $g_n\!\in\!{\rm GL}_m(F)$ such that we have: 
$K_{x_{\psi(n)}}<g_n{\rm GL}_m(\mathscr{O}_F)g_n^{-1}$. 
Using the definition of geometric convergence, we write: $t=\displaystyle \lim_{n\to\infty}k_n$ with $k_n\!\in\!K_{x_{\psi(n)}}$ for each $n\geq1$. 
Denoting by $\chi_M(x)$ the characteristic polynomial of a matrix $M\!\in\!M_{m\times m}(F)$, we have: $\chi_{k_n}(x)\!\in\!\mathscr{O}_F[x]$ for any $n\geq1$. 
Therefore, by passing to the limit we obtain: $\chi_t(x)\!\in\!\mathscr{O}_F[x]$. 
Let $v_0$ denote the minimal valuation over ${\rm Sp}(t)$, the set of eigenvalues of $t$ counted with multiplicities, and let $l$ denote the number of occurencies of $v_0$ in ${\rm Sp}(t)$. 
We set: $\chi_t(x)=x^m+a_1x^{m-1}+\dots+a_m$. 
Since $F$ is non-archimedean, we have: $v(a_l)=l.v_0$, and since $\chi_t(x)\!\in\!\mathscr{O}_F[x]$, we obtain: $v_0\geq 0$, i.e. 
${\rm Sp}(t)\subset\mathscr{O}_{\overline F}^\times$. 
\qed\vspace{2mmplus1mmminus1mm} 

We can now conclude the proof of Theorem \ref{th - CV groups}. 

\vspace{2mmplus1mmminus1mm} 

{\it Proof}.~
We keep the notation of the previous lemma. 
Let $t\!\in\!T^I\cap D$. 
Since $t$ is a diagonalizable matrix in ${\rm GL}_m(F)$, Lemma \ref{lemma - limit is distal} implies that 
the subgroup $\langle t \rangle$ is relatively compact in $G_F$. 
By the Bruhat-Tits fixed-point lemma \cite[Lemme 3.2.3]{BrT1}, 
this implies that $t$ fixes a point, say $y$, in the building $X$. 
Let us call $x$ the orthogonal projection of $y$ onto the closed convex subset $\langle\mathscr{Q}^I\rangle$ of $X$ 
\cite[II Proposition 2.4]{BriHae}. 
Since $L_I=\langle\mathscr{Q}^I\rangle$ is $t$-stable because $t\!\in\!T^I$, the uniqueness of 
$x$ implies that $t$ fixes $x$. 
Since $t$ acts as a translation on the affine space $L_I$, we deduce that
$t$ fixes $L_I$ pointwise. 
Finally, we have: $t\!\in\! K_{I,\underline d}$, and in view of the previous reductions 
this implies that $D_{I,\underline d}=K_{I,\underline d}\ltimes (U^I)_F$ is the only cluster value of 
$\{K_{v_n}\}_{n\geq1}$. 
This proves the desired convergence. 
\qed

\subsection{Compactification of the vertices}
\label{ss - Guivarc'h}
Let $\mathscr{K}(G_F)$ denote the space of maximal compact subgroups in $G_F$. 
Let $K:x\mapsto K_x$ denote the map assigning to a point $x\!\in\!X$ its stabilizer $K_x$, 
i.e. the parahoric subgroup ${\rm Stab}_{G_F}(x)$. 

\subsubsection{}
The map $(K\!\mid_{V_X})^{-1}$ establishes the one-to-one correspondence between maximal compact subgroups of $G_F$ and vertices of its Euclidean building $X$, given by the Bruhat-Tits fixed-point lemma \cite[Lemme 3.2.3]{BrT1}. 

\begin{prop}
\label{prop - Guivarc'h}
The restriction 

\vspace{2mmplus1mmminus1mm} 

\centerline{
$\begin{array}{rrrr} 
\hfill K\!\mid_{V_X}: \hfill & \hfill V_X \hfill & \hfill\to\hfill &
\hfill \mathscr{S}(G_F) \hfill\\ 
& \hfill v \hfill       & \hfill\mapsto\hfill & \hfill K_v \hfill 
\end{array}$
}

\vspace{2mmplus1mmminus1mm} 

of the above map $K$ to the vertices of the Bruhat-Tits building $X$ is a $G_F$-equivariant topological embedding of the set of vertices $V_X$ into the space of closed subgroups $\mathscr{S}(G_F)$ endowed with the Chabauty uniform structure. 
\end{prop}

We can rephrase the proposition by saying that the subset $\mathscr{K}(G_F)$ of maximal compact subgroups of $G_F$ is discrete for geometric convergence (\ref{ss - Chabauty}). 

\vspace{2mmplus1mmminus1mm} 

{\it Proof}. 
The $G_F$-equivariance of $K\!\mid_{V_X}$ is obvious, and so is the continuity since $V_X$ is discrete.  
The injectivity is also clear since $X^{K_v}=\{v\}$ for every vertex $v$. 
It remains to prove that $(K\!\mid_{V_X})^{-1}$ is continuous. 
Let $\{v_n \}_{n\geq1}$ be a sequence of vertices such that the corresponding sequence
$\{K_{v_n}\}_{n\geq1}$ of maximal compact subgroups converges to some maximal compact subgroup $K_v$. 
We have to show that $\displaystyle \lim_{n\to+\infty}v_n=v$. 

\vspace{2mmplus1mmminus1mm} 

First, $\{v_n \}_{n\geq1}$ is a bounded sequence in $V_X$ since otherwise we could extract a 
subsequence $\{v_{\varphi(n)} \}_{n\geq1}$ going to infinity. 
Then we could extract again a subsequence as in Lemma \ref{lemma - fundamental subsequence}: this would lead to a
contradiction  with $\displaystyle \lim_{n\to+\infty} K_{v_n}=K_v$ since by Theorem \ref{th - CV groups} the latter subsequence would converge to an unbounded limit group (\ref{ss - CV groups}). 

\vspace{2mmplus1mmminus1mm} 

Now assume $v'=\displaystyle \lim_{n\to+\infty}v_{\varphi(n)}$ for some increasing map 
$\varphi:{\bf N}\to{\bf N}$. 
By continuity we have: $\displaystyle \lim_{n\to+\infty}K_{v_{\varphi(n)}}=K_{v'}$. 
But the assumption $\displaystyle \lim_{n\to+\infty} K_{v_n}=K_v$ implies that $K_v=K_{v'}$, hence 
$v=v'$ by injectivity of $K\!\mid_{V_X}$. 
This shows that $v$ is the only cluster value of $\{v_n \}_{n\geq1}$, 
which finally proves the continuity of 
$(K\!\mid_{V_X})^{-1}$. 
\qed

\subsubsection{}
We can finally define the desired compactification of the set of vertices $V_X$. 

\vspace{2mmplus1mmminus1mm} {\bf Definition.~}\it 
\label{def - Guivarc'h}
The closure of $\mathscr{K}(G_F)$ in the compact space $\mathscr{S}(G_F)$ is called the {\rm group-theoretic compactification}~of the building $X$. 
We denote it by $\overline{V}^{\rm gp}_X$. 
\rm\vspace{2mmplus1mmminus1mm} 

We note that for an arbitrary linear algebraic group $G_{/F}$ and for $F$ of characteristic 0, we have: 
$\mathscr{K}(G_F) \neq \varnothing$ if, and only if, $G_{/F}$ is reductive 
\cite[Proposition 3.15]{PlatonovRapinchuk}. 

\subsubsection{}
\label{sss - limit groups} 
We define the {\it limit groups~}to be the cluster values of the sequences $\{K_{x_n}\}_{n\geq1}$ where $\{ x_n\}_{n\geq1}$ is a sequence in the building $X$ going to infinity. 
We denote by $\partial^{\rm gp}X$ the set of limit groups. 

\begin{cor}
\label{cor - limit groups} 
The set $\partial^{\rm gp}X$ consists of the closed subgroups $kD_{I,\underline d}k^{-1}$ when $k$ varies in $K$, $I$ varies over the proper subsets of the set of simple roots $S$, and $\underline d$ varies over the families of non-negative real numbers indexed by $I$. 
\end{cor}

{\it Proof}.~
This is an easy combination of Lemma \ref{lemma - fundamental subsequence} and Corollary \ref{cor - CV groups}. 
\qed\vspace{2mmplus1mmminus1mm} 

The assignment $x\mapsto K_x$ is not a continuous map from $X$ to $\mathscr{S}(G_F)$. 
Therefore geometrically, i.e. as far as compactifications for $X$ are involved, the only relevant limit groups are those arising from sequences of maximal parahoric (equivalently, maximal compact) subgroups. 

\section{Geometric description of the compactification}
\label{s - description}

We describe more precisely the compactification obtained in the previous section. 
We first compute all stabilizers and single out a closed orbit $G_F$-equivariantly homeomorphic to the flag variety $\mathscr{F}$. 
We also describe precisely the boundary of these compactifications. 
In the case when ${\rm rk}_F(G)\geq2$, this provides a major difference between the geometric compactification and the group-theoretic one. 
Finally we compare the group-theoretic compactification with the polyhedral one; since for vertices they are the same, we can extend the group-theoretic compactification $\overline{V}^{\rm gp}_X$ to a compactification of the full Bruhat-Tits building $X$. 

\subsection{Stabilizers and orbits}
\label{ss - stabilizers and orbits}
We compute the isotropy groups and describe a specific closed orbit for the $G_F$-action on the boundary of $\overline{V}^{\rm gp}_X$. 

\subsubsection{}
\label{sss - Z density}
Let us start with a lemma on Zariski closures of limit groups and of their normalizers. 
The result is used to compute stabilizers. 

\begin{lemma}
\label{lemma - Z density} 
Let $I$ be a subset of $S$ and let $\underline d$ be a family of non-negative real numbers indexed by $I$. 
The Zariski closure in $G$ of any limit group $D_{I,\underline d}$, hence of any $R_{I,\underline d}$, is equal
to the parabolic $F$-subgroup $P_I$. 
Therefore the family of the Zariski closures of the limit groups, or of their stabilizers, coincides with the family of the parabolic proper $F$-subgroups of $G$. 
\end{lemma} 

{\it Proof}.~
Let $\mathscr{B}$ be the standard Iwahori subgroup. 
The rational points of the $F$-subgroup $\overline{\mathscr{B}}^Z$ contain $\mathscr{B}$, so they are of positive measure for any Haar measure on $G_F$. 
By \cite[2.5.3]{Margulis}, this implies $\overline{\mathscr{B}}^Z=G$. 
This shows that if $H$ is a semisimple simply connected group over $F$, then the fixator of any facet in the Bruhat-Tits building of $H_{/F}$ is Zariski dense in $H$. 
It remains to apply this fact in various Levi factors. 
Indeed, given $D_{I,\underline d}$, we have: $\overline{G_I \cap D_{I,\underline d}}^Z=\overline{G_I \cap K_{I\underline d}}^Z=G_I$. 
Together with the Zariski density of the unique maximal compact subgroup of $(T^I)_F$ in $T^I$, this finally implies: 
$\overline{D_{I,\underline d}}^Z=P_I$. 
The equality $\overline{R_{I,\underline d}}^Z=P_I$ is then clear since $D_{I,\underline d}<R_{I,\underline d}<(P_I)_F$. 
\qed\vspace{2mmplus1mmminus1mm}

We see here a difference with the case of a semi-simple group over the real numbers. 
In the latter case, a compact semi-group is the group of real points of an algebraic ${\bf R}$-subgroup. 
This implies that in the case of symmetric spaces, (the stabilizer of) any limit group is the group of real points of an algebraic ${\bf R}$-subgroup (in a suitable proper parabolic ${\bf R}$-subgroup). 

\subsubsection{}
\label{sss - stabilizers} 
We can now compute the normalizer of each limit group for the $G_F$-action by conjugation on the space of closed subgroups $\mathscr{S}(G_F)$. 
This is slightly more general than computing the stabilizers of the points of the $G_F$-space $\overline{V}^{\rm gp}_X$. 

\begin{lemma}
\label{lemma - stabilizers} 
Let $I$ be a subset of $S$ and let $\underline d$ be a family of non-negative real numbers indexed by $I$. 
We have: $N_{G_F}(D_{I,\underline d})=R_{I,\underline d}$. 
\end{lemma} 

{\it Proof}.~
Let us set $R=N_{G_F}(D_{I,\underline d})$. 
On the one hand, since conjugation by any $g\!\in\!G_F$ is an algebraic automorphism of $G$, any $r\!\in\!R$ stabilizes the Zariski closure $\overline{D_{I,\underline d}}^Z$. 
So by Lemma \ref{lemma - Z density} and the fact that a parabolic subgroup is equal to its own normalizer \cite[11.16]{Borel}, 
we obtain: $R<(P_I)_F$. 
On the other hand, the group $(T^I)_F$ centralizes $K_{I,\underline d}$ and normalizes $U^I$, so $R$ contains 
$R_{I,\underline d}=D_{I,\underline d}\cdot (T^I)_F$. 
At this stage we have: $R_{I,\underline d}<R<(P_I)_F$. 
It remains to show that $R$ is not bigger than $R_{I,\underline d}$. 
Since $P_I=(G_I\cdot T^I)\ltimes U^I$ and $(T^I\ltimes U^I)_F<R$, it is enough to show that $R\cap G_I$ is equal to 
$K_{I,\underline d}\cap G_I$. 
The latter group is a parahoric subgroup of $(G_I)_F$, i.e. a parabolic subgroup of the affine Tits system of parahoric subgroups in 
$(G_I)_F$ \cite[\S 2]{BrT1}.
Therefore $K_{I,\underline d}\cap G_I$ is equal to its normalizer in $(G_I)_F$ \cite[IV.2 Proposition 4]{Lie456}. 
The group $R\cap G_I$ normalizes $D_{I,\underline d}\cap G_I$. 
Since $D_{I,\underline d}\cap G_I=K_{I,\underline d}\cap G_I$, this finally proves: $R\cap G_I<K_{I,\underline d}\cap G_I$. 
\qed

\subsubsection{}
\label{sss - closed orbit} 
The $G_F$-orbit described below will turn out to be the unique closed orbit in $\overline{V}^{\rm gp}_X$. 

\begin{lemma}
\label{lemma - closed orbit} 
The $G_F$-orbit of the limit group $D_\varnothing$ is closed and $G_F$-homeomorphic to the Furstenberg boundary $\mathscr{F}$. 
\end{lemma}

{\it Proof}.~
We have: $D_\varnothing=K_\varnothing \ltimes U_F$ with
$K_\varnothing=[M,M]_F\cdot T_{\rm cpt}$ and $U=\mathscr{R}_u(P)$, and where $T_{\rm cpt}$ is the unique maximal compact subgroup of $T_F$. 
The orbit map $G_F \to \mathscr{S}(G_F)$ attached to $D_\varnothing$ under the $G_F$-action by
conjugation, i.e. the map $g \mapsto gD_\varnothing g^{-1}$ factorizes through the quotient map 
$G_F \to G_F/P_F$, since by Lemma \ref{lemma - stabilizers} for $I=\varnothing$ we have: $N_{G_F}(D_\varnothing)=P_F$. 
The conclusion follows from the compactness of $\mathscr{F}$.
\qed

\subsection{Euclidean buildings in the boundary}
\label{ss - boundary}
In this subsection we fix $I$ a proper subset of the set of simple roots $S$ of $\Phi=\Phi(T,G)$. 
The choices of $T$, $P$ and $I$ determine a standard parabolic $F$-subgroup $P_I$, together with the
reductive Levi factor $M_I$ and the semisimple Levi factor $G_I=[M_I,M_I]$. 
Recall that we have: $P_I=M_I\ltimes U^I=(G_I\cdot T^I)\ltimes U^I$. 
We denote by $(T^I)_{\rm cpt}$ the unique maximal compact subgroup of the $F$-split torus $(T^I)_F$. 

\subsubsection{}
\label{sss - definitions in the boundary} 
On the one hand, we can introduce the Bruhat-Tits building $X_I$ of the semisimple $F$-group $G_I$.  
We can also apply the compactification procedure of \ref{ss - Guivarc'h} to $X_I$: we obtain the group-theoretic compactification $\overline V_{X_I}^{\rm gp}$. 
This is a closed subset of the compact set $\mathscr{S}\bigl((G_I)_F\bigr)$ of closed subgroups of $(G_I)_F$. 
In this situation limit groups are given, up to conjugation, by subsets $J$ of the set of simple roots $I$ of $G_I$ and families of non-negative real parameters $\underline d$ indexed by $J$. 
By Theorem \ref{th - CV groups}, a maximal limit group is of the form 
$D_{J,\underline d}\cap G_I=(K_{J,\underline d}\cap G_I)\ltimes(U^J\cap G_I)_F$, where $K_{J,\underline d}$ is the maximal compact subgroup of the reductive Levi factor $(M_J)_F$ which is determined by the parameters $\underline d$. 
We denote such a limit group by $D_{J\subset I,\underline d}$. 
To sum up, the boundary of $\overline V_{X_I}^{\rm gp}$ is the following set of maximal limit groups in $(G_I)_F$: 

\vspace{2mmplus1mmminus1mm} 

\centerline{
$\{k.D_{J\subset I,\underline d}.k^{-1} : k\!\in\!(K_o\cap G_I)$, $J\subset I$ and $\underline d$
non-negative real parameters indexed by $J\}$.} 

\vspace{2mmplus1mmminus1mm} 

On the other hand, by \cite[7.6]{BrT1} there is a non-essential realization of $X_I$ 
in the building $X$ of $G_{/F}$ (\ref{ss - Levi}). 
The apartment $A$ (resp. the Weyl chamber $\mathscr{Q}$, the alcove $c$) in $X$ determines an apartment 
$A_I$ (resp. a Weyl chamber $\mathscr{Q}_I$, an alcove $c_I$) in $X_I$. 
The vertices of the standard alcove $c_I$ are affine subspaces of the apartment $A$; 
they determine limit groups in the closure of the Weyl chamber $\mathscr{Q}$. 
We denote the latter groups by 
$D_{I,\underline d_0}$, $D_{I,\underline d_1}$, ... $D_{I,\underline d_{\mid I \mid}}$, and 
we denote by $Y_I$ the union of the $(G_I)_F$-orbits of the limit groups 
$D_{I,\underline d_i}$ in $\mathscr{S}(G_F)$ when $i$ ranges over $\{0;1;...\mid\! I\!\mid\}$. 

\vspace{2mmplus1mmminus1mm} 

We want to show that the closure of $Y_I$ in $\overline V_X^{\rm gp}$ is 
$(G_I)_F$-equivariantly homeomorphic to $\overline V_{X_I}^{\rm gp}$. 

\subsubsection{}
The connection is given by the map: 

\vspace{2mmplus1mmminus1mm} 

\centerline{
$\begin{array}{rrrr} 
\hfill\varphi_I:&\hfill\mathscr{S}\bigl((G_I)_F\bigr)\hfill & \hfill\to\hfill & \hfill\mathscr{S}\bigl((P_I)_F\bigr)\subset\mathscr{S}(G_F)\hfill\\ 
&\hfill H\hfill & \hfill\mapsto\hfill & \hfill (H\cdot(T^I)_{\rm cpt})\ltimes(U^I)_F,\hfill 
\end{array}$}

\vspace{2mmplus1mmminus1mm} 

where $\mathscr{S}\bigl((P_I)_F\bigr)$ is the set of closed subgroups in $(P_I)_F$. 
Note that $\varphi_I$ goes from the natural ambient compact space of the compactification 
$\overline V_{X_I}^{\rm gp}$ to the natural ambient compact space of the compactification 
$\overline V_X^{\rm gp}$. 

\begin{lemma}
\label{lemma - between Chabauty} 
The map $\varphi_I$ is continuous and $(G_I)_F$-equivariant. 
It is injective on the subset of closed subgroups of $(G_I)_F$ containing $Z(G_I)_F$, therefore it induces a homeomorphism 
from the latter space onto its image in $\mathscr{S}(G_F)$. 
\end{lemma} 

{\it Proof}.~
Injectivity. 
Let $q_I$ be the natural quotient map $P_I \to P_I/\mathscr{R}(P_I)$. 
Note that we have: ${\rm Ker}(q_I\!\mid_{G_F})=Z(G_I)_F$. 
For $H\!\in\!\mathscr{S}\bigl((G_I)_F\bigr)$, we have: $(q_I\circ\varphi_I)(H)=q_I\bigl(H\cdot(T^I)_{\rm cpt}\cdot(U^I)_F\bigr)=q_I(H)$. 
If $H$, $H'\!\in\!\mathscr{S}\bigl((G_I)_F\bigr)$ are such that $\varphi_I(H)=\varphi_I(H')$, then $q_I(H)=q_I(H')$, so we finally have: 
$H\cdot Z(G_I)_F=(q_I\!\mid_{G_F})^{-1}\bigl(q_I(H)\bigr)=(q_I\!\mid_{G_F})^{-1}\bigl(q_I(H')\bigr)=H'\cdot Z(G_I)_F$. 

\vspace{2mmplus1mmminus1mm} 

Continuity. 
Let us assume that we have: $\displaystyle \lim_{n\to\infty}H_n=H$ in $\mathscr{S}\bigl((G_I)_F\bigr)$. 
Let $g\!\in\!\varphi_I(H)$. 
Then we can write: $g=htu$, with $h\!\in\!H$, $t\!\in\!(T^I)_{\rm cpt}$ and $u\!\in\!(U^I)_F$. 
Since $\displaystyle \lim_{n\to\infty}H_n=H$, we can write: $\displaystyle h=\lim_{n\to\infty}h_n$ with $h_n\!\in\! H_n$ for each $n\geq 1$. 
This enables to write: $\displaystyle g=\lim_{n\to\infty}g_n$ with $g_n\!\in\! \varphi_I(H_n)$ for each $n\geq 1$. 
Now let $\{ n_j \}_{j\geq 1}$ be an increasing sequence of integers, and let $g_{n_j}\!\in\! \varphi_I(H_{n_j})$ converge in $G_F$ to some element $g$. 
We have to show that $g\!\in\! \varphi_I(H)$. 
We can write: $g_{n_j}=h_{n_j}.t_{n_j}.u_{n_j}$ with $h_{n_j}\!\in\!H_{n_j}$, $t_{n_j}\!\in\!(T^I)_{\rm cpt}$ and $u_{n_j}\!\in\!(U^I)_F$. 
Since for every $j\geq 1$ we have: $\varphi_I(H_{n_j})<(P_I)_F=(M_I)_F\ltimes (U^I)_F$, the convergence of $\{g_{n_j} \}_{j\geq 1}$ implies that of $\{u_{n_j} \}_{j\geq 1}$ in $(U^I)_F$, and (up to extracting) that of $\{t_{n_j} \}_{j\geq 1}$ in $(T^I)_{\rm cpt}$. 
This implies the convergence of $\{h_{n_j} \}_{j\geq 1}$. 
Since $\displaystyle \lim_{n\to\infty}H_n=H$, the limit belongs to $H$, and  we have: $g\!\in\! \varphi_I(H)$. 
This proves $\displaystyle \lim_{n\to\infty}\varphi_I(H_n)=\varphi_I(H)$, hence the continuity of $\varphi_I$. 

\vspace{2mmplus1mmminus1mm} 

Equivariance. 
Let $g\!\in\!(G_I)_F$ and let $H$ be a closed subgroup of $(G_I)_F$. 
Since $g$ normalizes $U^I$ and centralizes $T^I$, we have: 

\vspace{2mmplus1mmminus1mm} 

\centerline{
$\varphi_I(gHg^{-1})
=gHg^{-1} \cdot g(T^I)_{\rm cpt}g^{-1} \cdot g(U^I)_Fg^{-1}
=g\bigl(H\cdot (T^I)_{\rm cpt}\cdot (U^I)_F\bigr)g^{-1}
=g\varphi_I(H)g^{-1}$,} 

\vspace{2mmplus1mmminus1mm} 

which proves the desired equivariance. 
\qed

\subsubsection{}
\label{sss - BT in boundary}
We can restrict the previous lemma to sets of closed subgroups corresponding to compactifications. 
When ${\rm rk}_F(G)>1$, this shows that the group-theoretic compactification of $X$ is different from the geometric compactification of the Bruhat-Tits building $X$. 
In the latter compactification, the asymptotic boundary $\partial_\infty X$ is a geometric realization of the spherical building of 
$G_{/F}$; this spherical building reflects the combinatorics of the Tits system of the parabolic $F$-subgroups of $G$. 
In the former one, the boundary contains infinitely many Euclidean buildings: 

\begin{theorem}
\label{th - BT in boundary}
Let $I$ be a proper subset of the set of simple roots $S$ of $G$.  
\begin{enumerate} 
\item[(i)] The map $\varphi_I$ restricts to a $(G_I)_F$-equivariant homeomorphism from the group-theoretic
compactification $\overline V_{X_I}^{\rm gp}$ onto the closure $\overline{Y_I}^{\rm gp}$ of $Y_I$ in the boundary 
of $\overline V_X^{\rm gp}$. 
\item[(ii)] For any proper parabolic $F$-subgroup $Q$ in $G$, the group-theoretic compactification of the Bruhat-Tits building of the reductive $F$-group $Q/\mathscr{R}_u(Q)$ naturally sits in the boundary of the compactification $\overline V_X^{\rm gp}$. 
\item[(iii)] For any group $D\!\in\!\overline{V}^{\rm gp}_X$ there is a sequence 
$\{g_n\}_{n\geq 1}$ in $G_F$ such that 
$\displaystyle \lim_{n\to\infty}g_nDg_n{}^{-1}$ exists and lies in the closed orbit 
$\mathscr{F}=\{gD_\varnothing g^{-1}\}_{g\in G_F}$. 
\item[(iv)] The maximal Furstenberg boundary $\mathscr{F}\simeq G_F/P_F$ is the only closed $G_F$-orbit in $\overline{V}^{\rm gp}_X$. 
\end{enumerate}
\end{theorem} 

{\it Proof}.~
(i).~Let $J$ be a subset of $I$ and $\underline d$ be a family of non-negative real numbers indexed by $J$. 
We have: 
$\varphi_I(D_{J\subset I,\underline d})=\bigl((K_{J,\underline d}\cap G_I)\ltimes(U^J\cap G_I)_F\bigr) 
\cdot\bigl((T^I)_{\rm cpt}\ltimes(U^I)_F\bigr)$. 
Since $T^I$ centralizes $G_I$, this group is also 
$\bigl((K_{J,\underline d}\cap G_I)\cdot(T^I)_{\rm cpt}\bigr)\ltimes\bigl((U^J\cap G_I)_F\cdot(U^I)_F\bigr)$. 
Since $(K_{J,\underline d}\cap G_I)\cdot(T^I)_{\rm cpt}=K_{J,\underline d}$ and 
$(U^J\cap G_I)\cdot(U^I)=U^J$, this finally proves that $\varphi_I(D_{J\subset I,\underline d})$ is equal to the 
maximal limit group $D_{J,\underline d}$ of $G_F$. 
This shows that $\varphi_I(\overline V_{X_I}^{\rm gp})$ is contained in $\overline V_X^{\rm gp}\cap\mathscr{S}\bigl((P_I)_F\bigr)$, and that 
the preimage $\varphi_I^{-1}(D_{J,\underline d})$ is compact in $(G_I)_F$ if, and only if, $J=I$; in which case it is a maximal compact subgroup of $(G_I)_F$. 

\vspace{2mmplus1mmminus1mm} 

By definition, $Y_I$ is the union of the $(G_I)_F$-conjugacy classes of the maximal limit groups 
$D_{I,\underline d_0}$, $D_{I,\underline d_1}$, ... $D_{I,\underline d_{\mid I \mid}}$ (\ref{sss - definitions in the boundary}), and 
the groups $\varphi_I^{-1}(D_{I,\underline d_0})$, $\varphi_I^{-1}(D_{I,\underline d_1})$, ... $\varphi_I^{-1}(D_{I,\underline d_{\mid I \mid}})$ 
provide a complete system of representatives for the $(G_I)_F$-conjugacy classes of maximal compact subgroups in $(G_I)_F$. 
We have: $Y_I \subset \varphi_I(\overline V_{X_I}^{\rm gp})$, and by compactness of $\overline V_{X_I}^{\rm gp}$ and continuity of $\varphi_I$ (Lemma \ref{lemma - between Chabauty}), it follows that $\overline{Y_I}^{\rm gp}\subset\varphi_I(\overline V_{X_I}^{\rm gp})$. 
Therefore, in order to prove (i) it remains to prove the converse inclusion. 
Let $D\!\in\!\overline V_{X_I}^{\rm gp}$. 
By definition of a group-theoretic compactification, we have: 
$\displaystyle D=\lim_{n\to\infty}K'_n$, where $K'_n$ is a maximal compact subgroup of $(G_I)_F$. 
By the first remark of the paragraph, for each $n\geq1$ we can write: $K'_n=g_n\varphi_I^{-1}(D_{I,\underline d_{i(n)}})g_n^{-1}$ 
with $g_n\!\in\!(G_I)_F$ and $i(n)\!\in\!\{0;1;...\mid\!I\!\mid\}$. 
By equivariance and continuity of $\varphi_I$, we have: 
$\displaystyle \varphi_I(D)=\lim_{n\to\infty} g_nD_{I,\underline d_{i(n)}}g_n^{-1}$, which proves the desired inclusion. 

\vspace{2mmplus1mmminus1mm} 

(ii).~ Since the Bruhat-Tits building of a semisimple $F$-group is the building of its simply connected covering $F$-group, (ii) follows from (i) by conjugating by a suitable element in $K_o$. 

\vspace{2mmplus1mmminus1mm} 

(iii).~For a suitable element $k\!\in\!K_o$, the conjugate $kDk^{-1}$ lies in the group-theoretic compactification of the apartment $A$ attached to $T$. 
Let $t\!\in\!T_F$ be a regular element such that the vertex $t.o$ lies in the Weyl chamber $\mathscr{Q}$. 
By Theorem \ref{th - CV groups} we have: $\displaystyle \lim_{n\to\infty}t^nkDk^{-1}t^{-n}=D_\varnothing$, so 
we can take $g_n=t^nk$. 

\vspace{2mmplus1mmminus1mm} 

(iv) follows directly from Lemma \ref{lemma - stabilizers} and the previous paragraph. 
\qed\vspace{2mmplus1mmminus1mm} 

Let $g\!\in\!(G_I)_F$ and $u\!\in\!(U^I)_F$. 
We have: $u.g(U^I)_F.u^{-1}=g.(g^{-1}ug).(U^I)_F.u^{-1}$. 
But since $U^I$ is normalized by $G_I$, this implies that $u$ stabilizes the class $g(U^I)_F$. 
It follows that we have: 

\vspace{2mmplus1mmminus1mm} 

\centerline{$\displaystyle u.\bigl(gK_{I,\underline d}g^{-1} \ltimes (U^I)_F\bigr).u^{-1}
=u.\bigl(\bigcup_{h\in gK_{I,\underline d}g^{-1}} h.(U^I)_F\bigr).u^{-1}
=\bigcup_{h\in gK_{I,\underline d}g^{-1}} h.(U^I)_F
=gK_{I,\underline d}g^{-1} \ltimes (U^I)_F$.}Ê

\vspace{2mmplus1mmminus1mm} 

This proves that the action by conjugation of the unipotent radical $(U^I)_F$ on the limit groups of $G_F$ contained in 
$\varphi_I(\overline V_{X_I}^{\rm gp})$ is trivial. 
Even more simply, the fact that $T^I$ and $G_I$ commute with one another implies that the action by conjugation of the torus $(T^I)_F$ on the same subgroups is trivial. 
To sum up, we have: 

\begin{lemma}Ê
\label{lemma - type-preserving}Ê
The group $(P_I)_F$ acts on $\varphi_I(\overline V_{X_I}^{\rm gp})$ via the projection map $P_I\to G_I$. 
\qed
\end{lemma}Ê

The lemma implies in particular that the $(P_I)_F$-action on $\varphi_I(\overline V_{X_I}^{\rm gp})$ is type-preserving. 

\subsection{Polyhedral compactification}
\label{ss - polyhedral}
In this subsection, we use a compactification procedure defined in terms of the very definition of Bruhat-Tits buildings, 
i.e. by gluings. 
This construction is analogous to the case of symmetric spaces. 
We show in Theorem \ref{th - identification polyhedral} that, loosely speaking, the gluing procedure for the polyhedral compactification eventually amounts to filling in the vertices left empty by the group-theoretic procedure. 
Conversely, the Chabauty topology viewpoint provides a concrete approach to the polyhedral compactification. 

\subsubsection{}
\label{sss - summary polyhedral} 
The polyhedral compactification of Bruhat-Tits buildings is defined in \cite[\S 14]{Landvogt}. 

\vspace{2mmplus1mmminus1mm} 

The first step of the construction consists in compactifying an apartment, say $A$, by replacing Weyl chambers by {\it corners}~\cite[\S 2]{Landvogt}. 
This is done by using the combinatorics of the root system $\Phi$ associated to the maximal $F$-split torus $T$ defining  $A$ \cite[2.9]{Landvogt}. 
We call the so-obtained compactification the {\it polyhedral compactification~}of $A$, and we denote it by $\overline{A}^{\rm pol}$. 
In order to describe it topologically, let us denote by $\overline{\mathscr{Q}}^{\rm pol}$ the closure of the Weyl chamber $\mathscr{Q}$ in $\overline{A}^{\rm pol}$, and let $S$ be the set of simple roots in $\Phi$ defined by $\mathscr{Q}$. 
The points of $\overline{\mathscr{Q}}^{\rm pol}$ are in bijection with the families of parameters $\underline d=\{d_s\}_{s\in S}$ indexed by $S$ with values in $[0;+\infty]$. 
The topology on $[0;+\infty]$ extends the natural one on $[0;+\infty)$ by taking the intervals $[t,+\infty)$, $t\geq0$, as a basis of neighborhoods of $+\infty$. 
To $\underline d$ is attached the subset $I(\underline d)$ of $S$ by setting: $s\!\in\! I(\underline d) \Leftrightarrow d_s<+\infty$. 
For each subset $I$ of $S$, we set: $\mathscr{Q}_I=\{\underline d : I(\underline d)=I\}$. 
Set-theoretically, we have: $\overline{\mathscr{Q}}^{\rm pol}=\bigsqcup_{I\subset S} \mathscr{Q}_I$. 
The families $\underline d$ all of whose parameters are real numbers, i.e. those for which $I(\underline d)=S$, parametrize the points of $\overline{\mathscr{Q}}^X$. 
In this case the parameter $d_s$ corresponds to the distance to the face $\Pi^s$ of $\overline{\mathscr{Q}}^X$, where $\partial_\infty \Pi^s$ is the panel of type $s$ of $\partial_\infty \mathscr{Q}$ in the spherical 
building at infinity $\partial_\infty X$. 
A sequence $\{\underline d_n\}_{n\geq 1}$ converges to $\underline d$ if, and only if, for each $s\!\in\!S$ we have: $\displaystyle \lim_{n\to\infty}  d_{s,n}=d_s$ in $[0;+\infty]$. 
For instance, any sequence $\{\underline d_n\}_{n\geq 1}$ such that $\displaystyle \lim_{n\to\infty}  d_{s,n}=+\infty$ for each $s\!\in\!S$, converges to the unique point of $\mathscr{Q}_\varnothing$. 
Note that in the case of a sequence of vertices in $\overline{\mathscr{Q}}^X$, we obtain the same picture as the one given by convergence of vertices of $\overline{\mathscr{Q}}^X$ in the group-theoretic compactification $\overline{V}^{\rm gp}_X$ (Theorem \ref{th - CV groups}). 

\vspace{2mmplus1mmminus1mm} 

In the second step of the compactification, one attaches two groups $P_x$ and $U_x$ to each point 
$x\!\in\!\overline{A}^{\rm pol}$. 
Let us rephrase geometrically the definitions of $U_x$ and $P_x$ when $x\!\in\!\overline{\mathscr{Q}}^{\rm pol}$. 
We denote by $\underline d$ the parameters associated to $x$.  

\vspace{2mmplus1mmminus1mm} 

{\bf Definition.}~\it 
For each root $a\!\in\!\Phi$ and each $x\!\in\!\overline{A}^{\rm pol}$, we denote by $U_{a,x}$ the biggest subgroup of $(U_a)_F$ fixing $x$. 
\rm\vspace{2mmplus1mmminus1mm} 

Geometrically, the fixed-point set of $U_{a,x}$ in $A$ is the smallest half-apartement bounded by a wall of direction $\partial a$ whose closure in $\overline{A}^{\rm pol}$ contains $x$. 
If $a\!\in\!\Phi_{I(\underline d)}$ then $U_{a,x}$ is a non-trivial proper subgroup of the filtration of $(U_a)_F$ given by the valuation of the root datum associated to the vertex $o$. 
If $a\not\in\Phi_{I(\underline d)}$ then $U_{a,x}=\{1\}$ when $a$ is a negative root, and $U_{a,x}=(U_a)_F$ when $a$ is positive. 

\vspace{2mmplus1mmminus1mm} 

We  further denote by $N_x$ the fixator of $x$ in ${\rm Stab}_{G_F}(A)=N_G(T)_F$. 
If $x\!\in\!A$ then $N_x$ is a compact extension of the finite group $W_x$ by ${\rm Fix}_{G_F}(A)=Z_G(T)_{\rm cpt}$, where $W_x$ is the stabilizer of $x$ in the affine Weyl group action on $A$. 
If $x\!\in\!\overline{A}^{\rm pol}\setminus A$, then $I(\underline d)\subsetneq S$ and the finite parameters in $\underline d$ determine a unique point, say $x_I$, in the Weyl chamber $p_I(\mathscr{Q})$ in the standard apartment of the Euclidean building $X_I$ of $(G_I)_{/F}$. 
In this case, the group $N_x$ is the commutative product of the compact fixator of $x_I$ in  ${\rm Stab}_{(G_I)_F}\bigl(p_I(A)\bigr)$ by the torus $(T^{I(\underline d)})_F$. 
In other words, we reduce the case $x\!\in\!\overline{A}^{\rm pol}\setminus A$ to the case $x\!\in\!A$ by passing to a Levi factor, and then add a non compact factor corresponding translations orthogonal to the apartment of the Levi factor. 
The definitions below are taken from \cite[12.4]{Landvogt}. 

\vspace{2mmplus1mmminus1mm} 

{\bf Definition.}~\it 
The group $U_x$ is defined to be the group generated by the groups $U_{a,x}$ when $a$ varies in $\Phi$, and the group $P_x$ is defined to be the group generated by $U_x$ and $N_x$. 
\rm\vspace{2mmplus1mmminus1mm} 

For an arbitrary $x$, the group $P_x$ turns out to be the stabilizer of $x$ in $G_F$ once $\overline{X}^{\rm pol}$ is defined \cite[14.4 (i)]{Landvogt}. 
If $x\!\in\!A$ then $P_x$ is equal to the parahoric subgroup $K_x$. 
Otherwise, since the normalizer of $T_F$, i.e. the stabilizer of $A$, is transitive on the Weyl chambers with tip $o$, we may -- and shall -- assume that $x\!\in\!\overline{\mathscr{Q}}^{\rm pol}\setminus\mathscr{Q}$, and the situation is described as follows. 

\begin{lemma}
\label{lemma - devissage stabilizers}
Let $x\!\in\!\overline{\mathscr{Q}}^{\rm pol}\setminus\mathscr{Q}$ and let $\underline d$ be the parameters corresponding to $x$, with $I\subsetneq S$. 
Then: 
\begin{enumerate}
\item [(i)] the limit group $D_{I(\underline d),\underline d}$ is equal to the group $U_x \cdot Z_G(T)_{\rm cpt}$; 
\item [(ii)] its normalizer $R_{I(\underline d),\underline d}$ is equal to the group $P_x$. 
\end{enumerate}
\end{lemma} 

In this lemma the understatement is that in the index ${}_{I(\underline d),\underline d}$ we see $\underline d$ as a family of non-negative real numbers indexed by $I(\underline d)$, i.e. we forget the infinite parameters of $\underline d$. 

\vspace{2mmplus1mmminus1mm} 

{\it Proof}.~
(ii). By \cite[12.6]{Landvogt} we can write: $P_x=U^-_xU^+_xN_x$, and by \cite[12.5 (ii)]{Landvogt} we can write: $U^+_x=U^+_{I(\underline d),x}\cdot U^{I(\underline d),+}_x$, where $U^+_{I(\underline d),x}$ (resp. $U^{I(\underline d),+}_x$) is generated by the groups $U_{a,x}$ with $a\!\in\!\Phi_{I(\underline d)}^+$ (resp. $a\!\in\!\Phi^{I(\underline d),+}$). 
It follows from the discussion after the definition of $U_{a,x}$ that $U^{I(\underline d),+}$ is the full unipotent radical of the parabolic subgroup $(P_I)_F$. 
Morevover we have: $U^-_x=U^-_{I(\underline d),x}$ where $U^-_{I(\underline d),x}$ is generated by the groups $U_{a,x}$ with $a\!\in\!\Phi_{I(\underline d)}^-$. 
This implies that $U^-_{I(\underline d),x}U^+_{I(\underline d),x}N_x$ lies in the reductive Levi factor $(M_I)_F$ so we can write 
$P_x=(U^-_{I(\underline d),x}U^+_{I(\underline d),x}N_x)\cdot (U^{I(\underline d)})_F$. 
But it follows from the description of parahoric subgroups in $(G_I)_F$ \cite[6.4.9]{BrT1} and the description of $N_x$, that 
$U^-_{I(\underline d),x}U^+_{I(\underline d),x}N_x$ is a group, namely the product of $(T^{I(\underline d)})_F$ and of the parahoric  subgroup $G_I\cap K_{I(\underline d), \underline d}$ of $(G_I)_F$. 
This finally implies: $P_x=(K_{I(\underline d), \underline d}\cdot (T^{I(\underline d)})_F)\ltimes(U^{I(\underline d)})_F$. 

\vspace{2mmplus1mmminus1mm} 

(i). By \cite[12.5 (iii)]{Landvogt}, we have: $U_x=U^-_xU^+_x(N_x\cap U_x)$. 
Arguing as for (ii), we obtain: 
$U_x=\bigl(U^-_{I(\underline d),x}U^+_{I(\underline d),x}(N_x\cap U_x)\bigr)\cdot (U^{I(\underline d)})_F$. 
Multiplying by the pointwise fixator of $A$ provides: $U_x \cdot Z_G(T)_{\rm cpt}
=\bigl(U^-_{I(\underline d),x}U^+_{I(\underline d),x}(N_x\cap U_x)Z_G(T)_{\rm cpt}\bigr)\cdot (U^{I(\underline d)})_F$. 
By the description of $(N_x\cap U_x)$ given in \cite[12.5 (iv)]{Landvogt}, it follows that the group $U^-_{I(\underline d),x}U^+_{I(\underline d),x}(N_x\cap U_x)Z_G(T)_{\rm cpt}$ is the parahoric sugbroup $K_{I(\underline d), \underline d}$ of the reductive Levi factor $(M_I)_F$, so finally: $U_x \cdot Z_G(T)_{\rm cpt}=K_{I(\underline d), \underline d}\ltimes(U^{I(\underline d)})_F$. 
\qed\vspace{2mmplus1mmminus1mm} 

The last step to define the polyhedral compactification consists in extending the equivalence relation (R) used in \cite[7.4.2]{BrT1} to define $X$. 
The polyhedral compactification $\overline{X}^{\rm pol}$ is 
$\displaystyle {G_F \times \overline{A}^{\rm pol} \over \sim^*}$, where $\sim^*$ is a suitable extension of (R) involving the groups $U_x$ \cite[14.2]{Landvogt} and where $x$ runs over $\overline{A}^{\rm pol}$. 

\subsubsection{}
\label{sss - Landvogt is OK} 
As pointed out to us by A. Werner, in order to use the polyhedral compactification $\overline{X}^{\rm pol}$, we need 
to fix the proof of \cite[14.11]{Landvogt}. 
In terms of Bruhat cells, the mistake amounts to saying that a sequence of points of the big cell 
$\Omega_F=(U^-\cdot Z_G(T) \cdot U^+)_F$, converging in $G_F$, has its limit in $\Omega_F$ (while $\Omega_F$ is open and dense in $G_F$). 
Here is the statement. 

\begin{prop}
\label{prop - Landvogt is OK}
Let $\{x_n\}_{n\geq 1}$ be a sequence in the compactification $\overline{A}^{\rm pol}$ of the 
apartment $A$. 
Let $\{g_n\}_{n\geq 1}$ be a sequence in $G_F$ with $g_n\!\in\!P_{x_n}$ for each index $n\geq 1$. 
We assume that $\{g_n\}_{n\geq 1}$ converges to some $g\!\in\!G_F$ 
and that $\{x_n\}_{n\geq 1}$ converges to some $x\!\in\!\overline{A}^{\rm pol}$. 
Then, we have: $g\!\in\!P_x$. 
\end{prop} 

The proof uses the following reformulation of the proposition: 
let $\{x_n\}_{n\geq 1}$ be a sequence of $\overline{A}^{\rm pol}$ converging to $x$; 
then $P_x$ contains any Chabauty cluster value of $\{P_{x_n}\}_{n\geq 1}$. 

\vspace{2mmplus1mmminus1mm} 

{\it Proof}.~
Recall that we have a finite disjoint union decomposition: 
$\overline{\mathscr{Q}}^{\rm pol}=\bigsqcup_{I\subset S} \mathscr{Q}_I$ (\ref{sss - summary polyhedral}), so there is a subset $J$ of $S$ such that $x\!\in\!\mathscr{Q}_J$; and up to extracting a subsequence, we may -- and shall -- assume that for each $n\geq 1$, we have: $x_n\!\in\!\mathscr{Q}_I$ for some subset $I$ of $S$. 
In view of the topology on $\overline{\mathscr{Q}}^{\rm pol}$, we have: $J\subset I$. 
Let us denote by $\underline d_n$ (resp. by $\underline d$) the parameters of the point $x_n$ (resp. of $x$). 
We have: $I(\underline d_n)=I$ for each $n\geq 1$. 
By Lemma \ref{lemma - devissage stabilizers}, the sequence $\{P_{x_n}\}_{n\geq 1}$ is nothing else than the sequence 
$\{R_{I,\underline d_n}\}_{n\geq 1}$, and $P_x=R_{J,\underline d}$. 
We see $\mathscr{Q}_I$ as the standard Weyl chamber $p_I(\mathscr{Q})$ in the building $X_I$ of the semisimple Levi factor $G_I$. 

\vspace{2mmplus1mmminus1mm} 

Up to extracting a subsequence, we may -- and shall -- assume that there exists a closed subgroup $R<(P_I)_F$ such that $\{R_{I,\underline d_n}\}_{n\geq 1}$ converges to $R$ for the Chabauty topology on the closed subgroups of $(P_I)_F$. 
We have: $g\!\in\!R$, so it is enough to show that $R<R_{J,\underline d}$. 
Replacing $G$ by its Levi factor $G_I$, we are reduced to a convergence problem similar to that of \ref{ss - CV groups}, except that this time we consider sequences of extensions of parahoric subgroups of $(G_I)_F$ by the radical $\mathscr{R}(P_I)_F$. 
Moreover the sequence $\{x_n\}_{n\geq1}$ is not {\it a priori}~a $J$-canonical sequence in $\mathscr{Q}_I$ but since the facet of $p_J(\mathscr{Q}_I)=\mathscr{Q}_J$ containing $x$ lies in finitely many closures of facets, we are reduced to this case (possibly after extracting  again a subsequence). 
Let $\sigma$ be the facet of $\mathscr{Q}_J$ containing $x$, and let $\tau$ be the facet of $\mathscr{Q}_J$ defined by the $J$-canonical sequence $\{x_n\}_{n\geq1}$ and such that $\overline\tau\supset\sigma$. 
We are only interested in an upper bound for inclusion of $R$, i.e. in the analogue of Lemma \ref{lemma - upper bound cluster value}. 
Since $\mathscr{R}(P_I)_F=(T^I)_F\ltimes(U^I)_F$ acts trivially on the Furstenberg boundary $\mathscr{F}_I=(G_I)_F.\omega$, we can use probability measures on $\mathscr{F}_I$ as in the proof of Proposition \ref{prop - CV probas}. 
If $\underline d'$ is a family of parameters indexed by $J$ and determining a point in $\tau$, then arguing as in the proof of Lemma \ref{lemma - upper bound cluster value}, we obtain: $R<R_{J,\underline d'}$. 
The inclusion $\overline\tau\supset\sigma$ implies: $R_{J,\underline d'}<R_{J,\underline d}$, so we finally obtain: $R<R_{J,\underline d}$. 
\qed

\subsubsection{}
\label{sss - comparison with polyhedral} 
Let us now define a map from the polyhedral compactification to the group-theoretic one. 
Let $x\!\in\!\overline{X}^{\rm pol}$ and let $Q_x$ denote the parabolic $F$-subgroup $\overline{P_x}^Z$ 
(Lemma \ref{lemma - devissage stabilizers}). 
We denote by $X^*(Q_x)_F$ the abelian group of algebraic characters of $Q_x$ defined over $F$. 
For each $\chi\!\in\!X^*(Q_x)_F$, we set: 
$\overline\chi=\bigl(\mid\!-\!\mid_F \circ \, \chi\!\mid_{P_x}\bigr):P_x \to {\bf R}^\times_+$ . 

\vspace{2mmplus1mmminus1mm} 

{\bf Definition.}~\it
We denote by $D_x$ the intersection of all the kernels of the continuous characters $\overline\chi$ of $P_x$ when $\chi$ varies over $X^*(Q_x)_F$, i.e. $D_x=\bigcap_{\xi\in X^*(Q_x)_F} \overline\chi$. 
\rm\vspace{2mmplus1mmminus1mm} 

Assume that $x\!\in\!\overline{\mathscr{Q}}^{\rm pol}$ is defined by $I$ and the parameters $\underline d$, so that  $P_x=R_{I,\underline d}$ and $Q_x=P_I$. 
Then for $g\!\in\!P_I$ and $\chi\!\in\!X^*(P_I)_F$, we have: $\chi(g)=\chi(t)$ where $g=htu$ with $h\!\in\!G_I$, $t\!\in\!T^I$ and 
$u\!\in\!U^I$. 
If we choose $g$ in $R_{I,\underline d}$, i.e. if $h\!\in\!K_{I,\underline d}\cap G_I$, $t\!\in\!(T^I)_F$ and $u\!\in\!(U^I)_F$, then this equality implies that $g$ belongs to $D_x$ if, and only if, $t\!\in\!(T^I)_{\rm cpt}$. 
In other words, when $x\!\in\!\overline{\mathscr{Q}}^{\rm pol}$ the group $D_x$ coincides with $D_{I,\underline d}$. 

\vspace{2mmplus1mmminus1mm} 

{\bf Definition.}~\it
We denote by $\overline{V}^{\rm pol}_X$ the closure of the set of vertices $V_X$ in the polyhedral compactification 
$\overline{X}^{\rm pol}$. 
\rm\vspace{2mmplus1mmminus1mm} 

It follows from the Cartan decomposition that convergence of sequences in a compactification is basically described thanks to sequences in a given Weyl chamber. 
Therefore, in view of the comparison between the Chabauty convergence of parahoric subgroups as described by Theorem \ref{th - CV groups} and the very definition of the polyhedral compactification of an apartment (\ref{sss - summary polyhedral}), the identification below is not surprising. 

\begin{theorem}
\label{th - identification polyhedral}
Let $G$ be a semisimple simply connected group defined over a non-archimedean local field $F$. 
Let $X$ be the corresponding Bruhat-Tits building. 
Then the map $\mathscr{D}:x \mapsto D_x$ establishes a $G_F$-equivariant homeomorphism between the polyhedral compactification $\overline{V}^{\rm pol}_X$ and the group-theoretic compactification $\overline{V}^{\rm gp}_X$. 
\end{theorem} 

{\it Proof}.~
Since the stabilizer map $x\mapsto P_x={\rm Stab}_{G_F}(x)$ is $G_F$-equivariant, the map $\mathscr{D}$ is $G_F$-equivariant by definition of $D_x$ in $P_x$. 
Moreover a maximal limit group $D_{I,\underline d}$ is the image by $\mathscr{D}$ of the point in $\overline{\mathscr{Q}}^{\rm pol}$ defined by the parameters $(I,\underline d)$. 
By equivariance, this implies that $\mathscr{D}$ is surjective. 

\vspace{2mmplus1mmminus1mm} 

Let $x,y\!\in\!\overline{V}^{\rm pol}_X$ be such that $D_x=D_y$.
There exists $g\!\in\!G_F$ such that $x,y\!\in\!g.\overline{A}^{\rm pol}$ \cite[14.7]{Landvogt}, so we are reduced to the case when $x,y\!\in\!\overline{A}^{\rm pol}$ and $D_x=D_y$. 
In fact, using a suitable $n\!\in\!N_G(T)_F$, we may -- and shall -- even assume that $x\!\in\!\overline{\mathscr{Q}}^{\rm pol}$, $y\!\in\!\overline{A}^{\rm pol}$ and $D_x=D_y$. 
Let us denote by $I$ the subset of $S$ such that $x\!\in\!\mathscr{Q}_I$, and let $w$ be the element of the spherical Weyl group such that $y\!\in\!w.\overline{\mathscr{Q}}^{\rm pol}$. 
From $x\!\in\!\mathscr{Q}_I$, we deduce that $\overline{D_x}^Z$ is the standard parabolic subgroup $P_I$, and from $y\!\in\!w.\overline{\mathscr{Q}}^{\rm pol}$ we deduce that $\overline{D_y}^Z$ contains the minimal parabolic subgroup $wPw^{-1}$. 
The equality $\overline{D_x}^Z=\overline{D_y}^Z$ finally implies that $\overline{D_y}^Z=wP_Iw^{-1}$ and  that $w$ lies in the Weyl group $W_I$ of $P_I$. 
In particular, we have: $y\!\in\!W_I.\mathscr{Q}_I$, so $y$ lies in $p_I(A)$. 
By the description of the limit groups it follows that $x$ and $y$ are contained in the same facet of the apartment $p_I(A)$ in the building $X_I$. 
If $x$ and $y$ are both in $\overline{V}^{\rm pol}_X$, this implies $x=y$, which proves the injectivity of $\mathscr{D}$. 

\vspace{2mmplus1mmminus1mm} 

At this stage, we know that $\mathscr{D}$ is a $G_F$-equivariant bijection 
$\overline{V}^{\rm pol}_X\simeq \overline{V}^{\rm gp}_X$. 
By compactness, it remains to show that $\mathscr{D}$ is continuous. 
Let $\{ x_n \}_{n\geq 1}$ be a sequence in $\overline{V}^{\rm pol}_X$ converging to some point $x$.  
We can write $x_n=k_n.q_n$ with $k_n\!\in\!K_o$ and $q_n\!\in\!\overline{\mathscr{Q}}^{\rm pol}$ for each $n\geq 1$ \cite[14.18]{Landvogt}. 
We have to show that $D_x$ is the only Chabauty cluster value of $\{ D_{x_n} \}_{n\geq 1}$. 
Let $\displaystyle D=\lim_{j\to\infty} D_{x_{n_j}}$ be such a cluster value. 
We may -- and shall -- assume that $\{k_{n_j}\}_{j\geq1}$ and $\{q_{n_j}\}_{j\geq1}$ are such that 
$\displaystyle \lim_{j\to\infty}k_{n_j}=k$ in $K_o$ and $\displaystyle \lim_{j\to\infty}q_{n_j}=q$ in $\overline{\mathscr{Q}}^{\rm pol}$. 
By continuity of the $G_F$-action on $\overline{X}^{\rm pol}$ \cite[14.31]{Landvogt}, we have: $x=k.q$, and the equivariance of $\mathscr{D}$ implies that $D=kD_qk^{-1}$ and $D_{x_{n_j}}=k_{n_j}D_{q_{n_j}}k_{n_j}^{-1}$. 
By continuity of the $G_F$-action by conjugation on $\mathscr{S}(G_F)$, it is enough to show that $\{ D_{q_{n_j}} \}_{j\geq 1}$ converges to $D_q$ in the Chabauty topology. 

\vspace{2mmplus1mmminus1mm} 

Let us denote by $\underline d$ the parameters of $q$ and by $\underline d_{n_j}$ the parameters of $q_{n_j}$. 
By the finite disjoint union decomposition: 
$\overline{\mathscr{Q}}^{\rm pol}=\bigsqcup_{I\subset S} \mathscr{Q}_I$ (\ref{sss - summary polyhedral}), there exists $J\subset S$ such that $x\!\in\!\mathscr{Q}_J$. 
Up to extracting a subsequence, we may -- and shall -- assume that for each $j\geq 1$, we have: $q_{n_j}\!\in\!\mathscr{Q}_I$ for some $I\subset S$ such that $J\subset I$. 
We have: $D_q=D_{J,\underline d}$, and $D_{q_{n_j}}=D_{I,\underline d_{n_j}}$ for any index $j\geq1$. 
In particular, we have: $D_{q_{n_j}}<R_{I,\underline d_{n_j}}={\rm Stab}_{G_F}(q_{n_j})$, so applying Proposition \ref{prop - Landvogt is OK} we obtain an upper for $D$ with respect to inclusion of closed subgroups: $D<R_{J,\underline d}$. 
Then arguing geometrically as in Lemma \ref{lemma - lower bound cluster value} shows that we also have a lower bound: 
$D_{J,\underline d}<D$. 
Since $R_{J,\underline d}=D_{J,\underline d}\cdot (T^J)_F$, it remains to show that $D \cap T^J<K_{J,\underline d}$. 
First we apply Lemma \ref{lemma - limit is distal} to see that (in a linear $F$-embedding $G<{\rm GL}_m$) any element in each limit group $D_{I,\underline d_{n_j}}$ has all its eigenvalues of absolute value equal to 1. 
Then the proof of Lemma \ref{lemma - limit is distal} itself shows that the same property holds for any element of the group $D$. 
Arguing as in the final part in the proof of Theorem \ref{th - CV groups}, we conclude that $T^J\cap D<K_{J,\underline d}$. 
This finally proves that $D_q$ is the only cluster value of $\{D_{q_{n_j}}\}_{j\geq1}$. 
\qed\vspace{2mmplus1mmminus1mm} 

It was already known for the polyhedral compactification that the Bruhat-Tits buildings of the semisimple quotients of the proper parabolic $F$-subgroups of $G$ naturally sit in the boundary of $\overline{X}^{\rm pol}$ \cite[14.21]{Landvogt}. 
The combination of Theorem \ref{th - identification polyhedral} and of Theorem \ref{th - BT in boundary} enables to prove it in a natural way. 
The limit groups contained in a given parabolic subgroup $Q_{/F}$, once divided out by the unipotent radical $\mathscr{R}_u(Q)_F$, are the closed subgroups appearing in the group-theoretic compactification of the Euclidean building of $Q/\mathscr{R}_u(Q)$. 
Checking that both topologies coincide (the one from the big ambient compactification $\overline{V}^{\rm gp}_X$ and the one from the Chabauty topology on the closed subgroups of $Q_F$) amounts to computing convergence of sequences of limit groups by the same techniques as those used to compute convergence of sequences of maximal compact subgroups. 

\vspace{2mmplus1mmminus1mm}

{\bf Definition.}~\it 
Let $X$ be the Bruhat-Tits building of a simply connected semisimple $F$-group $G$.
\begin{enumerate}
\item [(i)] We call {\rm Euclidean}~or {\rm Bruhat-Tits building at infinity}~of $\overline{X}^{\rm pol}$, or of $\overline{V}^{\rm gp}_X$, the building of the semisimple quotient of some proper parabolic $F$-subgroup of $G$, embedded as in the previous subsection. 
\item [(ii)] We call {\rm facet}~of $\overline{X}^{\rm pol}$, or of $\overline{V}^{\rm gp}_X$, any facet in $X$ or in a Euclidean building at infinity of the compactification under consideration. 
\end{enumerate}
\rm\vspace{2mmplus1mmminus1mm} 

There is a criterion in terms of stabilizers to decide whether two points lie in the same Euclidean building, or in the same facet, of $\overline{X}^{\rm pol}$. 

\begin{prop}
\label{prop - commensurable}
Let $x$ and $y$ be points in $\overline{X}^{\rm pol}$. 
\begin{enumerate}
\item [(i)] The groups $P_x$ and $P_y$ are commensurable if, and only if, they have the same Zariski closure. 
Geometrically, this amounts for $x$ and $y$ to being in the same Bruhat-Tits building. 
\item [(ii)] The same holds with $P_x$ and $P_y$ replaced by $D_x$ and $D_y$, respectively. 
\item[(iii)] We have: $P_x=P_y$ if, and only if, $D_x=D_y$; geometrically, this amounts for $x$ and $y$ to being in the same facet of $\overline{X}^{\rm pol}$. 
\end{enumerate}
\end{prop} 

{\it Proof}.~(i). 
Combining Lemma \ref{lemma - devissage stabilizers} and Lemma \ref{lemma - Z density}, we know that the Zariski closures $\overline{P_x}^Z$ and $\overline{P_y}^Z$ are parabolic $F$-subgroups, hence are connected. 
If $P_x$ and $P_y$ are commensurable, then so are $\overline{P_x}^Z$ and $\overline{P_y}^Z$, and connectedness for the Zariski topology implies: $\overline{P_x}^Z=\overline{P_y}^Z$. 
Conversely, if $\overline{P_x}^Z=\overline{P_y}^Z$, then we denote by $Q$ this Zariski closure, and by $p:Q\to Q/\mathscr{R}(Q)$ the quotient map dividing by the radical. 
The groups $p(P_x)$ and $p(P_y)$ are open and compact in the same topological group, so they are commensurable, and it remains to note that $P_x=p^{-1}\bigl(p(P_x))$ and $P_y=p^{-1}\bigl(p(P_y))$. 
This proves (i), and (ii) is proved similarly. 

\vspace{2mmplus1mmminus1mm} 

(iii). The first equivalence follows from the equalities: $P_x=N_{G_F}(D_x)$ and $D_x={\rm Ker}(\Delta_{P_x})$. 
One implication of the remaining equivalence is clear, while the other is proved in the second paragraph of the proof of 
Theorem \ref{th - identification polyhedral} (injectivity of $\mathscr{D}$). 
\qed

\section{The case of trees}
\label{s - trees}
We define the group-theoretic compactification of a locally finite tree $X$ using the Chabauty topology on the closed subgroups of a sufficiently transitive automorphism group $G$. 
Trees provide an example for what was done in the previous sections, but results in this section also settle the initial step of induction arguments used in the next section. 
Note that here, the tree $X$ is only assumed to admit a big (i.e. locally $\infty$-transitive) automorphism group. 
In particular, it needn't come from a rank-one algebraic group over a local field (e.g. it may be a homogeneous tree of valency 7). 

\subsection{Combinatorics of locally $\infty$-transitive groups}
\label{ss - tree combinatorics}
The full automorphism group ${\rm Aut}(X)$ is endowed
with the topology of uniform convergence on finite subsets \cite[I.4]{FTN}. 
A basis of this topology consists of the subsets $U_Y(g) \subset {\rm Aut}(X)$, where: 

\vspace{2mmplus1mmminus1mm} 

\centerline{$Y \subset X$ is finite, \quad $g \!\in\! {\rm Aut}(X)$ \quad and \quad 
$U_Y(g)=\{h\!\in\!{\rm Aut}(X):g\!\mid_Y\,=\,h\!\mid_Y\}$.}

\vspace{2mmplus1mmminus1mm} 

The topology is locally compact and totally disconnected. 

\subsubsection{}
Rather than studying the only group ${\rm Aut}(X)$, we consider a wider class of closed non-discrete subgroups defined in \cite{BurMozTrees}. 
We denote by $S(v,n)$ the sphere of radius $n$ centered at $v$ and by $B(v,n)$ the ball of same radius and center. 
We denote by $K_v$ the stabilizer ${\rm Stab}_G(v)$. 

\vspace{2mmplus1mmminus1mm} {\bf Definition.~}\it 
A subgroup of tree automorphisms $G<{\rm Aut}(X)$ is called {\rm locally $\infty$-transitive~}if for any vertex $v \!\in\! X$ and any radius $n$, the group $K_v$ is transitive on $S(v,n)$. 
\rm\vspace{2mmplus1mmminus1mm} 

Each vertex $v$ has a {\it type~}$\tau_v \!\in\! \{ 0;1 \}$ so that two adjacent vertices $v$ and $w$ are such that $\tau_v\neq\tau_w$. 
From now on, we assume that the tree $X$ admits a closed locally $\infty$-transitive group $G$ of
automorphisms. 
This implies in particular that $X$ is semi-homogeneous, i.e. there are at most two possible valencies (one for each type of vertex).   
We also assume that the group $G$ is type-preserving, which can be
done after passing to a subgroup of index at most 2.  
The maximal compact subgroups in $G$ are the vertex fixators; they are open. 
We will use the following transitivity properties \cite[Lemma 3.1.1]{BurMozTrees}. 

\begin{lemma} 
\label{lemma - transitivity trees} 
Let $X$ be a locally finite tree and let $G<{\rm Aut}(X)$ be a closed subgroup of tree
automorphisms. The following are equivalent. 
\begin{enumerate} 
\item[(i)] The group $G$ is locally $\infty$-transitive. 

\item[(ii)] For every vertex $v \!\in\! X$, the group $K_v$ is transitive on
$\partial_\infty X$. 

\item[(iii)] The group $G$ is non-compact and transitive on $\partial_\infty X$. 

\item[(iv)] The group $G$ is $2$-transitive on $\partial_\infty X$. 
\qed
\end{enumerate} 
\end{lemma}  

\subsubsection{}
These transitivity properties have deep combinatorial consequences. 
To state this, we let $G$ be a closed locally $\infty$-transitive subgroup of ${\rm Aut}(X)$. 
We choose a geodesic line $L$ (defining two boundary points $\xi,\eta \!\in\!\partial_\infty X$) and an edge $E=[v;v']$ in $L$ (defining two adjacent vertices $v$ and $v'$).  
We set: $P_\xi={\rm Fix}_G(\xi)$, $N_L={\rm Stab}_G(L)$ and $\mathscr{B}_E={\rm Fix}_G(E)$. 

\begin{lemma}
\label{lemma - Tits system trees} 
There exists $s_v\!\in\! G$ $($resp. $s_{v'}\!\in\! G)$ fixing
$v$ $($resp. $v'$$)$ and switching $\xi$ and $\eta$, so that: 

\begin{enumerate} 
\item[(i)] the quadruple $(G, P_\xi, N_L,\{s_v\})$ is a spherical Tits system with Weyl group 
${\bf Z}/2{\bf Z}$; 
\item[(ii)] the quadruple $(G, \mathscr{B}_E, N_L, \{ s_v; s_{v'} \} )$ is an affine Tits system with Weyl group $D_\infty$. 
\end{enumerate}
\end{lemma}

Recall that $D_\infty$ is the infinite dihedral group. 

\vspace{2mmplus1mmminus1mm} 

{\it Proof}.~
We prove the existence of a suitable symmetry $s_v$ as above. 
For any radius $n \!\in\! {\bf N}$, by local $\infty$-transitivity there exists $g_n\!\in\!G$ whose restriction to $B(v,n)$ is a symmetry $s_v$ around $v$ stabilizing the diameter $[x;x']=L \cap B(v,n)$: first use transitivity on $S(v,n)$ to get $g_n'\!\in\!K_v$ sending $x$ on $x'$; then use the sphere centered at $x'$ and of radius twice bigger to get $g_n''\!\in\!K_v$ sending $g_n'x'$ on $x$. 
The automorphism $g_n=g_n''g_n'$ is an approximation of $s_v$ on $B(v,n)$. 
All the elements $g_n$ are in the compact subgroup $K_v$, so the sequence $\{g_n \}_{n \in {\bf N}}$ admits a cluster value which can be chosen for $s_v$. 
The same argument works for $v'$. 

\vspace{2mmplus1mmminus1mm} 

(i).~We prove the axioms (T1)-(T4) of a Tits system \cite[IV.2]{Lie456}. 
We already know that there is an element $s\!\in\!N$ switching $\xi$ and $\eta$ since we can choose $s_v$ or $s_{v'}$. 
(We can also use 2-transitivity of $G$ on $\partial_\infty X$.) 
Let $g\!\in\!G$. Assume that $g.\xi\neq\xi$. 
Then, by 2-transitivity of $G$ on $\partial_\infty X$, there exists $p\!\in\!P_\xi$ such that $p^{-1}g.\xi = \eta$, which implies $sp^{-1}g.\xi = \xi$. 
This proves the Bruhat decomposition $G=P_\xi\sqcup P_\xi sP_\xi$, which implies the first half of (T1): $G=\langle P_\xi,s\rangle$, as well as (T3). 
Axiom (T4) is clear since $sP_\xi s^{-1} = {\rm Fix}_G(\eta)\neq P_\xi$. 
The group $P_\xi \cap N_L$ is nothing else than the pointwise fixator of $\{ \xi;\eta \}$: it is normal in the global stabilizer $N_L$. 
This implies the second half of (T1), and (T2). 

\vspace{2mmplus1mmminus1mm} 

(ii).~The tree $X$ is a building with Weyl group $D_\infty$; its apartments are the geodesic lines and its chambers are the edges. 
The group $N_L$ contains the reflections $s_v$ and $s_{v'}$, so it is transitive on the edges in $L$. 
Combined with the 2-transitivity of $G$ on $\partial_\infty X$, this implies the transitivity of $G$ on pairs of edges at given distance. 
Therefore $G$ acts strongly transitively on $X$ with respect to  $L$, which implies (ii) by \cite[Theorem 5.2]{Ronan}. 
\qed

\subsubsection{}
\label{sss - subgroups trees} 
As a consequence $K_v$ is a parabolic subgroup (in the combinatorial sense) of the affine Tits system of (ii): $K_v = \mathscr{B}_E \sqcup \mathscr{B}_E s_v \mathscr{B}_E$. 
Moreover by irreducibility of both Weyl groups, $P_\xi$ and $K_v$ are maximal subgroups of $G$ \cite[IV.2]{Lie456}.  
After defining some additional subgroups, we obtain further decompositions of $G$ and of some subgroups. 

\vspace{2mmplus1mmminus1mm} 

{\bf Definition.~}\it 
Let $\xi\!\in\!\partial_\infty X$ and let $\tau$ be a hyperbolic translation of step $2$. 
\begin{enumerate} 
\item[(i)] We denote by $D_\xi$ the subgroup of $P_\xi$ stabilizing each horosphere centered at $\xi$. 
\item[(ii)] We set $T_\tau=\langle\tau\rangle$ and we denote by $ \overline T_\tau^+$ the semigroup $\{\tau^n\}_{n\geq 0}$. 
\end{enumerate} 
\rm\vspace{2mmplus1mmminus1mm} 

For Bruhat-Tits buildings, there is a dictionary between apartments and maximal split tori. 
The definition of the subgroup $T$ here (i.e. the analogue of a maximal split 
torus) not only depends on the choice of the geodesic line $L$, but also on that of $\tau$. 
Nevertheless, if $\xi$ (resp. $\eta$) denotes the attracting (resp.
repelling) point of $\tau$, we may use the notation $T_{\xi,\eta}$ instead of $T_\tau$. 
The choice of $\tau$ among other hyperbolic translations along $(\xi \eta)$ is usually harmless. 
In spite of these slight differences with the algebraic definitions, there are analogues of
well-known decompositions in Lie groups. 
In the algebraic case, the group $D_\xi$ is bigger than a unipotent radical; the difference is explained in further detail in \ref{sss - connection with BT}. 

\begin{prop}
\label{prop - decomposition} 
We denote by $\partial_\infty^2 T$ the product $\partial_\infty T \times \partial_\infty T$ minus its diagonal. 
\begin{enumerate} 
\item[(i)] For any $(\xi', \eta') \!\in\! \partial_\infty^2 T$,the group $G$ contains a hyperbolic translation $\tau_{\xi',\eta'}$ of step $2$, with attracting $($resp. repelling$)$ point $\xi'$ $($resp. $\eta'$$)$. 

\item[(ii)] The group $D_\xi$ is transitive on every horosphere centered at $\xi$ and $G$ has an {\rm Iwasawa decomposition:~} $G = K_v \cdot  T_\tau \cdot D_\xi$, for any hyperbolic translation $\tau$ with attracting $($resp. repelling$)$ point $\xi$ $($resp.
$\eta)$. 

\item[(iii)] The group $P_\xi$ is amenable; it has a semidirect product decomposition: $P_\xi =  T_\tau \ltimes D_\xi$. 

\item[(iv)] The group $G$ has a {\rm Cartan decomposition:~}$G = K_v \cdot \overline  T_\tau^+ \cdot K_v$. 
\end{enumerate} 
\end{prop}  

We note that (iii) and (iv) are proved and used in \cite{LubMoz} to prove the vanishing at $\infty$ of matrix
coefficients of some unitary representations of ${\rm Aut}(X)$ (Howe-Moore property). 

\vspace{2mmplus1mmminus1mm}

{\it Proof}.~
(i).~The automorphism $s_v \circ s_{v'}$ or its inverse is a required hyperbolic translation along the geodesic line $L$. 
The case of an arbitrary geodesic line $L'$ follows by conjugation since $G$ is 2-transitive on $\partial_\infty X$. 

\vspace{2mmplus1mmminus1mm} 

(ii).~Let $v$ and $v'$ be two points on the same horosphere centered at $\xi$.  
Denote by $\{ v_n \}_{n \geq 0}$ the set of vertices of $[v \xi) \cap [v' \xi)$, with $\beta_{v_0,\xi}(v_n)=-n$ where $\beta_{v_0,\xi}$ is the Busemann function associated to the ray $[v_0 \xi)$ \cite[II.8.17]{BriHae}. 
For each $n \geq 0$, $v$ and $v'$ are on the same sphere centered at $v_n$, so by local $\infty$-transitivity, there exists 
$g_n \!\in\! G$ mapping $[v_n;v]$ onto $[v_n;v']$. 
Each element $g_n$ fixes $[v_0;v_n]$ therefore lies in $K_{v_0}$.  
Any cluster value of $\{ g_n \}_{n \geq 0}$ is an element of $D_\xi$ sending $v$ to $v'$, so that $D_\xi$ is transitive on every horosphere centered at $\xi$. 
This gives the Iwasawa decomposition: let $g \!\in\! G$; by the previous result, $g.v$ can be sent by some $u \!\in\! D_\xi$ to a point of $(\xi \eta)$, and by type-preservation, a suitable power of $\tau$ sends $(ug).v$ to $v$. 

\vspace{2mmplus1mmminus1mm} 

(iii).~The amenability is proved in \cite{LubMoz}. 
Let $g \!\in\! P_\xi$ and let $S$ be a horosphere centered at $\xi$. 
By type-preservation a suitable power $\tau^n$ sends $g.S \cap (\eta \xi)$ onto $S \cap (\eta \xi)$, so $\tau^n g\!\in\!D_\xi$ (once an element of $P_\xi$ stabilizes a horosphere centered at $\xi$, it stabilizes all of them). 
The argument also shows that $\langle \tau \rangle$ normalizes $D_\xi$; this proves the semidirect product assertion since the only power $\tau^n$ stabilizing $S$ is $1$. 
 
\vspace{2mmplus1mmminus1mm} 

(iv) is proved by using successive approximations of an automorphism in order to send $v$ back onto itself, using the compactness of $K_v$ and its transitivity on all spheres around $v$. 
\qed

\subsection{Compactifications}
\label{ss - CompTrees}
We define three ways to compactify the tree $X$, by means of measures, closed subgroups and gluings, respectively. 
All these compactifications will be identified in the next section. 

\subsubsection{}
\label{sss - measure-theoretic compactification trees} 
We first deal with the measure-theoretic compactification of trees. 
First, as an easy consequence of the transitivity properties, of the Iwasawa decomposition of $G$ and of the
amenability of $P_\xi$, we have: \vspace{2mmplus1mmminus1mm} 

\begin{lemma} 
\label{lemma - F boundary trees} 
The geometric boundary $\partial_\infty X$ is the maximal Furstenberg boundary of $G$. 
\end{lemma} 

{\it Proof}.~The notion of Furstenberg boundary is recalled in \ref{sss - parabolic and asymptotic}. 
Minimality is satisfied since $\partial_\infty X$ is homogeneous under $G$.  
The dynamics of hyperbolic translations on $\partial_\infty X$ \cite[II.8.16]{GhysHarpe} and the Lebesgue dominated convergence theorem imply that: $\displaystyle \lim_{n \to \infty} \tau^n_*\mu=\delta_\xi$, for any probability measure $\mu$ on $\partial_\infty X$ and any hyperbolic translation $\tau$ with attracting point $\xi$, provided that the reppelling point of $\tau$ is not an atom of 
$\mu$.  
So by Proposition \ref{prop - decomposition} (i), the closure of the $G$-orbit of any probability measure $\mu$ on 
$\partial_\infty X$ contains a Dirac measure. 
This proves that $\partial_\infty X$ is a Furstenberg boundary of $G$. 
At last, by the Iwasawa decomposition of Proposition \ref{prop - decomposition} (ii), we can write $G=KTD$ with $K$ compact and $TD$ amenable, and it follows from this  \cite[4.4]{Furstenberg} that every Furstenberg boundary of $G$ is an equivariant image of 
$G/TD\simeq\partial_\infty X$.  
\qed\vspace{2mmplus1mmminus1mm} 

The above lemma says that $\partial_\infty X$ plays for $G$ the role of a maximal flag variety for a simple algebraic group. 
The next lemma is another analogy in this spirit. 

\begin{lemma}
\label{lemma - embedding proba trees}
Let $X$ be a semi-homogeneous tree and let $G$ be a closed locally $\infty$-transitive group of automorphisms. 
\begin{enumerate}
\item[(i)] To each vertex $v \!\in\! X$ is associated a unique probability measure $\mu_v$ on $\partial_\infty X$ whose
fixator is precisely the maximal compact subgroup $K_v$. 
\item[(ii)] The assignment $\mu : v \mapsto \mu_v$ defines an embedding of the discrete set of
vertices of $X$ into the space of probability measures $\mathscr{M}^1(\partial_\infty X)$. 
\end{enumerate}
\end{lemma}

{\it Proof}.~
(i).~By transitivity of $K_v$ on $\partial_\infty X$, there is a unique $K_v$-fixed probability measure 
$\mu_v\!\in\!\mathscr{M}^1(\partial_\infty X)$ \cite[Lemma 1.4]{Raghu}. 
Since $K_v$ is a maximal subgroup, if ${\rm Fix}_G(\mu_v)$ were strictly bigger than $K_v$, it would be the whole group $G$. 
This is impossible since $G$ contains hyperbolic translations, and any such $\tau$ satisfies: 
$\displaystyle \lim_{n \to \infty} \tau^n_*\mu_v = \delta_\xi$, where $\xi$ is the attracting point of $\tau$. 

\vspace{2mmplus1mmminus1mm} 

(ii).~By (i) there is a one-to-one $G$-equivariant correspondence between the measures $\mu_v$ 
and the maximal compact subgroups $K_v$, hence a one-to-one $G$-equivariant correspondence between the measures
$\mu_v$ and the vertices. 
By uniqueness of the measure attached to a vertex, we have $\mu_{g.v}=g_*\mu_v$. 
Assume now that there is a cluster value $\nu$ in the subset $\{ \mu_v \}_{v \in V_X}$ of $\mathscr{M}^1(\partial_\infty X)$, so that $\nu$ is the limit of an injective sequence $\{ \mu_{v_n} \}_{n \geq 1}$. 
This provides an injective sequence of vertices $\{ v_n \}_{n \geq 1}$, which has to go to $\infty$ by discreteness of $V_X$. 
By the Cartan decomposition of Proposition \ref{prop - decomposition} (iv), there is a geodesic ray $[v\xi)$, 
a subsequence $\{ v_{n_j} \}_{j \geq 1}$ in $[v\xi)$ going to $\infty$ 
and $\{ k_j \}_{j \geq 1}$ a sequence in $K_v$ converging to some $k$, such that 
$\displaystyle \lim_{j \to\infty}k_j^{-1} v_{n_j} = \xi$.  
This implies $\nu = \delta_{k.\xi}$, but the latter measure is not fixed by any maximal compact subgroup: contradiction. 
\qed\vspace{2mmplus1mmminus1mm} 

The two previous lemmas lead us to the following natural definition of the Furstenberg compactification for trees. 

\vspace{2mmplus1mmminus1mm} 

{\bf Definition.~}\it 
The closure of the image of the map $\mu$ is called the {\rm measure-theoretic compactification~}of the set of vertices $V_X$ of $X$. 
It is denoted by $\overline{V}_X^{\rm meas}$. 
\rm\vspace{2mmplus1mmminus1mm} 

\subsubsection{}
\label{sss - group-theoretic compactification trees} 
We now define the group-theoretic compactification, using the space of closed subgroups $\mathscr{S}(G)$, endowed with the Chabauty topology. 

\begin{prop}
\label{prop - groups trees} 
Let $R=[v\xi)$ be a geodesic ray in the tree $X$. 
\begin{enumerate}
\item [(i)] Let $\{v_n \}_{n \geq 1}$ be a sequence of vertices in $R=[v\xi)$ going to $\infty$. 
Then the sequence of maximal compact subgroups  
$\{K_{v_n}\}_{n \geq 1}$ converges in $\mathscr{S}(G)$ to the subgroup $D_\xi$. 
\item [(ii)] The set $\mathscr{K}(G)$ of maximal compact subgroups of $G$ is discrete in $\mathscr{S}(G)$, so the assignment 
$K\!\mid_{V_X}:v \mapsto K_v$ defines an embedding of the set of vertices $V_X$ into $\mathscr{S}(G)$. 
\end{enumerate}
\end{prop} 

{\it Proof}.~
(i).~By compactness of $\mathscr{S}(G)$, it is enough to show that any cluster value of $\{K_{v_n}\}_{n \geq 1}$ is equal to $D_\xi$. 
Let $\displaystyle D = \lim_{j \to \infty}K_{v_{n_j}} < G$ be such a closed subgroup. 
Choose a geodesic line $(\xi \eta)$ extending $R$ and a step 2 hyperbolic translation $\tau$ along $(\xi \eta)$, with attracting point $\xi$. 
At last, fix $v'$ a vertex of $(\xi \eta)$ adjacent to $v$. 
After passing to a subsequence, we may -- and shall -- assume that $K_{v_{n_j}} = \tau^{n_j} K_{v''} \tau^{-n_j}$ where 
$\{n_j \}_{j \geq 1}$ is a sequence of positive integers such that $\displaystyle \lim_{j \to \infty}n_j=\infty$ and where $v''=v$ or $v'$. 
Then by the Lebesgue's dominated convergence theorem, $\displaystyle \lim_{n \to \infty} \tau^n_*\mu_{v''}=\delta_\xi$ implies
that $D$ fixes $\delta_\xi$ hence $\xi$. 
This implies $D<P_\xi$. 
Conversely, let $g\!\in\!D$. 
Using $D<P_\xi$ and Proposition \ref{prop - decomposition} (iii), we can write $g = u \tau^N$, with $N \!\in\!{\bf Z}$ and 
$u \!\in\!D_\xi$. 
As an element of a limit group, $g$ can also be written $\displaystyle g = \lim_{j \to \infty} \tau^{n_j} k_j \tau^{-n_j}$, for a sequence 
$\{ k_j \}_{j \geq 1}$ of elements of $K_{v''}$. 
Therefore there exists $J \geq 1$ such that for any $j \geq J$, we have: 
$(u \tau^N).v'' = (\tau^{n_j} k_j \tau^{-n_j}).v''$. 
Since $u$ stabilizes any horosphere centered at $\xi$, there is a vertex $z$ in $(\xi \eta)$ such that $(u\tau^N).v''$ and $\tau^N.v''$ are on the same sphere centered at $z$. 
Hence, we may -- and shall -- choose $j$ large enough to have $d\bigl( \tau^{n_j}.v'', (u \tau^N).v'' \bigr)=2n_j-2N$.  
But the group $\tau^{n_j} K_{v''} \tau^{-n_j}$ stabilizes the spheres centered at $\tau^{n_j}.v''$, which implies that $d\bigl( \tau^{n_j}.v'', (\tau^{n_j} k_j \tau^{-n_j}).v''\bigr) = 2n_j$.  
Thus in order to have $(u \tau^N).v'' = (\tau^{n_j} k_j \tau^{-n_j}).v''$, we must have $N=0$, hence $g=u$.  
This shows that $D=D_\xi$. 
This proves (i), which together with the same argument as for Lemma \ref{lemma - embedding proba trees} (ii), implies (ii).  
\qed\vspace{2mmplus1mmminus1mm} 

{\bf Definition.~}\it 
The closure of the image of the map $K\!\mid_{V_X}$ is called the {\rm group-theoretic compactification~}of the set of vertices $V_X$ of $X$. 
It is denoted by $\overline{V}_X^{\rm gp}$. 
\rm

\subsubsection{}
\label{sss - polyhedral theoretic compactification trees} 
The last compactification to be defined is the polyhedral one. 
As for Bruhat-Tits buildings, we compactify the whole tree $X$ by extending an equivalence relation defining $X$ as a gluing 
(\ref{ss - polyhedral}). 
Taking the closure of the set of vertices gives a compact space to be compared with the previous compactifications. 
This is done in \ref{ss - identification and amenable, trees}. 

\vspace{2mmplus1mmminus1mm} 

Let us consider the closure $\overline L \subset X^{\rm geom}$ of the geodesic line $L=(\eta\xi)$ containing the standard edge 
$E=[v;v']$. 
The subspace $\overline L$ admits a $D_\infty$-action via the restriction map $N_L\to N_L \!\mid_{\overline L}$. 
For the sake of homogeneity of notation, for any $x\!\in\!\overline L$ we define $G_x$ to be $K_x$ if $x\!\in\!L$ and to be $D_x$ if $x\!\in\! \overline L \setminus L=\{\xi;\eta\}$. 
We define the binary relation $\sim$ by: $(g,x) \sim (h,y)$ if, and only if, there exists $n \!\in\! N_L$ such that $y = n.x$ and 
$g^{-1}hn\!\in\!G_x$.
It is easy to see that $\sim$ is an equivalence relation. 

\vspace{2mmplus1mmminus1mm}

 {\bf Definition.~}\it  
We define the quotient space $\displaystyle {G \times\overline L \over \sim}$ to be the {\rm polyhedral compactification~}of $X$. 
We denote it by $\overline{X}^{\rm pol}$. 
\rm\vspace{2mmplus1mmminus1mm}

We denote by $[g,x]$ the class of $(g,x)$ and by $\pi:G\times\overline L\to\overline{X}^{\rm pol}$ the natural projection. 
By definition, we have: $[g,x]=[gh,x]$ for any $h\!\in\!G_x$ and $[n,x]=[1,n.x]$ for any $n\!\in\!N$.
The group $G$ acts on $\overline{X}^{\rm pol}$ by setting: $h.[g,x]=[hg,x]$. 
We can also define the map: 

\vspace{2mmplus1mmminus1mm} 

\centerline{
$\begin{array}{rrrr} 
\hfill \phi: \hfill & \hfill G \times \overline L \hfill & \hfill\to\hfill &
\hfill \overline X^{\rm geom} \hfill\\ 
& \hfill (g,x) \hfill & \hfill\mapsto\hfill & \hfill g.x, \hfill 
\end{array}$
}

\vspace{2mmplus1mmminus1mm} 

where $g.x$ denotes the $G$-action on $ \overline X^{\rm geom}$. 
Now let $g$, $h\!\in\!G$ and $x$, $y\!\in\!\overline L$. 
If $(g,x)\sim(h,y)$, then $y = n.x$ and $g^{-1}hn \!\in\! G_x$ for some $n\!\in\!N_L$. 
Setting $g_x=g^{-1}hn$, we obtain: $h.y=hn.x=gg_x.x=g.x$. 
Conversely, if $h.y=g.x$ the equality $x=g^{-1}h.y$ implies that $y=n.x$ for some $n \!\in\! N_L$. 
Then $h.y=g.x$ writes $g.x=hn.x$, so that $g^{-1}hn.x=x$. 
Since $\phi$ is surjective, by factorizing it through $\pi$, we obtain a $G$-equivariant bijection 
$\overline\phi:\overline X^{\rm geom}\cong\overline{X}^{\rm pol}$. 

\begin{lemma} 
\label{lemma - polyhedral compactification trees} 
The space $\overline{X}^{\rm pol}$ is compact, so the factorization map $\overline\phi$ is a $G$-homeomorphism and 
$\overline{X}^{\rm pol}$ is a compactification of the tree $X$. 
\end{lemma} 

{\it Proof}.~Let us denote by $\overline R=[v\xi]$ the closure of the geodesic ray $R =[v \xi)$ in $\overline X^{\rm geom}$. 
By local $\infty$-transitivity of $G$, it is a fundamental domain for the action of $K_v$ on $\overline X^{\rm geom}$. 
Since $\overline\phi$ is a $G$-equivariant bijection, this shows that the restricted projection map 
$\pi: K_v \times \overline R \to \overline{X}^{\rm pol}$ is surjective. 
Hence, in order to conclude that $\overline{X}^{\rm pol}$ is compact, we need to show that it is
Hausdorff, i.e. we need to prove that the graph of $\sim$ is closed. 
Since $(g,x) \sim (g',x')$ is equivalent to $g.x=g'.x'$ in $\overline X^{\rm geom}$, this
comes from the continuity of the $G$-action on the geometric compactification. 
\qed

\subsection{Identification and amenable subgroups}
\label{ss - identification and amenable, trees} 
We identify all the previously defined compactifications, and we recall that we can use them to parametrize maximal amenable subgroups of the automorphism group $G$. 

\subsubsection{}
In the statement below, $\overline{V}_X^{\rm geom}$ (resp. $\overline{V}_X^{\rm pol}$) denotes the closure of the set of vertices in the geometric compactification $\overline X^{\rm geom}$ (resp. in the polyhedral compactification $\overline{X}^{\rm pol}$). 

\begin{prop} 
\label{prop - comparison trees} 
Let $X$ be a semi-homogeneous tree, with set of vertices $V_X$. 
Let $G$ be a closed locally $\infty$-transitive subgroup of ${\rm Aut}(X)$. 
Then, the following compactifications of $V_X$ are $G$-homeomorphic.  
\begin{enumerate} 
\item[(i)] The geometric compactification $\overline{V}_X^{\rm geom}=V_X \sqcup \partial_\infty X$. 
\item[(ii)] The polyhedral compactification $\overline{V}_X^{\rm pol}$. 
\item[(iii)] The group-theoretic compactification 
$\overline{V}_X^{\rm gp}=\{K_v\}_{v\in V_X}\sqcup\{D_\xi\}_{\xi\in\partial_\infty X}$.  
\item[(iv)] The measure-theoretic compactification 
$\overline{V}^{\rm meas}=\{\mu_v\}_{v\in V_X}\sqcup\{\delta_\xi\}_{\xi\in\partial_\infty X}$. 
\end{enumerate} 
\end{prop}  

{\it Proof}.~
By Lemma \ref{lemma - polyhedral compactification trees}, we already have a $G$-homeomorphism: 
$\overline\phi:\overline{V}_X^{\rm pol}\simeq\overline{V}_X^{\rm geom}$. 
Setting $\overline{V}_L=(L \cap V_X) \cup \{\xi;\eta\}$, we have: $\overline{V}_X^{\rm pol}=\pi(G\times\overline{V}_L)$. 
The isomorphism between (ii) and (iii) follows from factorizing the map: 

\vspace{2mmplus1mmminus1mm} 

\centerline{
$\begin{array}{rrrr} 
\hfill \nu: \hfill & \hfill G \times \overline{V}_L \hfill & \hfill\to\hfill &
\hfill \mathscr{M}^1(\partial_\infty X) \hfill\\ 
& \hfill (g,x) \hfill & \hfill\mapsto\hfill & \hfill g_*\nu_x, \hfill 
\end{array}$
}

\vspace{2mmplus1mmminus1mm} 

where $\nu_x$ is the probability measure $\mu_x$, $\delta_\xi$ or $\delta_\eta$ according to whether $x\!\in\!L$, $x=\xi$ or $x=\eta$, respectively. 
Finally, the isomorphism between (ii) and (iv) follows from factorizing the map: 

\vspace{2mmplus1mmminus1mm} 

\centerline{
$\begin{array}{rrrr} 
\hfill H: \hfill & \hfill G \times \overline{V}_L \hfill & \hfill\to\hfill &
\hfill \mathscr{S}(G) \hfill\\ 
& \hfill (g,x) \hfill & \hfill\mapsto\hfill & \hfill gH_xg^{-1}, \hfill 
\end{array}$
}

\vspace{2mmplus1mmminus1mm} 

where $H_x$ is the closed subgroup $K_x$, $D_\xi$ or $D_\eta$ according to whether $x\!\in\!L$, $x=\xi$ or $x=\eta$, respectively. 
\qed

\subsubsection{}
Classifying maximal amenable subgroups of tree autmorphism groups was done in \cite[I.8.1]{FTN} by elementary geometric arguments. 
In our context, we find more natural to prove it by a Furstenberg lemma about supports
of limit measures:  

\begin{lemma} 
\label{lemma - Furstenberg} 
Let $\{g_n \}_{n \geq 1}$ be an unbounded sequence of tree automorphisms. 
Assume there are two probability measures $\mu, \nu$ on $\partial_\infty X$ such that 
$\displaystyle \lim_{n \to \infty} {g_n}_*\mu = \nu$. 
Then the support of the limit measure $\nu$ contains at most two points. 
\end{lemma}

{\it Reference}.~This is \cite[4.3]{LubMozZim}, or 
\cite[Lemma 2.3]{BurMozCAT} for general CAT$(-1)$-spaces. 
\qed

\begin{prop}
\label{prop - amenable, trees} 
Let $H$ be an amenable subgroup of ${\rm Aut}(X)$. Then, either $H$ fixes a vertex $v \!\in\! X$,
either it fixes a boundary point $\xi \!\in\! \partial_\infty X$ or it stabilizes a geodesic
line $L \subset T$. 
\end{prop}

{\it Proof}. By amenability, $H$ fixes a measure $\mu\!\in\!\mathscr{M}^1(\partial_\infty X)$. 
If $H$ is compact, it fixes a vertex $v\!\in\!V_X$. 
Otherwise by lemma \ref{lemma - Furstenberg} the support of $\mu$ contains at most two points. 
It is stabilized by $H$, and we obtain the last two possibilities
according to whether $\mid\!{\rm supp}(\mu)\!\mid\,\,=1$ or $2$. 
\qed\vspace{2mmplus1mmminus1mm} 

Note that for the above two results, the automorphisms needn't be type-preserving. 

\subsubsection{}
\label{sss - connection with BT}
The connection with Bruhat-Tits theory is the following. 
Let $G$ be a simple algebraic group over a non-archimedean local field $F$ of $F$-rank 1. 
Then the Bruhat-Tits building $X$ of $G_{/F}$ is a semi-homogeneous tree. 
Its valencies are of the form $1+q_F^r$ where $q_F$ is the cardinal of the residue field $\kappa_F$ and $r\geq1$. 
In this situation, the groups geometrically defined in \ref{ss - tree combinatorics} have interpretations in terms of algebraic group theory \cite{BorelTits65}, \cite[6.1]{BrT2}, \cite{Borel}. 

\vspace{2mmplus1mmminus1mm} 

Let $L$ be an apartment, i.e. a geodesic line, in $X$ and let $\xi$ and $\eta$ be the ends of $L$. 
To this apartment is attached a maximal $F$-split torus $T$ of $G$. 
The $G_F$-action on $X$ naturally provides the following chain of inclusions of closed subgroups: 

\vspace{2mmplus1mmminus1mm} 

\centerline{${\rm Fix}_{G_F}(L) \subsetneq {\rm Fix}_{G_F}(\{\xi;\eta\})=P_\xi \cap P_\eta \subsetneq {\rm Stab}_{G_F}(L)$.}

\vspace{2mmplus1mmminus1mm} 

All these groups can be described algebraically. 
First, the groups $P_\xi$ and $P_\eta$ are the two opposite parabolic subgroups containing $T$. 
Their intersection ${\rm Fix}_{G_F}(\{\xi;\eta\})$ is the reductive Levi factor $M_F=Z_{G_F}(T_F)$ with anisotropic semisimple factor $M'=[M,M]$. 
The group $M'_F$ is compact. 
The stabilizer ${\rm Stab}_{G_F}(L)$ is the normalizer $N_{G_F}(T_F)$, and the fixator ${\rm Fix}_{G_F}(L)$ is equal to $M'_F\cdot T_{\rm cpt}$, where $T_{\rm cpt}$ is the unique maximal compact subgroup of $T_F$. 
The quotient group ${\rm Stab}_{G_F}(L)/{\rm Fix}_{G_F}(L)=N_{G_F}(T_F)/(M'_F\cdot T_{\rm cpt})$ is the affine Weyl group $D_\infty$ of $G_F$. 
The intersection $P_\xi \cap P_\eta$ is the subgroup of $N_{G_F}(T_F)$ which doesn't switch $\xi$ and $\eta$. 
The quotient group $N_{G_F}(T_F)/(P_\xi \cap P_\eta)$ is the spherical Weyl group ${\bf Z}/2{\bf Z}$ of $G_{/F}$, and the quotient group $(P_\xi \cap P_\eta)/Z_{G_F}(T_F)$ is free abelian of rank one. 
Geometrically, a generator of the latter group corresponds to a step 2 hyperbolic translation along $L$. 
This paragraph illustrates \ref{sss - apartment and standard subgroups}. 

\vspace{2mmplus1mmminus1mm} 

The algebraic situation provides another $G$-action, namely the (linear) adjoint representation 
${\rm Ad}:G \to {\rm GL}(\mathfrak{g})$. 
The Lie algebra $\mathfrak{g}_F$ admits a direct sum decomposition into three ${\rm Ad}(T_F)$-stable summands: 
$\mathfrak{g}_F=\mathfrak{g}_F^+ \oplus \mathfrak{m}_F \oplus \mathfrak{g}_F^-$, where $\mathfrak{m}_F$ is the subspace on which $T_F$ acts trivially. 
Note that $\mathfrak{m}_F$ is also the Lie algebra of the reductive anisotropic kernel $M_F$, and we have: 
$\mathfrak{m}_F=[\mathfrak{m}_F,\mathfrak{m}_F] \oplus \mathfrak{t}_F$, where $\mathfrak{t}_F={\rm Lie}(T_F)$. 
For the two other summands, there is a character $\alpha$ of $T$ defined over $F$ and such that any $t\!\in\!T_F$ acts via the adjoint action on $\mathfrak{g}_F^+$ (resp. $\mathfrak{g}_F^-$) by multiplication by $\alpha(t)$ (resp. $\alpha(t)^{-1}$). 
We can pick an element $t\!\in\!T_F$ inducing a step 2 hyperbolic translation along $L$. 
Up to replacing $t$ by its inverse, we may -- and shall -- assume that the attracting point of $t$ is $\xi$, and we also assume that the signs $\pm$ in $\mathfrak{g}_F^\pm$ have been chosen so that $\mid\! \alpha(t) \!\mid_F>1$. 
Let $U_\xi$ (resp. $U_\eta$) be the unipotent root group with Lie algebra $\mathfrak{g}_F^+$ (resp. $\mathfrak{g}_F^-$). 
We have: $U_\xi<P_\xi$ since ${\rm Lie}(P_\xi)=\mathfrak{m}_F\oplus \mathfrak{g}_F^+$. 
The adjoint action of $t$ on $\mathfrak{g}_F^+$ is expanding, whereas it is contracting on $\mathfrak{g}_F^-$. 
This can be seen geometrically as follows.  
For each vertex $v\!\in\!L$ we can define the subgroup $U_{\xi,v}=\{g\!\in\!U_\xi : g$ fixes the geodesic ray $[v\xi)$ pointwise$\}$. 
This provides a filtration on $U_\xi$ closely related to the valuated root datum structure put on $G_F$ by Bruhat-Tits theory. 
The smaller is $n\!\in\!{\bf Z}$, the bigger is the geodesic ray $[t^n.v\xi)$ fixed by the group $t^nU_{\xi,v}t^{-n}=U_{\xi,t^n.v}$, and vice versa. 
This paragraph illustrates \ref{sss - valuated root data and half spaces}. 

\vspace{2mmplus1mmminus1mm} 

The lemma below relates the contraction property of the adjoint action of an element $g\!\in\!P_\xi$ and the action on the horospheres centered at $\xi$. 
The latter action is used in \ref{sss - subgroups trees} to define the subgroup $D_\xi$ of $P_\xi$. 
Recall that a group $H$ acts {\it distally~}on vector space $V$ over a local field $F$ via a linear representation $\varphi$ if the eigenvalues of any element $\varphi(h)$, $h\!\in\!H$, have absolute value 1 (see also \ref{sss - def distality} for more details). 

\begin{lemma} 
\label{lemma - distality and limit groups for trees}
In the above setting and with the above notation, we have: \begin{enumerate}
\item [(i)] The limit group $D_\xi$ is the biggest subgroup of $P_\xi$ with distal adjoint action on $\mathfrak{g}_F$. 
\item [(ii)] We have the semidirect product decomposition: $D_\xi=(M'_F\cdot T_{\rm cpt})\ltimes U_\xi$, which can also be written: $D_\xi={\rm Fix}_{G_F}(L)\ltimes U_\xi$
\end{enumerate}
\end{lemma} 

{\it Proof}.~(i) follows from Lemma \ref{lemma - rank one (distal)}, and (ii) follows from the fact that the limit group $D_\xi$ can be computed in two ways: from the algebraic viewpoint by Theorem \ref{th - CV groups} and from the geometric viewpoint by Proposition \ref{prop - groups trees}. 
\qed\vspace{2mmplus1mmminus1mm} 

Finally, we note that it is not hard to check that the $G_F$-action on $X$ is locally $\infty$-transitive. 
First, we can invoke the general fact that $G_F$ has an affine Tits system providing a Euclidean building on which it acts strongly transitively \cite[\S 5]{Ronan}. 
Since ${\rm rk}_F(G)=1$, the affine Weyl group of this Tits system is the infinite dihedral group and the building under consideration is a tree 
\cite[2.7]{TitsCorvallis}. 
Strong transitivity of the $G_F$-action in this case amounts to transitivity on pairs of vertices at given distance from one another (with respect of types), and this implies what we need. 
We can also say that given $v$ a vertex in $X$, the stabilizer $K_v$ is equal to the $\mathscr{O}_F$-points of some group scheme over $\mathscr{O}_F$ whose reduction modulo $\varpi_F$ is a finite $\kappa_F$-group acting strongly transitively on the neighbours of $v$ \cite[5.1.32]{BrT2}. 
This implies (2-)transitivity of the $G_F$-action on spheres of radius one. 
For bigger radii, one uses moreover that $U_{\xi,v}$ fixes $[v\xi)$ and acts transitively on the vertices at given distance from $v$ and different from the one in $[v\xi)$. 
This folding argument also proves that $(\eta\xi)$ is a fundamental domain for the $U_\xi$-action on the tree $X$ (Iwasawa decomposition). 

\section{Geometric parametrization of remarkable subgroups}
\label{s - parametrization}

Back to the algebraic situation, where $G$ is a simply connected semisimple $F$-group of arbitrary positive $F$-rank, we use the previously defined compactifications of the Bruhat-Tits building $X$ to parametrize two classes of remarkable subgroups of $G_F$. 
The first class consists of the amenable closed subgroups with connected Zariski closure. 
The second class consists of the subgroups acting without any contraction on the Lie algebra $\mathfrak{g}_F$ of $G_F$ (via the adjoint representation). 

\subsection{Amenable subgroups} 
\label{ss - amenable}
Compactifications of Bruhat-Tits buildings can be used to parametrize amenable subgroups in $G_F$. 
In the case of real semisimple Lie groups, this was proved by C.C. Moore, see \cite{MooreCompactification} and \cite{MooreAmenable}. 

\subsubsection{}
\label{sss - statement amenable}
A survey on amenable groups is for instance \cite[I.5]{Margulis}. 
The reference \cite[4.1]{Zimmer} will be sufficient for our purposes. 
Here is our geometric classification result. 

\begin{theorem}
\label{th - parametrization amenable}
Let $G$ be a semisimple simply connected algebraic group defined over a locally compact
non-archimedean local field $F$. 
Let $X$ be the Bruhat-Tits building of $G_{/F}$ and let $\overline{X}^{\rm pol}$ be its polyhedral compactification. 
\begin{enumerate} 
\item[(i)] Any closed, amenable, Zariski connected subgroup of $G_F$ fixes a facet in $\overline{X}^{\rm pol}$. 
\item[(ii)] Conversely, the stabilizer of any facet in $\overline{X}^{\rm pol}$ is an amenable Zariski connected subgroup. 
\item[(iii)] The closed amenable Zariski connected subgroups of $G_F$, maximal for these properties, are the vertex fixators for
the $G_F$-action on the compactification $\overline{X}^{\rm pol}$. 
\end{enumerate} 
\end{theorem}

Since maximal compact subgroups in $G_F$ are Zariski dense in $G$, hence connected, this theorem is an extension of the one-to-one correspondence between maximal compact subgroups of $G_F$ and vertices in $X$ \cite[Chap. 3]{BrT1}. 
Note that for any minimal parabolic $F$-subgroup $Q$, the group $Q_F$ is amenable and Zariski connected. 
Since any subgroup of $G_F$ containing $Q_F$ is a parabolic subgroup with non-compact semisimple quotient, $Q_F$ is maximal for these properties. 
Moreover any semisimple Levi factor of such a $Q$ is anisotropic over $F$, i.e. its group of $F$-rational points is compact and its Bruhat-Tits building is a point \cite[Theorem 1]{Prasad}, \cite[5.1.27]{BrT2}. 
This point, appearing in the boundary of $\overline{X}^{\rm pol}$, is of course considered as a facet. 
Maximal compact and minimal parabolic subgroups provide the two extreme cases (at least with respect to the dimension of the Zariski closure) of the above geometric parametrization. 

\begin{cor} 
\label{cor - finite orbit for amenable}Ê
In the above situation, any closed amenable subgroup of $G_F$ has a finite orbit in the compactification $\overline{X}^{\rm
pol}$. 
\end{cor}Ê

Note that since there is no non-positively curved distance on the compactification $\overline{X}^{\rm pol}$, an amenable
subgroup may not have a fixed point in $\overline{X}^{\rm pol}$.  
This is illustrated by the example of the normalizer of a maximal $F$-split torus (\ref{sss - Z connectedness necessary}). 

\vspace{2mmplus1mmminus1mm} 

{\it Proof}.~Let $R$ be a closed amenable subgroup of $G_F$. 
If $H$ denotes the identity component of the Zariski closure $\overline{R}^Z$, the intersection $R^\circ=H\cap R$
is a finite index normal subgroup of $R$, which is Zariski connected. 
Therefore, $R^\circ$ fixes a point $x$ in $\overline{X}^{\rm pol}$ by Theorem \ref{th - parametrization amenable}. 
The orbit $R.x$ has at most $[R:R^\circ]$ elements. 
\qed 

\subsubsection{}
The end of this subsection is devoted to proving the result. 
The first step consists of several reductions. 

\vspace{2mmplus1mmminus1mm} 

{\it Proof}.~
(ii). The groups under consideration are the conjugates $kR_{I,\underline d}k^{-1}$ 
where $k\!\in\!K_o$, $I$ is a set of simple roots and $\underline d$ is a 
family of non-negative real numbers indexed by $I$ (\ref{ss - stabilizers and orbits}). 
The amenablility of $R_{I,\underline d}$ is clear since this group is a 
compact-by-solvable extension of topological groups. 
Moreover we have: $\overline{R_{I,\underline d}}^Z=P_I$ for 
any subset of simple roots $I$ and any family of parameters $\underline d$; 
and parabolic subgroups are Zariski connected \cite[Theorem 11.16]{Borel}. 

\vspace{2mmplus1mmminus1mm} 

We concentrate on (i), and prove it 
by induction on ${\rm rk}_F(G)$, the $F$-rank of the group $G_{/F}$. 
The induction hypothesis is the statement of (i) when ${\rm rk}_F(G) \leq n$. 
The case $n=1$, where the Bruhat-Tits buildings are trees, is dealt with in 
Proposition \ref{prop - amenable, trees}, see also \ref{sss - connection with BT}. 
We henceforth assume that ${\rm rk}_F(G) \geq 2$. 
Let $R$ be a closed, amenable, Zariski connected subgroup of $G_F$. 
If $R$ is compact, we are also done since by the Bruhat-Tits fixed point lemma, $R$ fixes a point in $X$ \cite[3.2.3]{BrT1}. 
We henceforth assume that $R$ is non-compact. 
Here is the most important reduction, which uses the geometric structure of the boundary of $\overline{X}^{\rm pol}$ 
(Theorem \ref{th - BT in boundary} and Lemma \ref{lemma - type-preserving}). 

\begin{lemma}
\label{lemma - in parabolic is enough (amenable)}
It suffices to show that $R$ is contained in some proper parabolic $F$-subgroup $P$. 
\end{lemma}

{\it Proof}.~
By Theorem \ref{th - BT in boundary}, the group-theoretic compactification $\overline{Y}^{\rm pol}$ of the Bruhat-Tits building $Y$ of the semisimple $F$-group $P/\mathscr{R}(P)$ naturally lies in the boundary of $\overline{X}^{\rm pol}$. 
The group $R$ acts on it via its image under the quotient map $q:P\to P/\mathscr{R}(P)$. 
Moreover the closure of the image of an amenable group by a continuous group homomorphism is again amenable \cite[Lemma
4.1.13]{Zimmer}.  
Since ${\rm rk}_F\bigl(P/\mathscr{R}(P)\bigr)<{\rm rk}_F(G)$ we can use the induction hypothesis to have an $R$-fixed point
in $\overline{Y}^{\rm pol}$, hence in the boundary of $\overline{X}^{\rm pol}$. 
\qed

\subsubsection{}
Our goal now is to prove that $R<P_F$ for some parabolic $F$-subgroup $P$ of $G$. 

\vspace{2mmplus1mmminus1mm} 

We choose an embedding of $F$-algebraic groups $G<{\rm GL}(V)$ where $V$ is a
finite-dimensional $F$-vector space. 
We see $R$ as a closed, amenable, Zariski connected subgroup of ${\rm GL}(V_F)$. 
An inductive use of \cite[Corollary 3.7]{Cornulier} implies the
existence of an $R$-stable flag 

\vspace{2mmplus1mmminus1mm} 

\centerline{$\{0 \}=V_0\subsetneq V_1 \subsetneq V_2\subsetneq ... \subsetneq V_d=V_F$,}

\vspace{2mmplus1mmminus1mm} 

where the image of the natural map $R \to {\rm PGL}(V_i/V_{i-1})$ is relatively compact for each 
$i\geq1$. 
We can therefore write: $R<KTU$, where $Q=\overline{KTU}^Z$ is a proper parabolic subgroup
of ${\rm GL}(V)$ defined over $F$, and $TU=\mathscr{R}(Q)_F$ (the subgroup $K$ is a maximal
compact subgroup of the semisimple Levi factor defined by the $R$-invariant flag) \cite[8.4.6, exercise 4]{Springer}. 

\vspace{2mmplus1mmminus1mm} 

Let us consider the commutative square: 

\vspace{2mmplus1mmminus1mm} 

\centerline{
$\begin{array}{rrr} 
\hfill Q\hfill&\hfill{\buildrel\pi\over\longrightarrow}\hfill&\hfill Q/\mathscr{R}(Q)\hfill\\ 
\hfill \cup \hfill & & \hfill \cup \hfill\\ 
\hfill G\cap Q\hfill&\hfill{\buildrel\pi\over\longrightarrow}\hfill 
&\hfill(G\cap Q)/\bigl(G\cap \mathscr{R}(Q)\bigr),\hfill
\end{array}$
}

\vspace{2mmplus1mmminus1mm} 

where the horizontal maps $\pi$ are quotient maps in the sense of \cite[\S 6]{Borel}. 

\begin{lemma}
\label{lemma - image controlled (amenable)}
The closure $\overline{\pi(R)}$, taken in $\bigl(Q/\mathscr{R}(Q)\bigr)_F$, is compact and contained in $\pi\bigl((G\cap Q)_F\bigr)$. 
\end{lemma}

{\it Proof}.~
The maps $\pi$ are separable, so $\pi\bigl((G\cap Q)_F\bigr)$ is open and closed in 
$\Bigl((G \cap Q)/\bigl(G \cap \mathscr{R}(Q)\bigr)\Bigr)_F$ \cite[II.3.18]{BorelTits73}, which itself is closed in 
$\bigl(Q/\mathscr{R}(Q)\bigr)_F$ \cite[I.2.1.3(i)]{Margulis}. 
Therefore $\pi\bigl((G\cap Q)_F\bigr)$ is closed in $\bigl(Q/\mathscr{R}(Q)\bigr)_F$, and since $R<(G\cap Q)_F$, we have: $\overline{\pi(R)}<\pi\bigl((G\cap Q)_F\bigr)$. 
We know that $R<KTU$, so $\overline{\pi(R)}<\pi(K)$. 
This proves the compactness assertion. 
\qed\vspace{2mmplus1mmminus1mm} 

The statement of the next lemma uses the notion of $F$-trigonalizability for $F$-subgroups of ${\rm GL}(V)$
\cite[15.3]{Borel}. 

\begin{lemma}
\label{lemma - triangular normalized}
The group $R$ normalizes a connected non-trivial $F$-trigonalizable $F$-subgroup $H<G$. 
\end{lemma}

{\it Proof}.~Let us denote by $p: (G\cap Q)_F \to (G\cap Q)_F/\bigl(G\cap \mathscr{R}(Q)\bigr)_F$ the restriction $\pi\!\!\mid_{(G\cap Q)_F}$. 
It is a continuous and surjective homomorphism of topological groups. 
Its image contains $\overline{\pi(R)}$ by Lemma \ref{lemma - image controlled (amenable)}. 
Let us set $\widetilde R=p^{-1}(\overline{\pi(R)})$. 
Since $p$ is surjective, we have: $p(\widetilde R)=\overline{\pi(R)}$. 
Since the groups are locally compact \cite[II.3.18]{BorelTits73}, we have an isomorphism of topological groups: 
$\widetilde R/(\widetilde R \cap TU) \simeq \overline{\pi(R)}$, where the first group is endowed with the quotient topology \cite[VII Appendice I Lemme 2]{Integration78}. 
Since $R$ is not compact, neither is $\widetilde R$, and the previous isomorphism implies that $\widetilde R \cap TU$ cannot be compact either. 
In particular, $\widetilde R \cap TU$ is infinite, so the identity component $H=\bigl(\overline{\widetilde R \cap TU}^Z\bigr)^\circ$ of the Zariski closure $\overline{\widetilde R \cap TU}^Z$ is an $F$-trigonalizable group of positive dimension. 
\qed\vspace{2mmplus1mmminus1mm} 

We conclude the proof of the theorem by the following. 

\begin{lemma}
\label{lemma - in parabolic (amenable)}
The group $R$ lies in a proper parabolic $F$-subgroup of $G$. 
\end{lemma}

{\it Proof}.~We note that $\mathscr{R}_u(H)$ is defined over $F$ \cite[15.4 (ii)]{Borel} and we distinguish
two cases. 

\vspace{2mmplus1mmminus1mm} 

First case:~the unipotent radical $U=\mathscr{R}_u(H)$ is non-trivial. 
If the characteristic of $F$ is zero, by \cite[Introduction]{BorelTits71} there is 
a parabolic subgroup $P=\mathscr{P}(U)$ such that 
$U<\mathscr{R}_u(P)$ and $N(U)<P$. 
In the case ${\rm char}(F)=p>0$, the residue field is anyway finite, hence perfect, so 
we have $[F:F^p]\leq p$. 
Therefore we can use \cite[Theorem 2]{Gille}: there is a parabolic subgroup $P'$ of $G$,
defined over $F$ and such that $U<\mathscr{R}_u(P')$. 
Denoting by $F_s$ the separable closure of $F$, we can choose a Borel subgroup $B_{/F_s}$ contained in $P'$
and  defined over $F_s$. 
Then $\mathscr{R}_u(B)$ contains $\mathscr{R}_u(P')$, which implies that $\mathscr{R}_u(H)$ is $F_s$-embeddable. 
Therefore we can use \cite[Introduction]{BorelTits71} also in this case to conclude that there is 
a parabolic $F$-subgroup $P=\mathscr{P}(U)$ such that 
$U<\mathscr{R}_u(P)$ and $N_{G}\bigl(\mathscr{R}_u(H)\bigr) \cap G_{F_s}<P$. 
Whatever the characteristic of $F$, we obtain: $R<\mathscr{P}(U)$. 

\vspace{2mmplus1mmminus1mm} 

Second case:~we have $\mathscr{R}_u(H)=\{1\}$. Then $H$ is an $F$-split torus \cite[15.4]{Borel}. 
Since it is normalized by $R$ which is Zariski connected, it is actually centralized by
$R$ \cite[8.10]{Borel}.  
But the centralizer of an $F$-split torus in a reductive $F$-group is a Levi factor of some
parabolic $F$-subgroup \cite[20.4]{Borel}. 
We also find in this case that $H$, hence $R$, lies in some proper parabolic
$F$-subgroup $P$. 
\qed

\subsection{Distal subgroups} 
\label{ss - distal}
Compactifying the Bruhat-Tits building $X$ can also be used to parametrize maximal distal subgroups in $G_F$. 
For this we need the very definition of the group-theoretic compactification of $V_X$, i.e. we need to use the compactification in which points are closed subgroups of $G_F$. 
We don't pass to stabilizers as in the previous subsection. 

\subsubsection{}
\label{sss - def distality} 
The notion of distality comes from topological dynamics; it is a natural extension of isometric actions on metric spaces \cite{FurstenbergDistal}. 
A very general definition in the context of group actions on uniform structures is given in \cite{Abels81}. 
Roughly speaking, requiring distality of a group action amounts to avoiding contractions (equivalently, expansions). 
The situations are different according to whether one considers distal actions on compact spaces or on vector spaces. 
On projective spaces, the notion opposite distality, i.e. proximality (\ref{sss - measure-theoretic compactification trees}), is a key tool in the proof of Tits' alternative \cite{TitsFree}; it is also an ingredient of the definition of a Furstenberg boundary \cite[VI.1]{Margulis}. 
 
\vspace{2mmplus1mmminus1mm} 
 
Here we are only interested in actions on vector spaces over valuated fields, arising from linear representations. 
In this case, a $H$-action arising from a linear representation $\rho:H\to{\rm GL}(V_F)$ is called {\it distal}~if the closure of any $\rho(H)$-orbit in $V_F$ is minimal (\ref{sss - measure-theoretic compactification trees}). 
This is equivalent to the fact that the eigenvalues of $\langle\rho(h)\rangle$, for any $h\!\in\!H$, are all of absolute value 1. 
(The eigenvalues of an endomorphism of some $F$-vector space are taken in a suitable finite extension of $F$; moreover the absolute value is implicitly extended, in a unique way, to the algebraic closure  $\overline F$.) 
One implication is obvious from computing iterations of endomorphisms \cite[II.1]{Margulis}, and the other implication is proved in \cite{ConGui} or \cite{Abels78} (see \cite{Cornulier} for fixing the confusion between irreducible and absolutely irreducible $H$-modules in \cite{ConGui}). 
\vspace{2mmplus1mmminus1mm} 

Finally, we denote by ${\rm Ad}:G \to {\rm GL}(\mathfrak{g})$ the adjoint representation of $G_{/F}$. 

\begin{theorem}
\label{th - parametrization distal}
Let $G$ be a semisimple simply connected algebraic group defined over a locally compact
non-archimedean local field $F$. 
Let $X$ be the Bruhat-Tits building of $G_{/F}$ and let $\overline V_X^{\rm gp}$ be the group-theoretic compactification of its vertices. 
\begin{enumerate} 
\item[(i)] Any subgroup of $G_F$ with distal adjoint action on $\mathfrak{g}_F$ is contained in a point of $\overline V_X^{\rm gp}$. 
\item[(ii)] The adjoint action of any limit group of $G_F$ is distal. 
\item[(iii)] The subgroups of $G_F$ with distal adjoint action on $\mathfrak{g}_F$ and maximal for these properties are the groups of $\overline V_X^{\rm gp}$; in particular they are closed and Zariski connected. 
\end{enumerate} 
\end{theorem}

We note that Zariski connectedness is not required in the assumptions of the above theorem. 

\subsubsection{}
We can now prove the parametrization of distal subgroups in a non-archimedean semisimple Lie group. 

\vspace{2mmplus1mmminus1mm} 

{\it Proof}.~
(iii).~It is a straightforward consequence of (i), (ii) and Lemma  \ref{lemma - stabilizers} asserting that $\overline{D_{I,\underline d}}^Z=P_I$ for each subset $I$ of simple roots and each family of parameters $\underline d$. 

\vspace{2mmplus1mmminus1mm} 

(ii).~The groups under consideration are the conjugates $kD_{I,\underline d}k^{-1}$ where $k\!\in\!K_o$, $I$ is a set of simple roots and $\underline d$ is a family of non-negative real numbers indexed by $I$ (\ref{ss - stabilizers and orbits}). 
The distality of $D_{I,\underline d}$ is clear since its adjoint image is a compact-by-unipotent extension of linear groups. 

\vspace{2mmplus1mmminus1mm} 

(i).~The proof goes by induction on ${\rm rk}_F(G)$, the $F$-rank of the group $G_{/F}$. 
The induction hypothesis is the statement of (i) when ${\rm rk}_F(G) \leq n$. 
The case $n=1$ is dealt with in the following lemma. 

\begin{lemma}
\label{lemma - rank one (distal)}
Assume ${\rm rk}_F(G)=1$. 
Let $D$ be a non-compact subgroup in $G_F$ whose adjoint action is distal. 
Then $D$ lies in a limit group. 
\end{lemma}

{\it Proof}.~
If $D$ is bounded, it fixes a point in the Bruhat-Tits tree $X$ of $G_{/F}$.  
We henceforth assume that $D$ is unbounded. 
By iterating Furstenberg's lemma \cite[Theorem 3.10]{Cornulier}, there exists an ${\rm Ad}(D)$-stable flag: 

\vspace{2mmplus1mmminus1mm} 

\centerline{$\{0 \}=V_0\subsetneq V_1 \subsetneq V_2\subsetneq ... \subsetneq V_d=\mathfrak{g}_F$,}

\vspace{2mmplus1mmminus1mm} 

such that the image of the natural map ${\rm Ad}(D)\to{\rm GL}(V_i/V_{i-1})$ 
is relatively compact for each $i\geq1$. 
This implies that the closure $\overline{{\rm Ad}(D)}$ of ${\rm Ad}(D)$ 
in ${\rm GL}(\mathfrak{g}_F)$ is an amenable subgroup. 
Since ${\rm Ker}({\rm Ad}\!\!\mid_{G_F})$ is finite, 
${\rm Ad}\!\!\mid_{G_F}$ is a proper map, so the group 
$({\rm Ad}\!\!\mid_{G_F})^{-1}\bigl(\overline{{\rm Ad}(D)}\bigr)$ is
amenable too, and so is $\overline D$ as a closed subgroup of the latter group. 

\vspace{2mmplus1mmminus1mm} 

We henceforth use the notation of \ref{sss - connection with BT}. 
We claim that  $\overline D$ fixes a point, say $\xi$, in the ideal boundary $\partial_\infty X$. 
Otherwise, by Proposition \ref{prop - amenable, trees} the group $\overline D$ would stabilize a geodesic line $L$ and switch its two ends. 
In other words, it would be contained in the normalizer $N'_F$ of a maximal $F$-split torus $T'_F$ and its image under the natural map 
$N'_F \to {\bf Z}/2{\bf Z}$ would be non-trivial. 
Since distality prevents $D$ from containing a hyperbolic translation, we conclude that $D$ would be contained in the extension of ${\bf Z}/2{\bf Z}$ by the maximal compact subgroup of $T'_F$, hence would be bounded: contradiction. 

\vspace{2mmplus1mmminus1mm} 

It follows from the previous paragraph that $D$ is contained in a proper (equivalently, minimal) parabolic subgroup $P_\xi$ of $G_F$. 
Let us choose a maximal (i.e. one-dimensional) $F$-split torus $T$ in $P_\xi$ (equivalently, a geodesic line in $X$ with one end equal to $\xi$). 
This provides a decomposition: $P_\xi = (M'_F \cdot T_F) \ltimes U_\xi$, where $U_\xi=\mathscr{R}_u(P_\xi)_F$ and 
$M'$ is the semisimple Levi factor attached to $T$ such that $M'_F \cdot T_F=Z_{G_F}(T_F)$. 
The distality of $D$ in the adjoint action implies that the $T_F$-part of the decomposition 
of any element in $D$ must lie in the maximal compact subgroup $T_{\rm cpt}$ of $T_F$. 
This finally proves that $D$ is contained in the limit group $D_\xi$. 
\qed\vspace{2mmplus1mmminus1mm} 

We henceforth assume that ${\rm rk}_F(G) \geq 2$. 
Let $D$ be a subgroup of $G_F$ whose adjoint action on the Lie algebra $\mathfrak{g}_F$ of $G_F$ is distal. 
If $D$ is bounded, it fixes a point in $X$ \cite[VI.4]{Brown}, so we henceforth assume that $D$ is unbounded. 

\begin{lemma}
\label{lemma - in parabolic is enough (distal)}
It suffices to show that $D$ is contained in some proper parabolic $F$-subgroup $P$ of $G$. 
\end{lemma}

{\it Proof}.~
Let $D$ be a subgroup of $G_F$ with distal adjoint action, and such that some proper parabolic $F$-subgroup $P$ contains it. 
Distality is preserved by conjugation, so we may -- and shall -- assume that $P$ is the standard 
parabolic subgroup $P_I$ attached to the set of simple roots $I$. 
Let us denote by $q_I:P_I\to P_I/\mathscr{R}(P_I)$ the natural surjection. 
We also introduce $\mathfrak{p}_I$ the Lie algebra of $P_I$ and $\mathscr{R}(\mathfrak{p}_I)$ the Lie algebra of $\mathscr{R}(P_I)$. 
The distality of the adjoint action of $D$ on $\mathfrak{g}_F$ implies the distality of the adjoint action of 
$(q_I\!\mid_{G_F})^{-1}\bigl(q(D)\bigr)$ on $(\mathfrak{p}_I)_F/\mathscr{R}(\mathfrak{p}_I)_F\simeq {\rm Lie}(G_I)_F$. 
Since $G_I$ is a simply connected semisimple $F$-group of $F$-rank smaller than ${\rm rk}_F(G)$, we are in
position to  apply the induction hypothesis. 
We deduce from it that there is a maximal limit group in $(G_I)_F$ containing $(q_I\!\mid_{G_F})^{-1}\bigl(q(D)\bigr)$. 
In view of the description of the limit groups in $(G_I)_F$, 
this says that there exist $k\!\in\!K_o\cap G_I$, a subset of simple roots $J$ in $I$ and a family $\underline
d$ of non-negative real numbers indexed by $J$ such that 
$(q_I\!\mid_{G_F})^{-1}\bigl(q(D)\bigr)<k\bigl(K_{J,\underline d}\ltimes(U^J\cap G_I)\bigr)k^{-1}$, 
hence $q(k^{-1}Dk)<q\bigl(K_{J,\underline d}\ltimes(U^J\cap G_I)\bigr)$. 
Since $U^J=(U^J\cap G_I)\ltimes U^I$, this finally implies: 
$k^{-1}Dk<K_{J,\underline d}\ltimes U^J=D_{J,\underline d}$. 
\qed

\subsubsection{}
Our goal now is to prove that $D<P_F$ for some proper parabolic $F$-subgroup $P$ of $G$. 
Again by \cite[Theorem 3.10]{Cornulier} there exists an ${\rm Ad}(D)$-stable flag 

\vspace{2mmplus1mmminus1mm} 

\centerline{$\{0 \}=V_0\subsetneq V_1 \subsetneq V_2\subsetneq ... \subsetneq V_d=\mathfrak{g}_F$,}

\vspace{2mmplus1mmminus1mm} 

where the image of the natural map ${\rm Ad}(D)\to{\rm GL}(V_i/V_{i-1})$ is relatively compact for each 
$i\geq1$ (see also \cite{ConGui}). 
We can therefore write: ${\rm Ad}(D)<KU$, where $Q=\overline{KU}^Z$ is a proper parabolic $F$-subgroup
of ${\rm GL}(\mathfrak{g})$, $U=\mathscr{R}_u(Q)_F$ and $K$ is a maximal
compact subgroup of the reductive Levi factor defined by the ${\rm Ad}(D)$-invariant flag. 
Let us consider the commutative square: 

\vspace{2mmplus1mmminus1mm} 

\centerline{
$\begin{array}{rrr} 
\hfill Q\hfill&\hfill{\buildrel\pi\over\longrightarrow}\hfill&\hfill Q/\mathscr{R}_u(Q)\hfill\\ 
\hfill \cup \hfill & & \hfill \cup \hfill\\ 
\hfill {\rm Ad}(G)\cap Q\hfill&\hfill{\buildrel\pi\over\longrightarrow}\hfill 
&\hfill({\rm Ad}(G)\cap Q)/\bigl({\rm Ad}(G)\cap \mathscr{R}_u(Q)\bigr),\hfill
\end{array}$
}

\vspace{2mmplus1mmminus1mm} 

where the horizontal maps $\pi$ are quotient maps in the sense of \cite[\S 6]{Borel}. 
The diagram is similar to the one in the previous subsection, except that we replaced $\mathscr{R}(Q)$ by
$\mathscr{R}_u(Q)$. 
We conclude the proof of the theorem thanks to the last point of the following lemma, most of whose proof
mimicks the proofs of Lemmas \ref{lemma - image controlled (amenable)}, \ref{lemma - triangular normalized}
and \ref{lemma - in parabolic (amenable)}. 

\begin{lemma}
\label{lemma - similar to amenable}
With the above notation, the following holds. 
\begin{enumerate}
\item[(i)] The closure $\overline{\pi\bigl({\rm Ad}(D)\bigr)}$ in
$\bigl(Q/\mathscr{R}_u(Q)\bigr)_F$, is compact and  contained in $\pi\bigl({\rm Ad}(G)\cap Q_F\bigr)$. 
\item[(ii)] The group ${\rm Ad}(D)$ normalizes a non-compact unipotent subgroup $V$ of ${\rm Ad}(G_F)$. 
\item[(iii)] The group $D$ lies in a proper parabolic $F$-subgroup of $G$. 
\end{enumerate} 
\end{lemma}

{\it Proof}.~
(i).~Once $G$ is replaced by ${\rm Ad}(G)$ and $\mathscr{R}(Q)$ is replaced by $\mathscr{R}_u(Q)$, 
use the same arguments as for Lemma \ref{lemma - image controlled (amenable)}. 

\vspace{2mmplus1mmminus1mm} 

(ii).~Let us denote by $p: (G\cap Q)_F \to (G\cap Q)_F/\bigl(G\cap \mathscr{R}_u(Q)\bigr)_F$
the restriction $\pi\!\!\mid_{({\rm Ad}(G)\cap Q)_F}$. 
By (i), we have: 
$\overline{\pi\bigl({\rm Ad}(D)\bigr)}<{\rm Im}(p)$, so we can set: 
$\widetilde D=p^{-1}\Bigl(\overline{\pi\bigl({\rm Ad}(D)\bigr)}\Bigr)$. 
Then, as for Lemma \ref{lemma - triangular normalized}, we obtain that 
$V={\rm Ker}(p) \cap \widetilde D$ is a non-compact unipotent subgroup in ${\rm Ad}(G_F)$, 
normalized by ${\rm Ad}(D)$. 

\vspace{2mmplus1mmminus1mm} 

(iii).~If the characteristic of $F$ is 0, we can use \cite[Introduction]{BorelTits71} in ${\rm Ad}(G)$ to
obtain a parabolic $F$-subgroup $P'=\mathscr{P}(V)$ such that $V<\mathscr{R}_u(P')$ and $N_{{\rm Ad}(G)}(V)<P'$. 
This proves (iii) in this case because ${\rm Ad}^{-1}(P')$ is a parabolic $F$-subgroup of $G$ containing
$D$. 
We henceforth assume that ${\rm char}(F)=p>0$. 
Then the group $V$ is an infinite group of finite exponent, and so is its preimage 
$({\rm Ad}\!\!\mid_{G_F})^{-1}(V)$ since ${\rm Ker}({\rm Ad})$ is finite. 
Then it follows from \cite[Lemma VIII.3.7]{Margulis} that the identity component of the Zariski closure of 
$({\rm Ad}\!\!\mid_{G_F})^{-1}(V)$, say $\widetilde V$, is a unipotent group of positive dimension. 
As in the first case of the proof of Lemma \ref{lemma - in parabolic (amenable)}, we can combine 
\cite[Theorem 2]{Gille} and \cite[Introduction]{BorelTits71} 
to obtain a parabolic $F$-subgroup $P=\mathscr{P}(\widetilde V)$ such that
$\widetilde V<\mathscr{R}_u(P)$ and $N_{G}(\widetilde V)\cap G_{F_s}<P$. 
\qed

\subsection{Discussion of the hypotheses} 
\label{ss - hypotheses}
We discuss  the hypotheses of our last two main theorems on parametrization of remarkable subgroups. 
We show that Zariski connectedness is necessary to properly parametrize amenable subgroups. 
We also suggest simplified proofs for both theorems when the local ground field has characterisitic 0. 
Finally, we discuss the scope of all our results; in particular, we explain that our results and proofs sometimes improve the case of symmetric spaces. 

\subsubsection{}
\label{sss - Z connectedness necessary} 
The Zariski connectedness assumption is necessary to classify amenable subgroups in $G_F$ in terms of fixed facets (Theorem \ref{th - parametrization amenable}). 

\begin{lemma}Ê
\label{lemma - Z connected necessary}Ê
Let $T$ be the standard maximal $F$-split torus in $G$. 
\begin{enumerate}
\item [(i)] The fixed-point set for $T_F$ acting on $\overline{V}_X^{\rm gp}$ is the intersection of the closure of $A$ in $\overline{X}^{\rm pol}$ with the unique closed orbit in $\overline{V}_X^{\rm gp}$. 
\item [(ii)] The group $N_G(T)_F$ doesn't have any fixed point in $\overline{X}^{\rm pol}$. 
\end{enumerate}
\end{lemma} 

Note that $N_G(T)_F$ is amenable since it is the extension of the spherical Weyl group of $G_F$ by the abelian group $T_F$. 
This provides an amenable subgroup of $G_F$ with a finite orbit, but without any fixed facet in $\overline{X}^{\rm pol}$. 

\vspace{2mmplus1mmminus1mm} 

{\it Proof}.~Ê
(i). Since $T_F$ is non-compact, it cannot have any fixed point in the building $X$. 
Let $D$ be a limit group normalized by $T_F$. 
The Zariski closure $Q=\overline{D}^Z$ is a parabolic $F$-subgroup normalized by $T_F$, i.e. a fixed point for $T_F$ acting on the spherical building at infinity $\partial_\infty X$. 
These are the parabolic $F$-subgroups containing $T$, i.e. the stabilizers of the facets in the boundary $\partial_\infty A$. 
Moreover $Q$ has to be a minimal parabolic $F$-subgroup because $T_F$ cannot stabilize any maximal compact subgroup in the Levi factor of $Q$ associated to $T$ (i.e. doesn't stabilize any limit group $D$ such that $\overline{D}^Z=Q$) unless this Levi factor is itself compact. 
This shows that the limit groups normalized by $T_F$ are the groups $wD_\varnothing w^{-1}$, where $w$ ranges over the spherical Weyl group $W$. 

\vspace{2mmplus1mmminus1mm} 

(ii). It remains to note that if the normalizer $N_G(T)_F$ had a fixed facet in $\overline{X}^{\rm pol}$, this facet would contain a vertex fixed by $T_F$. 
But the vertices of the previous paragraph are obviously permuted by $W$. 
\qed

\subsubsection{}
The proof of each theorem of the last two subsections is easier when the characteristic of the local field
$F$ is 0. 
This is due to the fact that when ${\rm char}(F)=0$, which we henceforth assume, the unipotent radical
of an algebraic group defined over $F$, is itself defined over $F$. 
This remark is applied below to groups defined as Zariski closures. 

\vspace{2mmplus1mmminus1mm} 

Let $R$ be a closed, non-compact, amenable subgroup in $G_F$. 
We assume that the Zariski closure $\overline R^Z$, which we denote by $H$, is connected. 
Using Lemma \ref{lemma - in parabolic is enough (amenable)}, we shall show that 
$R$ lies in a proper parabolic $F$-subgroup of $G$.
If $\mathscr{R}_u(H)\neq\{1\}$, by \cite{BorelTits71} 
there is a parabolic $F$-subgroup $P$ such that $\mathscr{R}_u(H)<\mathscr{R}_u(P)$ and 
$R<N_G\bigl(\mathscr{R}_u(H)\bigr)<P$.   
Otherwise $H$ is a reductive $F$-group, so we can choose a faithful, irreducible representation 
$H \to {\rm GL}(V)$ defined over $F$, 
and we denote by $q$ the natural surjection ${\rm GL}(V_F)\to{\rm PGL}(V_F)$. 
Then $q(H_F)$ is closed in ${\rm PGL}(V_F)$, hence contains $\overline {q(R)}$, 
and we set $\widetilde R=(q\!\!\mid_{H_F})^{-1}\bigl(\overline{q(R)}\bigr)$. 
By Zariski density, the representation $V$ is also irreducible as an $R$-module, and it follows 
from Zariski connectedness of $R$ and the Furstenberg lemma that 
$\overline{q(R)}=q(\widetilde R)$ is compact \cite[Corollary 3.2.2]{Zimmer}. 
Since we have a homeomorphism 
$\overline{q(R)}\simeq\widetilde R/\bigl(\widetilde R\cap{\rm Ker}(q)\bigr)$, 
this implies that $R$ normalizes a non-compact group of scalar matrices in $G_F$. 
Therefore $R$ centralizes a non-trivial $F$-split torus in $G$, so it lies in (the Levi
factor of) some proper parabolic $F$-subgroup. 

\vspace{2mmplus1mmminus1mm} 

Let $D$ be a subgroup of $G_F$ with distal adjoint action on $\mathfrak{g}_F$. 
We denote by $H$ the identity component of $\overline D^Z$, and by $D^\circ$ the subgroup $D\cap H$ of finite index in $D$. 
If $\mathscr{R}_u(H)\neq\{1\}$, by \cite[Introduction]{BorelTits71} there is a parabolic $F$-subgroup $P$ such that $\mathscr{R}_u(H)<\mathscr{R}_u(P)$ and $D<N_G(H)<P$, and we can use Lemma \ref{lemma - in parabolic is enough (distal)}. 
Otherwise ${\rm Ad}(H)$ is reductive, ${\rm Ad}(D^\circ)$ is Zariski dense in ${\rm Ad}(H)$ and distal on $\mathfrak{g}_F$. 
By \cite[Lemma 1]{Prasad} the group ${\rm Ad}(D^\circ)$ is bounded, and since ${\rm Ad}\!\!\mid_{G_F}$ is proper, 
this implies  that $D^\circ$, hence $D$, is bounded. 

\subsubsection{}
\label{sss - symmetric spaces} 
We finish this section by mentioning two problems concerning the archimedean case of symmetric spaces. 
First the proofs in this section, together with their simplifications in characteristic 0, enable to prove analogues of Theorem \ref{th - parametrization amenable} and Theorem \ref{th - parametrization distal} in the case of Lie groups obtained as rational points of semisimple ${\bf R}$-groups. 
This applies to the connected components of the isometry groups of Riemannian symmetric spaces of non-compact type, since the latter groups are then semisimple center-free. 
Note that it is not clear that the existing litterature on compactifications of symmetric spaces contains the analogues of Theorem \ref{th - parametrization amenable} and Theorem \ref{th - parametrization distal}. 
In general, it would be interesting to check whether the present paper provides substantial simplifications to the more classical real case. 
It is clear that not all results go through: e.g. maximal compact subgroups in the real case are real points of Zariski closed ${\bf R}$-subgroups, while they are Zariski dense in the non archimedean case. 
Still, the geometric idea underlying the induction for the last two theorems (i.e. using the action of Levi factors on symmetric spaces of smaller rank in the boundary) may be useful. 

\vspace{2mmplus1mmminus1mm} 

Our second remark is that it may be possible to prove the geometric parametrization of amenable subgroups by using \cite{AdBa}. 
This question can be asked in both the real and the non archimedean cases. 

\section{The example of the special linear group}
\label{s - GL} 
We illustrate some of the previous results and techniques to the case of the special linear group. 
This is the opportunity to recall the concrete viewpoint of non-archimedean additive norms in order to introduce the corresponding Euclidean building. 
The convergence of maximal compact subgroups in the Chabauty topology can be checked by matrix computation. 
We draw a picture describing the boundary of a Weyl chamber in terms of upper triangular-by-block matrices (i.e. in terms of flags). 
Our approach is elementary, and we refer the reader to \cite[Planche I]{Lie456} for the connection with root systems. 

\subsection{The building in terms of additive norms and lattices}Ê
\label{ss - building of SL(n)} 
We recall the concrete definition of the building of ${\rm SL}_n(F)$ by analogy with the case of the symmetric space of 
${\rm SL}_n({\bf R})$. 

\subsubsection{}
\label{sss - additive norms} 
Let us fix $E$ a vector space of finite dimension $n$ over the local field $F$. 
As in the real case, we are interested in logarithms of norms \cite[2.9]{TitsCorvallis}: 

\vspace{2mmplus1mmminus1mm} 

{\bf Definition}.~\it 
An {\rm additive norm}~on $E$ is a map $\gamma:E\to{\bf R}\cup\{+\infty\}$ satisfying: 
\begin{enumerate} 
\item[(AN1)]  for any  $x \! \in \! E$, we have: $\gamma(x)=+\infty$ if, and only if, $x=0$;
\item[(AN2)] for any  $x \! \in \! E$ and $\lambda \! \in \! F$, we have: $\gamma(\lambda x)=\gamma (x)+v_F(\lambda)$; 
\item[(AN3)] for any $x, y \! \in \! E$, we have: $\gamma(x+y) \geq \inf \{ \gamma (x);\gamma (y) \}$. 
\end{enumerate}
\rm\vspace{2mmplus1mmminus1mm} 

A basic result is the analogue of Gram-Schmidt reduction: for any ultrametric norm $\parallel \! - \! \parallel$, there exist a basis $\{ e_i \}_{1 \leq i \leq n}$ and positive real numbers $\{ C_i \}_{1 \leq i \leq n}$, such that for any $x \! \in \! E$, we have: 
$\parallel \! x  \! \parallel =  \sup_{i\in\{1;2;...\, n\}} \{ C_i\,.\mid\!\!\lambda_i\!\!\mid\}$, where 
$x = \sum_{i=1}^n \lambda_i e_i$; we then say that the basis $\{ e_i \}_{1 \leq i \leq n}$ is {\it adapted}~to the norm $\parallel \! - \! \parallel$. 
By a result due to A. Weil, there always exists a basis simultaneously adapted to any pair of norms \cite[Proposition 1.3]{GoldmanIwahori}. 
Let us fix now a basis $\mathbb{B}=\{ e_i \}_{1 \leq i \leq n}$. 
For each $\{ c_i \}_{1 \leq i \leq n}\!\in\!{\bf R}^n$, we denote by $\parallel \! - \! \parallel_{\mathbb{B}, \{c_i\}}$ the ultrametric norm $\sum_{i=1}^n \lambda_i e_i \mapsto \sup_i \{ q^{c_i} \cdot \! \mid \! \lambda_i \! \mid \}$, and by $\gamma_{\mathbb{B}, \{c_i\}}$ the additive norm: $-\log_q \circ \parallel \! - \! \parallel_{\mathbb{B}, \{c_i\}}$, also defined by:  
$\sum_{i=1}^n \lambda_i e_i \mapsto \inf_i \{ v_F(\lambda_i)-c_i \}$. 
The set $\mathscr{N}_\mathbb{B} = \{ \gamma_{\mathbb{B}, \{c_i\}}: \{ c_i \} \in {\bf R}^n \}$ is an $n$-dimensional affine space for the action  
${\bf R}^n \times \mathscr{N}_\mathbb{B} \to \mathscr{N}_\mathbb{B}$ defined by 
$(\{d_i\},\gamma_{\mathbb{B}, \{c_i\}}) \mapsto \gamma_{\mathbb{B}, \{c_i+d_i\}}$. 
 
\vspace{2mmplus1mmminus1mm} 

{\bf Definition}.~\it 
We call the set of additive norms on $E$ the {\rm Goldman-Iwahori space~}$E$, and we denote it by $\mathscr{N}_E$. 
We denote by $X_E$ the quotient of  $\mathscr{N}_E$ in which two additive norms are identified whenever their difference is constant, and we call it the {\rm Bruhat-Tits building~}of ${\rm SL}(E)$
\rm\vspace{2mmplus1mmminus1mm}

The notion of an $\mathscr{O}_F$-lattice, i.e. of a free sub-$\mathscr{O}_F$-module generating $E$ over $F$, 
distinguishes in $\mathscr{N}_E$ some norms whose classes are the vertices of a simplicial structure on $X_E$. 
To any $\mathscr{O}_F$-lattice $M$ is associated an additive norm $\gamma_M$ by setting 
$\gamma_M (x)= \sup \{ n \! \in \! {\bf Z}:  x \! \in \! \varpi_F^n M \} $. 
We have: $\gamma_M(E)={\bf Z} \cup \{ \infty \}$, and conversely, if $\gamma \! \in \! \mathscr{N}_E$ takes integral values then $\gamma^{-1}({\bf N} \cup \{ \infty \})$ is an $\mathscr{O}_F$-lattice, which we denote by $M_\gamma$. 
The correspondence $M_\gamma \leftrightarrow \gamma_M$ is a bijection between the set $\mathscr{L}_E$ of $\mathscr{O}_F$-lattices in $E$ and the set of additive norms with integral values on $E$. 
It is equivariant for the natural left ${\rm GL}(E)$-actions on $\mathscr{L}_E$ and $\mathscr{N}_E$ (by precomposition in the latter case). 

\vspace{2mmplus1mmminus1mm} 

Let now $\gamma$ be an additive norm. 
We choose a basis $\mathbb{B} = \{ e_i \}$ adapted to it. 
Permuting the indices $i$ (resp. multiplying the vectors $e_i$ by powers of the uniformizer $\varpi_F$) corresponds to actions by monomial (resp. diagonal) matrices. 
Using these operations, we send $\gamma$ onto 
$\gamma_{\mathbb{B}, \{ c_i \}}$ with  $c_i \! \in \! [0;1[$ and $c_1 \leq c_2 \leq ... \leq c_n$. 
We set $c_0=0$ and $c_{n+1}=1$ and define the $\mathscr{O}_F$-lattices 
$M_i(\gamma)= \gamma^{-1}([-c_i; \infty])$ and the associated additive norms
$\gamma^{(i)}=\gamma_{M_i(\gamma)}$. 
Whatever the choice of the basis  $\mathbb{B}$, these matrix operations lead to the same ordered sequence of real numbers $(c_i) \! \in \! [0;1[^n$ and the same $\mathscr{O}_F$-lattices. 
Moreover the sequence $(c_i) \! \in \! [0;1[^n$ only depends on the ${\rm GL}(E)$-orbit of $\gamma$. 
If $\gamma \! \in \! \mathscr{N}_{\mathbb{B}'}$ for another basis $\mathbb{B}'$, then $\gamma^{(i)} \! \in \! \mathscr{N}_{\mathbb{B}'}$ for each $i$,  and we have: $\gamma = \sum_{i=0}^n (c_{i+1} - c_i) \gamma^{(i)}$.
Conversely, a family of $\mathscr{O}_F$-lattices $\{M_0; M_1;...\, M_n\}$  comes from an additive norm $\gamma$ $($i.e., $M_i=M_i(\gamma)$ for each $i)$ if and only if: 
$\varpi_F M_n \subset M_0 \subset M_1 \subset\,  ... \subset M_n$ and for any $i$ with  $M_i \neq M_{i+1}$, 
$ i = {\rm dim}_\kappa\bigl(M_i/\varpi_F M_n\bigr)$.
This says that any additive norm is the barycenter of a well-defined system of \og weighted $\mathscr{O}_F$-lattices\fg. 
Moding out by additive constants endows $X_E$ with a simplicial structure and a compatible ${\rm GL}(E)$-action. 
The above facts on $X_E$ are proved in \cite{GoldmanIwahori}. 

\subsubsection{}
\label{sss - BT building structure for GL}Ê
We henceforth see $\mathbb{B} = \{ e_i \}$ as an {\it ordered}~basis. 
We call the convex hull of the homothety classes of the above lattices $M_i$ the {\it closed facet~}associated to $\gamma$.  
The set of the homothety classes of additive norms  with the same flag of $\mathscr{O}_F$-lattices
$M_i$ is called the {\it open facet~}associated to $\gamma$. 
A facet of maximal dimension is called an {\it alcove}. 
The {\it apartment~}$A_\mathbb{B}$ associated to $\mathbb{B}$ is the set of the classes of the additive norms to which $\mathbb{B}$  is adapted. 
The set of vertices in $A_\mathbb{B}$ is denoted by $V_\mathbb{B}$. 
Let $\mathbb{E} \simeq {\bf R}^{n-1}$ be the quotient of the vector space ${\bf R}^n$ by the vector all of whose coordinates are $1$.  
Since the apartment $A_\mathbb{B}$ is the quotient of $\mathscr{N}_\mathbb{B}$ by additive constants, we have a map 
$\mathbb{E} \times A_\mathbb{B} \to A_\mathbb{B}$ defined by 
$\bigl( [\{d_i\}],[\gamma_{\mathbb{B}, \{c_i\}}] \bigr) \mapsto [\gamma_{\mathbb{B}, \{c_i+d_i\}}]$. 
This endows $A_\mathbb{B}$ with the structure of an $(n-1)$-dimensional affine space. 
We call {\it wall}~of $A_\mathbb{B}$ an affine hyperplane of the form: $\{c_i-c_j = r\}$ for some $i \neq j$ and $r \! \in \! {\bf Z}$; we call {\it wall}~of $X_E$ any ${\rm SL}(E)$-transform of a wall of $A_\mathbb{B}$. 
To any $\underline \nu= \{ \nu_i \}_{1 \leq i \leq n} \! \in \! {\bf Z}^n$ we associate the $\mathscr{O}_F$-lattice 
$L_{\underline \nu}=\bigoplus_{1 \leq i \leq n} \mathscr{O}_F \varpi_F^{\nu_i} e_i$,
and denote by $[L_{\underline \nu}]$ its homothety class. 
The subset of vertices in $A_\mathbb{B}$ is $V_\mathbb{B}=\{ [L_{\underline \nu}]: \underline \nu \! \in \! {\bf Z}^n \}$. 
Less canonically, identifying the sequences $\underline \nu$ such that $\nu_1=0$ with ${\bf Z}^{n-1}$, the map defined by: $\underline \nu \mapsto [L_{\underline \nu}]$ defines a bijection $ {\bf Z}^{n-1} \simeq V_\mathbb{B}$. 
We denote by $o$ the vertex associated to the null sequence. 

\vspace{2mmplus1mmminus1mm} 

The subset of $A_\mathbb{B}$ consisting of the classes of the additive norms $[\gamma_{\mathbb{B}, \{c_i\}}]$ with $c_i < c_{i+1}$ for each $i <n$ is called the {\it sector}~(or {\it Weyl chamber}) attached to $\mathbb{B}$. It is denoted by $\mathscr{Q}$ and we denote by $V_\mathscr{Q}$ the set of vertices contained in the closure $\overline{\mathscr{Q}}^X$ of $\mathscr{Q}$ in $X_E$. 
Given a subset $I$ of $\{ 1;2;...\, n-1\}$, the {\it sector face $\mathscr{Q}^I$}~is the subset of the classes $[\gamma_{\mathbb{B}, \{c_i\}}]$ in $\overline{\mathscr{Q}}^X$  satisfying $c_i=c_{i+1}$ for all $i \! \in \! I$. 
We call a codimension one sector face, say defined by $c_i=c_{i+1}$, a {\it sector panel}~and we denote it by $\Pi^i$. 
The closure $\overline{\mathscr{Q}}^X$ is a simplicial cone whose faces are the sector faces $\mathscr{Q}^I$;  the subset $V_\mathscr{Q}$ consists of the classes of lattices $[L_{\underline \nu}]$ with $\nu_1 \leq \nu_2 \leq ... \leq \nu_n$ and the convex hull of $V_\mathscr{Q}$ is $\overline{\mathscr{Q}}^X$. 
The non-decreasing finite sequences correspond to $V_\mathscr{Q}$ in the identification $V_\mathbb{B} \simeq {\bf Z}^{n-1}$. 
At last, we have an action of the symmetric group $\mathscr{S}_n$ on the additive norms by permuting the indices. 
It is compatible with moding out by the additive constants, so it defines an action on $A_\mathbb{B}$ for which 
$\overline{\mathscr{Q}}^X$ is a fundamental domain. 
The main result about $X_E$ is that $A_\mathbb{B}$ is a geometric realization of the tiling of ${\bf R}^{n-1}$ by regular simplices, and that 
the family of subcomplexes $A_\mathbb{B}$ when $\mathbb{B}$ ranges over the bases of $E$, is the apartment system of a Euclidean building structure on the space $X_E$.
These facts justify {\it a posteriori}~our use of the building terminology in this paragraph and the previous one \cite{BruhatTitsGL}.  

\subsubsection{}
\label{sss - building at infinity and Levi buildings for GL} 
Let us now illustrate \ref{sss - parabolic and asymptotic}, which deals with the spherical building at infinity; a detailed reference is
 \cite[V.8 and VI.9F]{Brown}. 
For instance, the above sector $\mathscr{Q}$ defines a chamber at infinity $\partial_\infty\mathscr{Q}$ of the spherical building $\partial_\infty X_E$; the upper triangular standard Borel subgroup $B$ is also ${\rm Fix}_{{\rm SL}(E)}(\partial_\infty\mathscr{Q})$. 
Similarly, the standard torus $T$ of determinant 1 matrices diagonal with respect to the basis $\mathbb{B}$ defining $A_\mathbb{B}$, is also ${\rm Fix}_{{\rm SL}(E)}(\partial_\infty A_\mathbb{B})$. 
Its normalizer $N$, generated by $T$ and the monomial matrices: 

\vspace{2mmplus1mmminus1mm} 

\centerline{
$N_i=\begin{pmatrix}Ê
{\rm id}_{i-1}&0&0\\
0&
\begin{pmatrix}Ê
0&1 \\ 
-1&0
\end{pmatrix}
&0\\ 
0&0&{\rm id}_{n-i-1}
\end{pmatrix}$
} 

\vspace{2mmplus1mmminus1mm} 

for $i \! \in \! \{1;2;...\, n-1\}$ is also ${\rm Stab}_{{\rm SL}(E)}(\partial_\infty A)$. 
Furthermore, let $I$ be a subset of $\{1;2;...\,n-1\}$. 
We have the equivalence relation $i \sim_I j$ on $\{1;2;...\,n-1\}$ for which $i < j$ are {\it $I$-equivalent~}if, and only if, $\{i;i+1;... j-1\} \subset I$. 
We write $\{1;2;...\,n-1\}$ as a disjoint union $\bigsqcup_{j=1}^m  I_j$ of consecutive intervals of integers, and we set $d_j = \mid \! I_j \! \mid$, so that $\sum_{j=1}^m d_j =n$. 
The standard parabolic subgroup $P_I$ is the fixator of $\partial_\infty\mathscr{Q}^I$, and it is also the group of determinant $1$ matrices which are upper triangular by block and whose $j$-th block has size $d_j$. 
The standard reductive Levi subgroup $M_I$ is the fixator of $\partial_\infty{\rm Vect}(\mathscr{Q}^I)$, i.e. of the asymptotic classes of the geodesic rays contained in ${\rm Vect}(\mathscr{Q}^I)$: this is also the fixator of the union of the facet at infinity 
$\partial_\infty\mathscr{Q}^I$ and of its opposite in $\partial_\infty A_\mathbb{B}$, or the group of determinant $1$ matrices which are
diagonal by block and whose $j$-th block has size $d_j$. 
To illustrate completely \ref{ss - Bruhat-Tits} and \ref{ss - Levi}, let us mention that the standard semisimple Levi factor $G_I=[M_I,M_I]$ consists of the diagonal by block matrices whose $j$-th block has size $d_j$ and determinant equal to $1$, that $T^I$ is the subgroup of the torus $T$ consisting of the scalar by block matrices whose $j$-th (scalar) diagonal block has size $d_j$ and finally that the unipotent radical $U^I$consists of the upper triangular by block matrices whose $j$-th diagonal block is the $d_j \times d_j$ identity matrix. 

\vspace{2mmplus1mmminus1mm} 

In order to illustrate simply \ref{sss - Levi BT buildings}, let us consider the case where $n=3$ and $I$ reduces to $\{1\}$. 
In other words, we are interested in the non-essential realization of the Bruhat-Tits building of the upper left ${\rm GL}_2$ block inside the Bruhat-Tits building of ${\rm SL}_3(F)$. 
In this case, the affine subspace $L_{\{1\}}$ of $A_{\mathbb{B}}$ is the line $\{c_1=c_2\}$. 
It is easy to check that its ${\rm GL}_2(F)$-transforms in $A_{\mathbb{B}}$ are the straight lines $\{c_1-c_2=r\}$ when $r$ ranges over ${\bf Z}$. 
We can use the elementary unipotent subgroups $U_{12}$ and $U_{21}$ to construct ${\rm GL}_2(F).A_{\mathbb{B}}$. 
The so-obtained space is the product of a tree $T$ by $L_{\{1\}}$, the vertices of $T$ being the ${\rm GL}_2(F)$-transforms of $L_{\{1\}}$, and the edges being the minimal strips between two such lines. 
After \og shrinking the inessential direction $L_{\{1\}}$\fg, we obtain the Bruhat-Tits tree of ${\rm GL}_2(F)$. 

\subsection{Convergence and compactification}Ê
\label{ss - convergence and compactification} 
After describing concretely the parahoric subgroups of ${\rm SL}(E)$ and other subgroups related to the affine Tits system, we indicate how convergence of canonical sequences of parahoric subgroups in the Chabauty topology can be proved by elementary matrix computation. 

\subsubsection{}
\label{sss - Tits systems for GL} 
We denote by $c$ the alcove in $\overline{\mathscr{Q}}^X$ whose closure contains $o$: it is the alcove corresponding to the flag of $\mathscr{O}_F$-lattices $M_i$ where $M_n=\bigoplus_j \mathscr{O}_F e_j$, 
$M_i= \bigoplus_{j \leq n-i} \varpi_F\mathscr{O}_Fe_j \oplus \bigoplus_{j >n- i} \mathscr{O}_F e_j$ for $0<i<n$ and  $M_0=\varpi_FM_n$. 
We set: $K= {\rm Fix}_{{\rm SL}(E)}(o)$ and $\mathscr{B}= {\rm Fix}_{{\rm SL}(E)}(c)$: these groups are the standard maximal compact subgroup and Iwahori subgroup, respectively. 
We identify ${\rm SL}(E)$ to ${\rm SL}_n(F)$ via $\mathbb{B}$, so that $K$ corresponds to ${\rm SL}_n(\mathscr{O}_F)$ and $\mathscr{B}$ to the subgroup of ${\rm SL}_n(\mathscr{O}_F)$ reducing to the upper triangular matrices of ${\rm SL}_n(\kappa)$ modulo $\varpi_F$. 
Since ${\rm SL}(E)$ acts transitively on the alcoves of $X_E$, the Iwahori subgroups are the conjugates of $\mathscr{B}$. 
The standard parahoric subgroups are also defined as subgroups of ${\rm SL}_n(\mathscr{O}_F)$ with the condition to be a parabolic subgroup 
modulo $\varpi_F$. 

\vspace{2mmplus1mmminus1mm} 

We denote by $K_I$ the intersection of the standard reductive Levi factor $M_I$ with the maximal compact subgroup $K$, and we denote respectively by $D_I$ and $R_I$ the semi-direct products $K_I \ltimes U^I$ and $(K_I\cdot T_I) \ltimes U^I$. 
In matrix notations,  this gives: 

\vspace{2mmplus1mmminus1mm} 

\centerline{
$D_I=\{Êg \! \in \! 
\begin{pmatrix}
{\rm GL}_{d_1}(\mathscr{O}_F) & \ast & \cdots & \ast  \\
0 & {\rm GL}_{d_2}(\mathscr{O}_F) & \ast & \cdots \\
\cdots & 0 & \cdots & \ast \\
0 & \cdots & 0 & {\rm GL}_{d_m}(\mathscr{O}_F)
\end{pmatrix} 
: {\rm det}(g)=1 \}$ 
}

\vspace{2mmplus1mmminus1mm} 

and 

\vspace{2mmplus1mmminus1mm} 

\centerline{
$R_I=\{Êg \! \in \! 
\begin{pmatrix}
k^\times\cdot{\rm GL}_{d_1}(\mathscr{O}_F) & \ast & \cdots & \ast \\
0 & k^\times\cdot{\rm GL}_{d_2}(\mathscr{O}_F) & \ast & \cdots \\
\cdots & 0 & \cdots & \ast \\
0 & \cdots & 0 & k^\times\cdot{\rm GL}_{d_m}(\mathscr{O}_F)
\end{pmatrix}
 : {\rm det}(g)=1 \}$.}

\vspace{2mmplus1mmminus1mm} 

From this, the Zarisiki density of $D_I$ in $P_I$ is obvious. 
The  subgroup $T_I \cap K$ is infinite: its elements are scalar by block matrices with coefficients in $\mathscr{O}_F^\times$. 
The group $K_I$ consists of the determinant 1 matrices in $\prod_{j=1}^m {\rm GL}_{d_j}(\mathscr{O}_F)$, and $G_I \cap K_I$ is a maximal compact subgroup of $G_I$, naturally isomorphic to $\prod_{j=1}^m {\rm SL}_{d_j}(\mathscr{O}_F)$. 
In order to sum up the main combinatorial properties of ${\rm SL}(E)$ in terms of Tits systems, we need to introduce the further \og almost monomial\fg matrix: 

\vspace{2mmplus1mmminus1mm} 

\centerline{
$N_0=\begin{pmatrix}
0& \cdots \quad 0 &-\varpi_F^{-1}\\
\begin{matrix}
\cdots \\ 0
\end{matrix}Ê
&{\rm id}_{n-2}& 
\begin{matrix}
0 \\ \cdots
\end{matrix}\\
\varpi_F& 0 \quad \cdots &0
\end{pmatrix}.$
}

\vspace{2mmplus1mmminus1mm} 

It follows from Gauss reduction that $\bigl( {\rm SL}(E), B, N, \{ s_i \}_{1 \leq i \leq n-1} \bigr)$ is a Tits system with associated Coxeter system $(\mathscr{S}_n, \{ s_i \}_{1 \leq i \leq n-1})$ \cite[IV.2.2]{Lie456}. 
Moreover the spherical building at infinity $\partial_\infty X_E$ is a geometric realization of the combinatorial
building associated to this Tits system  \cite[VI.9F]{Brown}. 
From the interpretation of $X_E$ in terms of additive norms, we see that the ${\rm SL}(E)$-action on $X$ by precomposition is strongly transitive \cite[V.1]{Brown}. 
This implies that $\bigl( {\rm SL}(E), \mathscr{B}, N, \{ s_i \}_{0 \leq i \leq n-1} \bigr)$ is a Tits system whose Weyl group is an affine reflection group with linear part the symmetric group $\mathscr{S}_n$ \cite[V.1F, 1G, 2A]{Brown}. 

\subsubsection{}
\label{sss - parametrization for GL} 
We use more precise combinatorics for the group ${\rm SL}(E) \simeq {\rm SL}_n(F)$, seen as a subset of
the $n \times n$ matrices ${\rm M}_n(F)$ with basis $\{ E_{ij} \}_{1 \leq i,j \leq n}$, where $E_{ij}$ is the elementary matrix 
$[\delta_{k,i} \cdot \delta_{l,j}]_{1 \leq k,l \leq n}$. 
Given any ordered sequence of scalars  $\underline\lambda \! \in \! (k^\times)^n$, we denote by ${\rm Diag}(\lambda_1, ... \lambda_n)$ the corresponding diagonal matrix with respect to the ordered basis $\mathbb{B}$. 
Given any $\underline \nu \! \in \! {\bf Z}^n$, we denote by  $\varpi_F^{\underline \nu}$ the diagonal matrix ${\rm Diag}(\varpi_F^{\nu_1}, ... \varpi_F^{\nu_n})$. 
We also introduce the group $\Lambda=\{ \varpi_F^{\underline \nu}: \nu_1=0 \}$ and the semigroup 
$\overline \Lambda^+=\{ \varpi_F^{\underline \nu}: 0 = \nu_1 \leq \nu_2 \leq ... \leq \nu_n \}$. 
For any $\underline \nu \! \in \! {\bf Z}^n$, we have: $\varpi_F^{\underline \nu}.L_0=L_{\underline \nu}$; this shows that $\Lambda$ is simply transitive on the 
vertices of $A_\mathbb{B}$, and the orbit map $\varpi_F^{\underline\nu} \mapsto [L_\nu]$ for the origin $o=[L_0]$ provides a bijection 
$\Lambda \simeq V_\mathbb{B}$ which identifies the vertices in $\overline{\mathscr{Q}}^X$ and the semigroup $\overline \Lambda^+$. 
To keep on using subgroups of ${\rm SL}_n(F)$ exclusively, we also introduce the discrete subgroup $T_0 = \{ \varpi_F^{\underline \nu}: \sum_i \nu_i=0\}$ and the discrete semigroup $\overline T^+_0 = \{ \varpi_F^{\underline \nu}: \sum_i \nu_i=0$ and 
$\nu_1\leq \nu_2 \leq ... \leq \nu_n\}$, which both lie in $T$. 

\vspace{2mmplus1mmminus1mm} 

Given $I \subset \{ 1;2;...\, n-1\}$ we can also define $\mathscr{A}^I$ to be the set of affine subspaces in the
apartment $A_\mathbb{B}$ which are intersections of $\mid \! I \! \mid$ distinct walls parallel to ${\rm Vect}(\mathscr{Q}^I)$; e.g., ${\rm Vect}(\mathscr{Q}^I)$ belongs to $\mathscr{A}^I$, and $\mathscr{A}^{\{ 1;2;...\, n-1\}}$ is the 
set $V_\mathbb{B}$ of vertices in $A_\mathbb{B}$. 
The group $\Lambda$ acts transitively on $\mathscr{A}^I$; moreover any affine subspace of $\mathscr{A}^I$ intersecting $\overline{\mathscr{Q}}^X$ can be written $\varpi_F^{\underline \nu}.{\rm Vect}(\mathscr{Q}^I)$ for some $\varpi_F^{\underline \nu}$ in $\overline\Lambda^+$. 
The set of $T_0$-orbits in $\mathscr{A}^I$ is finite, and there is a finite subset $\{ÊE_1;E_2;...\, E_m \}$ of $\mathscr{A}^I$ such that 
any affine subspace of $\mathscr{A}^I$ intersecting the semigroup can be written $\varpi_F^{\underline \nu}.E_j$ for some $j$ and some 
$\varpi_F^{\underline \nu}$ in the semigroup $\overline T_0^+$. 
The affine subspaces $E_j$ as above can themselves be written $t_j.{\rm Vect}(\mathscr{Q}^I)$ with 
$t_j \! \in \! \Lambda$. 
Therefore any affine subspace $E$ of $\mathscr{A}^I$ can be written $E=(t_0t_j).{\rm Vect}(\mathscr{Q}^I)$ with 
$t_0 \! \in \! T_0$, $t_j \! \in \! \Lambda$, and $t_0$ can be chosen in $\overline \Lambda^+$
whenever $E$ meets $\overline{\mathscr{Q}}^X$. 

\vspace{2mmplus1mmminus1mm} 

We can also deal with the action of elementary unipotent matrices on the apartment $A_\mathbb{B}$. 
Let $U_{ab}=u_{ab}(F)$ be the image of the  homomorphism $u_{ab}:(F,+)\to{\rm SL}(E)$ defined by $\lambda\mapsto{\rm id}+\lambda E_{ab}$. 
Then the geometric interpretation of the valuation of the additive parameter $\lambda$ is that the fixed-point set of $u=u_{ab}(\lambda)$ in 
$A_\mathbb{B}$ is the half-space  $D_u = \{ c_a - c_b \leq v_F(\lambda) \}$. 

\vspace{2mmplus1mmminus1mm} 

At last, in the case of the special linear group the Cartan decomposition ${\rm SL}(E) = K \cdot \overline T_0^+ \cdot K$, as well as the Iwasawa decomposition ${\rm SL}(E) = K \cdot T_0 \cdot U^-$ (both with respect to $\mathbb{B}$), can be proved by direct matrix computation. 
It can also be checked that the apartment $A_\mathbb{B}$ is a fundamental domain for the action of the Iwahori subgroup $\mathscr{B}$ on the Bruhat-Tits building $X_E$, and that the closure $\overline{\mathscr{Q}}^X$ is a fundamental domain for the action of the maximal compact subgroup $K$ on $X_E$. 

\subsubsection{}
\label{sss - convergence for GL} 
Let us now show that in the case of ${\rm SL}_n(F)$, Chabauty convergence can be proved by elementary matrix computation. 
Let $I \subsetneq \{1;2;...\, n-1\}$ and $\{ v_n \}_{n \geq 1}$ be an $I$-canonical sequence of vertices in $\overline{\mathscr{Q}}$. 
In view of the properties of the $T$-action on the space $\mathscr{A}^I$ (\ref{sss - parametrization for GL}), it is enough to consider a sequence where the vertices all lie in the sector face $\mathscr{Q}^I$ and for which the distances to the vector panels $\Pi^i$ for $i\not\in I$ explode. 
We denote by $\{ K_{v_n} \}_{n \geq 1}$ the associated sequence of maximal compact subgroups. 
Let us illustrate some points (Lemmas \ref{lemma - lower bound cluster value} and \ref{lemma - upper bound cluster value}) of the proof of Theorem \ref{th - CV groups}, which says that $\{ K_{v_n} \}_{n \geq 1}$ converges to $D_I = K_I \ltimes U^I$. 
Using the previous parametrization of $V_\mathscr{Q}$ (\ref{sss - parametrization for GL}), we write 
$v_n=\varpi_F^{\underline \nu(n)}.o$, with $\nu_{i+1}(n) = \nu_i(n)$ for each $n \geq 1$ whenever $i \! \in \! I$ and $\nu_{i+1}(n) - \nu_i(n) \to \infty$ as $n \to \infty$ otherwise. 
The corresponding sequence of compact subgroups is $\{ \varpi_F^{\underline \nu(n)} K \varpi_F^{-\underline \nu(n)} \}_{n \geq 1}$, and we have to show that  it converges to $D_I$. 
It is enough to show that any cluster value of $\{ \varpi_F^{\underline \nu(n)} K \varpi_F^{-\underline \nu(n)} \}_{n \geq 1}$ is equal to $K_I \ltimes U^I$. 

\vspace{2mmplus1mmminus1mm} 

Let $D$ be such a cluster value. 
Then $D$ is the set of limits of converging sequences $\{Êg _n\}_{n \geq 1}$ with 
$g_n \! \in \! \varpi_F^{\underline \nu(n)} K \varpi_F^{-\underline \nu(n)}$ for each $n\geq 1$. 
Let us write $g_n=\varpi_F^{\underline \nu(n)} k_n \varpi_F^{-\underline \nu(n)}$ with 
$k_n \! \in \! {\rm SL}_n(\mathscr{O}_F)$. 
We consider the conjugates $\varpi_F^{\underline \nu(n)} E_{ij} \varpi_F^{-\underline \nu(n)}$ of the elements
of the natural basis $\{ E_{ij} \}_{1 \leq i,j \leq n}$ of the $n \times n$ matrices. 
We have: 
$\varpi_F^{\underline \nu(n)} E_{ij} \varpi_F^{-\underline \nu(n)}=
\varpi_F^{\nu_i(n)-\nu_j(n)} E_{ij}$. 
If $i$ is not equivalent to $j$ for $\sim_I$ and if $i>j$, then this shows that 
$\displaystyle \lim_{n \to \infty} \varpi_F^{\underline \nu(n)} E_{ij} \varpi_F^{-\underline \nu(n)}=0$, which
implies that $D < P_I$. 
If $i$ is equivalent to $j$, then for any $n\geq 1$ the diagonal matrix 
$\varpi_F^{\underline \nu(n)}$ centralizes $E_{ij}$. 
Varying the equivalent $i$ and $j$, we obtain that $\varpi_F^{\underline \nu(n)}$ centralizes $K_I$. 
This implies that we have: $K_I < D$.
Now we consider the case when $i$ is not equivalent to $j$ and $i<j$, and we pick 
$\lambda \! \in \! k$, which we write $\lambda=u\varpi_F^v$ with $u \! \in \! \mathscr{O}_F$ and 
$v \! \in \! {\bf Z}$. 
For each $n \geq 1$, we set $k_n={\rm id}+\varpi_F^{v-(\nu_i-\nu_j)(n)}E_{ij}$. 
There is $M \geq 1$ such that $\nu_i(n)-\nu_j(n) \leq v$ for any $n\geq M$, so $k_n$ belongs to $K$ for
$n \geq M$. 
This shows that any element ${\rm id}+\lambda E_{ij}$ is the limit of an eventually constant sequence 
$\{ g_n \}_{n \geq 1}$ with $g_n \! \in \! \varpi_F^{\underline \nu(n)} K \varpi_F^{-\underline \nu(n)}$. 
Therefore we have: $U^I < D$.

\subsection{Boundary, parametrizations and identifications}Ê
\label{ss - boundary and identifications} 
We describe the limit groups in the boundary of the standard sector $\mathscr{Q}$, and illustrate the geometric parametrization of maximal amenable and distal subgroups. 
Then we announce a study of non maximal Furstenberg compactifications in the ${\rm SL}_n$ case, which is related to A. Werner's compactification in terms of seminorms. 

\subsubsection{}
\label{sss - boundary of standard Weyl chamber} 
A special feature in the case of the group ${\rm SL}_n(F)$ is that the bigger group ${\rm GL}_n(F)$ also acts on the Euclidean building. 
This is clear from the additive norm viewpoint; the drawback of the ${\rm GL}(E)$-action on $X_E$ is that it is not type-preserving, but the big advantadge is that it is vertex-transitive. 
For instance, it is clear that the semigroup $\overline{\Lambda}^+$ acts transitively on the vertices of $V_\mathscr{Q}$. 
In fact it follows from the description of the $\Lambda$-action on the sets $\mathscr{A}^I$ of affine subspaces obtained as suitable intersections of walls (\ref{sss - parametrization for GL}), that for any $I \subsetneq \{1;2;...\, n-1\}$ each limit group of an $I$-canonical sequence is a suitable $\overline{\Lambda}^+$-transform of the group $D_I=K_I\ltimes U^I$. 
In the case where $n=3$, this says that there are three kinds of limit groups in the closure of the sector $\mathscr{Q}$: 
\begin{enumerate}
\item [1.] the single limit group obtained as the semidirect product of the unique maximal compact subgroup ${\rm Diag}(\mathscr{O}_F^\times)$ of the diagonal matrices by the unipotent upper triangular matrices $U$; 
\item [2.] the groups obtained as the semidirect product of a maximal compact subgroup in the upper left ${\rm GL}_2$ diagonal block by the upper triangular unipotent group $U^{\{1\}}$ of the matrices with two additive parameters on the last column (the lower right diagonal coefficient is the inverse of the determinant of the ${\rm GL}_2$ block); 
\item [3.] the groups obtained similarly after replacing the upper left ${\rm GL}_2$ block by the lower right one, and the last column by the first line. 
\end{enumerate}

This leads to the picture below. 
The group $D_{\{1\}}$ on the picture is the limit group of any sequence going to infinity and staying in the sector panel $\mathscr{Q}^{\{1\}}$. 
The second class of $\{1\}$-canonical sequences on the picture is represented by the first vertical dashed line on the right of $\mathscr{Q}^{\{1\}}$. 
The corresponding limit group is: 

\vspace{2mmplus1mmminus1mm} 

\centerline{$
\begin{pmatrix}
A\in{\rm GL}_2\begin{pmatrix}
\mathscr{O}_F^\times&\mathscr{O}_F^\times\\
\varpi_F\mathscr{O}_F^\times&\mathscr{O}_F^\times
\end{pmatrix} &
\begin{matrix}*\in F\\ *\in F\end{matrix}\\
0\qquad0 & {\rm det}(A)^{-1}\in\mathscr{O}_F^\times
\end{pmatrix}$.}

\vspace{2mmplus1mmminus1mm} 

The difference between $D_{\{1\}}$ and the latter group is that the upper left block is not a maximal compact subgroup of the corresponding Levi factor, but an Iwahori subgroup. 
The last class of $\{1\}$-canonical sequences on the picture corresponds to the rightest vertical dashed line. 
The corresponding limit group is obtained by taking 
${\rm GL}_2\begin{pmatrix}
\mathscr{O}_F^\times&(\varpi_F)^{-1}\mathscr{O}_F^\times\\
\varpi_F\mathscr{O}_F^\times&\mathscr{O}_F^\times
\end{pmatrix}$ 
as upper left diagonal block. 
Varying the rays in $\mathscr{Q}$ parallel to $\mathscr{Q}^{\{1\}}$ and taking the limit groups, we describe the vertices in a geodesic ray of the Bruhat-Tits tree of the upper left Levi factor ${\rm GL}_2(F)$ of ${\rm SL}_3(F)$. 

\vspace{1cm}Ê

\begin{picture}(0,0)%
\includegraphics{sector.pstex}%
\end{picture}%
\setlength{\unitlength}{2072sp}%
\begingroup\makeatletter\ifx\SetFigFont\undefined%
\gdef\SetFigFont#1#2#3#4#5{%
  \reset@font\fontsize{#1}{#2pt}%
  \fontfamily{#3}\fontseries{#4}\fontshape{#5}%
  \selectfont}%
\fi\endgroup%
\begin{picture}(10432,8121)(31,-7384)
\put(8146,-6316){\makebox(0,0)[lb]{\smash{{\SetFigFont{11}{13.2}{\familydefault}{\mddefault}{\updefault}$\{2\}$-canonical sequence}}}}
\put(6976,-7261){\makebox(0,0)[lb]{\smash{{\SetFigFont{11}{13.2}{\familydefault}{\mddefault}{\updefault}sector $\mathscr{Q}$}}}}
\put(8281,-961){\makebox(0,0)[lb]{\smash{{\SetFigFont{11}{13.2}{\familydefault}{\mddefault}{\updefault}$D_\varnothing=\begin{pmatrix}a\in \mathscr{O}_F^\times&*\in F&*\in F\\0&b\in\mathscr{O}_F^\times&*\in F\\0&0&(ab)^{-1}\end{pmatrix}$}}}}
\put(2701,434){\makebox(0,0)[lb]{\smash{{\SetFigFont{11}{13.2}{\familydefault}{\mddefault}{\updefault}$D_{\{1\}}=\begin{pmatrix}A\in{\rm GL}_2(\mathscr{O}_F)&\begin{matrix}*\in F\\ *\in F\end{matrix}\\0\qquad 0&{\rm det}(A)^{-1}\in\mathscr{O}_F^\times\end{pmatrix}$}}}}
\put(9676,-3481){\makebox(0,0)[lb]{\smash{{\SetFigFont{11}{13.2}{\familydefault}{\mddefault}{\updefault}$D_{\{2\}}=\begin{pmatrix}{\rm det}(A)^{-1}\in\mathscr{O}_F^\times&*\in F\quad *\in F\\\begin{matrix}0\\0\end{matrix}&A\in{\rm GL}_2(\mathscr{O}_F^\times)\end{pmatrix}$}}}}
\put( 46,-5461){\makebox(0,0)[lb]{\smash{{\SetFigFont{11}{13.2}{\familydefault}{\mddefault}{\updefault}sector panel $\Pi^1$}}}}
\put(1126,-6856){\makebox(0,0)[lb]{\smash{{\SetFigFont{11}{13.2}{\familydefault}{\mddefault}{\updefault}alcove $c$}}}}
\put(4951,-3661){\makebox(0,0)[lb]{\smash{{\SetFigFont{11}{13.2}{\familydefault}{\mddefault}{\updefault}$\{1\}$-canonical sequence}}}}
\end{picture}%

\subsubsection{}
\label{sss - amenable subgroups and stable flags} 
The limit groups described in \ref{sss - boundary of standard Weyl chamber} are the groups $tD_It^{-1}$ for $t\!\in\!\overline{\Lambda}^+$ and $I\subsetneq\{1;2;...\, n-1\}$. 
These groups have the common property to stabilize a flag on the subquotients of which they act via a compact group (their unipotent part acts trivially on it). 
In fact, the ${\rm SL}_n$ case together with an embedding of $F$-algebraic groups, is used to prove the classification of amenable (\ref{ss - amenable}) and distal (\ref{ss - distal}) subgroups of arbitrary semisimple groups over $F$, so strictly speaking, for Theorem \ref{th - parametrization amenable} and Theorem \ref{th - parametrization distal}, the ${\rm SL}_n$ case is a necessary first step more than merely a concrete example. 

\vspace{2mmplus1mmminus1mm}

Let us simply mention that in terms of flags, these theorems say that a distal (resp. a Zariski connected closed amenable) subgroup of 
${\rm SL}_n(F)$ stabilizes a flag on the subquotients of which the corresponding linear (resp. projective) actions are via a compact group \cite{Cornulier}. 
In both cases, the proof of these statements goes by induction on the dimension $n$. 
The main tool for distality is to use Burnside's density theorem 
\cite[\S 4.2 Th\'eor\`eme 1]{Algebre8} 
combined with the non-degeneracy of the trace form. 
These ideas elaborate on the proof of the finiteness of torsion finitely generated linear groups, and already appear in \cite{TitsFree}, \cite{ConGui}, \cite{Prasad} for instance. 
Roughly speaking, the case of amenablility is proved by replacing Burnside's density theorem by Furstenberg's lemma on stabilizers of probability measures on projective spaces \cite[3.2]{Zimmer}. 

\subsubsection{}
\label{sss - other compactifications by seminorms} 
Let us finish this paper by mentioning that in the ${\rm SL}_n$ case, a concrete interpretation of the polyhedral compactification is given, at least for vertices \cite[15]{Landvogt}. 
In order to extend the description of the vertices in $X_E$ as $\mathscr{O}_F$-lattices (\ref{sss - additive norms}), one has to introduce {\it generalized $\mathscr{O}_F$-lattices}, and the notion of homothety has to be extended too \cite[Definition 15.1]{Landvogt}. 
More recently, A. Werner defined two concrete ways to compactify the Bruhat-Tits building of ${\rm SL}_n(F)$; the first procedure uses sublattices in $E$ \cite{WernerSublattices}, and the second one uses seminorms on $E$ \cite{WernerSeminorms}. 
They both lead to compactifications which are different from the polyhedral or the geometric one. 
In a next paper, we will define a family of measure-theoretic (i.e. Furstenberg) compactifications for ${\rm SL}_n(F)$. 
There is one compactification for each choice of a conjugacy class of proper parabolic subgroups. 
The minimal parabolic subgroups lead to the group-theoretic compactification (up to ${\rm SL}_n(F)$-equivariant homeomorphism): it is the maximal measure-theoretic compactification. 
We will also investigate the connection between A. Werner's compactifications and intermediate measure-theoretic compactifications. 

\bibliographystyle{amsalpha}
\bibliography{Groupes}

\vspace{1cm} 

Institut de Recherches en Math\'ematiques Avanc\'ees 
\hfill
Institut Camille Jordan 

UMR 6625 CNRS - Rennes 1
\hfill
UMR 5208 CNRS - Lyon 1

Universit\'e de Rennes 1 \hfill
Universit\'e de Lyon 1

Campus de Beaulieu 
\hfill
21 rue Claude Bernard

35042 Rennes Cedex, France
\hfill
69622 Villeurbanne Cedex, France

{\tt yves.guivarch@math.univ-rennes1.fr} \hfill
{\tt remy@igd.univ-lyon1.fr}

\addtolength{\parindent}{-1.6pt}

\end{document}